\documentclass[twocolumn, envcountsect, envcountsame]{svjour3}  
\smartqed  
\usepackage{graphicx}
\usepackage{mathptmx}      
\usepackage[T1]{fontenc}

\usepackage{amsmath}
\usepackage{amssymb}
\usepackage{amsfonts}
\usepackage{mathdots}
\usepackage{oldgerm}
\usepackage{mathrsfs}
\usepackage{enumerate}
\usepackage{siunitx}
\usepackage{comment}
\usepackage{supertabular}

\usepackage{subfigure}

\usepackage{xcolor}
\usepackage{colortbl}
\definecolor{hage}{rgb}{0.4,0.6,1}

\usepackage{tikz}
\usetikzlibrary{fit}
\usetikzlibrary{calc}
\pgfdeclarelayer{background}
\pgfsetlayers{background,main}
\colorlet{inbox}{lightgray!20}
\colorlet{outbox}{lightgray!50}

\makeatletter
\renewcommand*\env@matrix[1][*\c@MaxMatrixCols c]{%
  \hskip -\arraycolsep
  \let\@ifnextchar\new@ifnextchar
  \array{#1}}

\newcommand{\cC}{\mathcal{C}}

\newcommand{\cW}{\mathcal{W}}
\newcommand{\cG}{\mathcal{G}}

\newcommand{\cK}{\mathcal{K}}

\newcommand{\cD}{\mathcal{D}}

\newcommand{\cF}{\mathcal{F}}

\newcommand{\RN}[1]{\uppercase\expandafter{\romannumeral#1}}

\newcommand{\N}{{\mathbb{N}}}

\newcommand{\R}{{\mathbb{R}}}
\newcommand{\C}{{\mathbb{C}}}

\newcommand{\Gl}{\mathbf{Gl}}

\newcommand{\vt}{\vartheta}
\newcommand{\vp}{\varphi}

\DeclareMathOperator{\rk}{rk}

\DeclareMathOperator{\im}{im}

\newcommand{\setdef}[2]{\left\{\ #1\ \left|\ \vphantom{#1} #2\ \right.\right\}}
\newcommand{\ddt}{\tfrac{\text{\normalfont d}}{\text{\normalfont d}t}}

\newcommand{\diag}{\text{\rm diag\,}}
\newcommand{\ds}[1]{{\rm \, d} #1 \,}

\DeclareMathOperator{\spn}{span}

\DeclareMathOperator{\loc}{loc}

\newenvironment{smallpmatrix}
{\left(\begin{smallmatrix}}
{\end{smallmatrix}\right)}

\newenvironment{smallbmatrix}
{\left[\begin{smallmatrix}}
{\end{smallmatrix}\right]}

\begin{document}

\title{Tracking control for underactuated non-minimum phase multibody systems\thanks{This work was supported by the German Research Foundation (Deutsche Forschungsgemeinschaft) via the grants BE 6263/1-1 and SE 1685/6-1.}}



\author{Thomas Berger \and Svenja Dr\"ucker \and Lukas Lanza \and Timo Reis \and Robert Seifried}


\institute{Corresponding author: Thomas Berger  \at Tel.: +49 5251 60-3779
            \and
            Thomas Berger, Lukas Lanza \at
             Institut f\"ur Mathematik, Universit\"at Paderborn, Warburger Str.~100, 33098~Paderborn, Germany \\
              \email{\{thomas.berger, lanza\}@math.upb.de}
           \and
           Svenja Dr\"ucker, Robert Seifried \at
           Institute of Mechanics and Ocean Engineering, Hamburg University of Technology, Ei\ss endorfer Stra\ss e 42, 21073 Hamburg, Germany\\
           \email{\{svenja.druecker, robert.seifried\}@tuhh.de}
           \and
            Timo Reis \at
              Department of Mathematics, Universit\"at Hamburg, Bundesstra\ss e 55, 20146 Hamburg, Germany \\
              \email{timo.reis@math.uni-hamburg.de}           
}

\date{Received: date / Accepted: date}

\maketitle

\begin{abstract}
We consider tracking control for multibody systems which are modelled using holonomic and nonholonomic constraints. Furthermore, the systems may be underactuated and contain kinematic loops and are thus described by a set of differential-algebraic equations that cannot be reformulated as ordinary differential equations in general. We propose a control strategy which combines a feedforward controller based on the servo-constraints approach with a feedback controller based on a recent funnel control design. As an important tool for both approaches we present a new procedure to derive the internal dynamics of a multibody system. Furthermore, we present a feasible set of coordinates for the internal dynamics avoiding the effort involved with the computation of the Byrnes-Isidori form. The control design is demonstrated by a simulation for a nonlinear non-minimum phase multi-input, multi-output robotic manipulator with kinematic loop.

\end{abstract}

\section{Introduction}\label{Sec:Intro}

In the present paper, we propose a combined feedforward and feedback tracking control strategy for underactuated non-minimum phase multibody systems. We follow a popular approach to two degree of freedom controller design as proposed e.g.\ in~\cite{SkogPost05}. The feedforward control input design is based on a reference model of the system such that, if the model truthfully captures reality, exact tracking of a given reference signal by the output is achieved. We utilize the so-called servo-constraints approach for feedforward control design.

In order to compensate (inevitable) modeling errors, uncertainties, disturbances, noise, etc.\ an additional feedback loop is used to stabilize the system around the given reference signal. Since a robust feedback controller is desired and the output must respect prescribed error margins around the reference signal we use the funnel controller first proposed in~\cite{IlchRyan02b}. Since the funnel controller presented in~\cite{IlchRyan02b} is not feasible for non-minimum phase systems we use an extension recently developed in~\cite{Berg20a}.

An important tool both in feedforward and feedback control design is the Byrnes-Isidori form, which allows a decoupling of the internal dynamics of the system. However, a calculation of the Byrnes-Isidori form and the accompanying nonlinear transformation requires a lot of computational effort in general. The approach presented in the present paper avoids this computation. In the feedforward control design, the servo-constraints constitute an approach which does not require the Byrnes-Isidori form for the solution of the inverse model. For the feedback control design, we present a new approach to choose a set of variables for the internal dynamics of the system directly in terms of the system parameters -- circumventing the Byrnes-Isidori form. The feedback controller is then based on this representation of the internal dynamics.

More details of the considered system class and the proposed control methodology are given in the following. Furthermore, we recall the concept of vector relative degree and the Byrnes-Isidori form. {We like to emphasize that the proposed control design is potentially feasible in any open set where the vector relative degree is well-defined.}

{Before continuing, we like to summarize the main contributions of the present paper in the following.
\begin{enumerate}[a)]
  \item A new method for computing the internal dynamics of a multibody system which avoids the DAE formulation is presented (Section~\ref{Sec:Computing-ID}).
  \item To circumvent the computational effort involved with the Byrnes-Isidori form of the auxiliary ODE arising in~a), a feasible set of coordinates which directly yields the internal dynamics is derived (Section~\ref{Sec:intdyn-coord}).
  \item A feedforward control strategy based on servo-constraints for model inversion of non-minimum phase systems is used. For the required solution of the arising boundary value problem we apply newly proposed boundary conditions, which significantly simplify the inversion process (Section~\ref{Sec:Servo}).
  \item We present a new funnel control design for nonlinear non-minimum phase systems which only have a vector relative degree (Section~\ref{Sec:FunCon}).
  \item We demonstrate that the combination of the feedforward and feedback control strategies is able to achieve tracking with prescribed performance for a nonlinear, non-minimum phase robotic manipulator with kinematic loop (described by a DAE that cannot be reformulated as an ODE) -- the controller performance of the combination is favorable compared to the individual controllers (Section~\ref{Sec:Robot}).
\end{enumerate}}

\subsection{Nomenclature}\label{Ssec:Nomencl}

\noindent{
\hspace*{-3mm}
\begin{supertabular}{p{60pt}p{160pt}}
{$\N$, $\R$, $\C$} & {the set of natural, real and complex numbers, resp.} \\
$\R_{\ge 0}$ &$=[0,\infty)$, \\
{$\C_+$ \big($\C_-$\big)} & {the open set of complex numbers with positive (negative) real part} \\
{$\R^{n\times m}$}
  &   {the set of real $n\times m$ matrices}\\
$ \Gl_n$ & the group of invertible matrices in $\R^{n\times n}$\\
{$\sigma(A)$} & {the spectrum of a matrix $A\in\R^{n\times n}$} \\
$\mathcal{L}^\infty_{\loc}(I\!\to\!\R^n)$
  &  the set of locally essentially bounded functions $f:I\!\to\!\R^n$, $I\subseteq\R$ an interval\\
$\mathcal{L}^\infty(I\!\to\!\R^n)$
  &  the set of essentially bounded functions $f:I\!\to\!\R^n$\\
$\|f\|_\infty$ & = ${\rm ess\ sup}_{t\in I} \|f(t)\|$\\
$\mathcal{W}^{k,\infty}(I\!\to\!\R^n)$
  &  the set of $k$-times weakly differentiable functions $f:I\!\to\!\R^n$ such that $f,\ldots, f^{(k)}\in \mathcal{L}^\infty(I\!\to\!\R^n)$\\
$\mathcal{C}^k(V\!\to\!\R^n) $
  &  the set of $k$-times continuously differentiable functions $f:V\!\to\! \R^n$, $V\subseteq\R^m$\\
$\mathcal{C}(V\!\to\!\R^n) $
  &$= \mathcal{C}^0(V\!\to\!\R^n) $\\
$\left.f\right|_{W}$ & restriction of the function $f:V\!\to\!\R^n$ to $W\subseteq V$\\
\end{supertabular}

\subsection{System class}\label{Ssec:SysClass}

We consider multibody systems which possibly contain kinematic loops as well as holonomic and nonholonomic constraints. The equations of motion are given by
\begin{equation}\label{eq:MBS}
\begin{aligned}
    \dot q(t) &= v(t),\\
    M(q(t)) \dot v(t) &= f\big(q(t),v(t)\big) + J(q(t))^\top \mu(t)\\
    &\quad  + G(q(t))^\top \lambda(t) + B(q(t))\, u(t),\\
    0 &= J(q(t)) v(t) + j(q(t)),\\
    0 &= g(q(t)),\\
    y(t) &=h\big(q(t),v(t)\big),
\end{aligned}
\end{equation}
with
\begin{itemize}
  \item the generalized coordinates~$q: I\to \R^n$ and generalized velocities~$v: I\to \R^n$ (in the case of no constraints), where $I\subseteq \R_{\ge 0}$ is some interval,
  \item the generalized mass matrix $M:\R^n\to\R^{n\times n}$,
  \item the generalized forces $f:\R^n\times \R^n\to\R^n$,
  \item the holonomic constraints $g:\R^n\to\R^\ell$ and $G:\R^n\to\R^{\ell\times n}$,
  \item the nonholonomic constraints ${J:\R^n\to\R^{p\times n}}$ and ${j:\R^n\to\R^p}$ (which may incorporate a change of physical units),
  \item the input distribution matrix $B:\R^n\to\R^{n\times m}$,
  \item the output measurement function $h:\R^n\times\R^n\to\R^m$.
\end{itemize}
The functions $u:\R\to\R^m$ are the inputs that influence the multibody system~\eqref{eq:MBS} in affine form.
We explicitly allow for ${m < n - \ell}$, which means that, if no nonholonomic constraints are present, underactuated multibody systems are encompassed by this formulation.
The affine input has the interpretation of a force or torque, however standard industrial actuators are typically velocity controlled and not force controlled, which may lead to different system properties~\cite[Sec.~4.4]{OttoSeif18a}. In the present paper, we assume that the actuators are force/torque controlled and that the multibody system is given in the form~\eqref{eq:MBS}.

The functions $y:\R\to\R^m$ are the outputs associated with the multibody system~\eqref{eq:MBS} and typically represent appropriate measurements. The generalized forces~$f$ usually encompass several terms, including applied forces, Coriolis, centrifugal and gyroscopic forces. The functions $\lambda:\R\to\R^\ell$ and $\mu:\R\to\R^p$ are the Lagrange multipliers corresponding to the holonomic and nonholonomic constraints, resp. The functions in~\eqref{eq:MBS} are assumed to be sufficiently smooth and to satisfy the natural assumptions $\ell + p + m \le n$ and
\begin{equation}\label{eq:assumptions}
  \forall\,q\in\R^n:\quad M(q) = M(q)^\top>0\quad \text{and}\quad g'(q) = G(q).
\end{equation}
Let us stress that invertibility of~$M(q)$ is not sufficient, we require that it is positive definite. This follows for example from deriving the equations of motion using the direct application of d'Alemberts principle in the form of Lagrange.

System~\eqref{eq:MBS} is a differential-algebraic system and its treatment needs special care. For an overview of important concepts for differential-algebraic control systems we refer to~\cite{Berg14a,KunkMehr06,Sime13} and the series of survey articles~\cite{IlchReis13b}.

\subsection{Vector relative degree}\label{Ssec:RelDeg}

An important property of the system~\eqref{eq:MBS} is its  vector relative degree, which, roughly speaking, is the collection of numbers of derivatives of each output component needed so that the input appears explicitly. We briefly recall the required concepts from~\cite{Isid95} and, to this end, consider a general nonlinear system affine in the control given by
\begin{equation}\label{eq:System}
\begin{aligned}
 \dot{x}(t) &= f(x(t)) + g(x(t)) u(t), \quad x(0) = x^0 \in \R^d, \\
y(t) &= h(x(t)),
\end{aligned}
\end{equation}
where ${f: \R^d \to \R^d}$,  ${g: \R^d \to \R^{d \times m}}$ and ${h: \R^d \to \R^m}$ are sufficiently smooth and $x:I\to \R^n$,  where $I\subseteq \R_{\ge 0}$ is some interval. Denote the components of~$h$ by $h_1,\ldots,h_m$, and recall the definition of the Lie derivative of~$h_i$, $i\in\{1,\ldots,m\}$, along a vector field~$f$ at a point~$z \in U \subseteq \R^d$,~$U$ open:
\begin{align*}
\left(L_f h_i \right)(z) := h_i'(z) f(z),
\end{align*}
where~$h_i'$ is the Jacobian of~$h_i$, i.e., the transpose of the gradient of $h_i$.
We may gradually define~$L_f^k h_i = L_f(L_f ^{k-1} h_i)$ with $L_f^0 h_i = h_i$. Furthermore, denoting with $g_j(z)$ the columns of $g(z)$ for $j=1,\ldots,m$, we define
\[
    \left(L_g h_i \right)(z) := [(L_{g_1} h_i)(z), \ldots,(L_{g_m} h_i)(z)] \in\R^{1\times m}.
\]
Now, in virtue of~\cite{Isid95}, the system~\eqref{eq:System} is said to have \textit{vector relative degree}~$(r_1,\ldots,r_m)\in \mathbb{N}^{1\times m}$ on~$U$, if for all~$z \in U$ we have:
\begin{align*}
&\forall\,i\in\{1,\ldots,m\}\ \forall\, k \in \{ 0,...,r_i-2 \}: \ (L_g L_f^k h_i)(z) = 0 \\
&\text{and}\quad  \Gamma(z) = \begin{bmatrix} \Gamma_1(z)\\ \vdots\\ \Gamma_m(z)\end{bmatrix} :=  \begin{bmatrix} (L_g L_f^{r_1-1} h_1)(z)\\ \vdots\\ (L_g L_f^{r_m-1} h_m)(z)\end{bmatrix} \in \Gl_m,
\end{align*}
where $\Gamma:U\to\Gl_m$ denotes the high-frequency gain matrix. If $r_1 = \ldots = r_m =: r\in\N$, then system~\eqref{eq:System} is said to have \textit{strict relative degree}~$r$.

If system~\eqref{eq:System} has vector relative degree ${(r_1,\ldots,r_m) \in \mathbb{N}^{1\times m}}$ on an open set $U \subseteq \R^d$, then there exists a (local) diffeomorphism~${\Phi: U \to W \subseteq \R^d}$, $W$ open, such that the coordinate transformation $\begin{smallpmatrix} \xi(t) \\ \eta(t) \end{smallpmatrix} = \Phi(x(t))$, $\xi(t) =  (\xi_{ij}(t))_{i=1,\ldots,m; j=1,\ldots,r_i}\in\R^{\bar r}$, $\eta(t)\in\R^{d-\bar r}$, where $\bar r= \sum_{i=1}^m r_i$, puts system~\eqref{eq:System} into Byrnes-Isidori form
\begin{align}
\label{eq:BIF}
y_i(t) &= \xi_{i1}(t), \nonumber \\
\dot{\xi}_{i1}(t) &= \xi_{i2}(t),  \nonumber \\
& \ \ \vdots \\
\dot{\xi}_{i,r_i-1}(t) &= \xi_{ir_i}(t), \nonumber \\
\dot{\xi}_{ir_i}(t) &= (L_f^{r_i} h_i)\big(\Phi^{-1}(\xi(t),\eta(t))\big)\nonumber\\
&\quad + \Gamma_i\big(\Phi^{-1}(\xi(t),\eta(t))\big) u(t),\qquad\quad i=1,\ldots, m, \nonumber\\
\dot{\eta}(t) &= q(\xi(t),\eta(t)) + p(\xi(t),\eta(t))u(t). \nonumber
\end{align}
The last equation in~\eqref{eq:BIF} represents the internal dynamics of system~\eqref{eq:System}. The diffeomorphism~$\Phi$ can be represented as
\begin{align}
\label{eq:BIF-trafo}
\Phi(x) =
\begin{pmatrix}
h_1(x) \\
(L_fh_1)(x) \\
\vdots \\
(L_f^{r_1-1}h_1)(x) \\
 h_2(x) \\
 \vdots\\
(L_f^{r_m-1}h_m)(x)\\
\phi_{\bar r+1}(x) \\
\vdots \\
\phi_{d}(x)
\end{pmatrix},
\end{align}
where $\phi_i:U\to\R$, $i=\bar r +1,\ldots,d$ are such that $\Phi'(z)\in\Gl_d$ for all $z\in U$. By a~distribution, we mean a~mapping from an open set $U\subseteq\R^d$ to the set of all subspaces of $\R^d$.
If the distribution~$x\mapsto\cG(x) := \im (g(x))$ in~\eqref{eq:System} is involutive, then the functions~$\phi_i$ in~\eqref{eq:BIF-trafo} can additionally be chosen such that
\begin{align}
\label{eq:involutivity}
\forall\, i=\bar r +1,\ldots,d \ \ \forall\, z\in U:\ (L_g \phi_i) (z) = 0,
\end{align}
by which $p(\cdot) = 0$ in~\eqref{eq:BIF}, cf.~\cite[Prop.~5.1.2]{Isid95}. Recall from~\cite[Sec.~1.3]{Isid95} that~
a~distribution $\cG$ is involutive, if for all smooth vector fields $g_1, g_2: U \to \R^d$ with $g_i(x)\in \cG(x)$ for all $x \in U$ and $i=1,2$ we have that the Lie bracket $[g_1,g_2](x) = g_1'(x) g_2(x) - g_2'(x) g_1(x)$ satisfies ${[g_1,g_2](x) \in \cG(x)}$ for all $x\in U$.

\subsection{Control methodology}\label{Ssec:ContrMeth}

As in our preliminary work~\cite{BergOtto19} we follow the popular controller design methodology for mechanical systems undergoing large motions, where a feedforward controller is combined with a feedback controller, see e.g.~\cite{SkogPost05}. Both controllers are designed independently for tracking of a reference trajectory $y_{\rm ref}:\R_{\ge 0}\to\R^m$. The feedforward control input~$u_{\rm ff}$ is designed using a reference (inverse) model of the system, while the feedback control input~$u_{\rm fb}$ is applied to the actual system. The latter may deviate from the reference model; the situation is depicted in Figure~\ref{fig:blockdiagram}.


\begin{figure}[h!t]
  \centering
  \resizebox{0.47\textwidth}{!}{
\begin{tikzpicture}[thick,node distance = 12ex, box/.style={fill=white,rectangle, draw=black}, blackdot/.style={inner sep = 0, minimum size=3pt,shape=circle,fill,draw=black},plus/.style={fill=white,circle,inner sep = 0, minimum size=5pt,thick,draw},metabox/.style={inner sep = 3ex,rectangle,draw,dotted,fill=gray!20!white}]

    \node (sys)     [box,minimum size=8ex] {$\begin{array}{c} \text{Multibody}\\ \text{System}\end{array}$};
    \node (ufork)   [plus, left of = sys, xshift=-1ex] {$+$};
    \node (invsys)  [box,minimum size=8ex,left of = ufork,xshift=-1ex] {$\begin{array}{c} \text{Inverse}\\ \text{Model}\end{array}$};
    \node (yref)    [blackdot, left of = invsys, xshift=4ex] {};
    \node (yref1)   [minimum size=0pt, inner sep = 0pt, left of = yref, xshift=8ex] {};
    \node (fun)     [box, below of = ufork, yshift=2ex, xshift=0ex, minimum size=8ex]  {$\begin{array}{c} \text{Feedback}\\ \text{Controller}\end{array}$};
    \node (right1)  [blackdot,right of = sys,xshift=-2ex] {};
    \node (end)     [minimum size=0pt, inner sep = 0pt, right of = right1, xshift=-8ex] {};

    \draw[->] (ufork)   -- (sys)  {};  
    \draw[->]  (sys)     -- (end)  node[pos=0.4,above] {$y(t)$};
    \draw[->] (invsys)   -- (ufork)    node[midway,above] {$u_{\rm ff}(t)$};
    \draw[->] (fun)   -- (ufork)    node[pos=0.5,left] {$u_{\rm fb}(t)$};
    \draw[->] (yref1) -- (invsys)    node[pos=0.5,above] {$y_{\rm ref}(t)$};
    \draw[->]  (yref)     |- (fun)  {};
    \draw[->]  (right1)     |- (fun)  {};

  \end{tikzpicture}
  }
\caption{Two degree of freedom control approach for multibody systems.}
\label{fig:blockdiagram}
\end{figure}
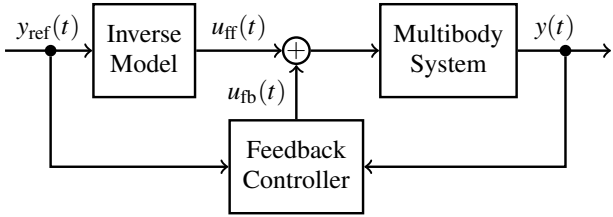

As in~\cite{BergOtto19} the feedforward control input $u_{\rm ff}$ is computed using the method of servo-constraints introduced in~\cite{Blaj92}, cf.\ also~\cite{BlajKolo04,BlajKolo11}. In this framework the equations of motion of the multibody system~\eqref{eq:MBS} are extended by constraints which enforce the output to coincide with the reference trajectory~$y_{\rm ref}$. The resulting set of differential-algebraic equations (DAEs) is  solved numerically for the inverse system, by which the feedforward control input $u_{\rm ff}$ is obtained. The details of the method are presented in Section~\ref{Sec:Servo}.

Servo-constraints have been successfully used for the control of rotary cranes~\cite{BetsQuas09}, overhead gantry cranes~\cite{BlajKolo04,OttoSeif18}, (infinitely long) mass-spring chains~\cite{AltmHeil17,FumaMasa10} and flexible robots~\cite{BurkSeifEber19}. The addition of the servo-constraints to the system may result in a higher differentiation index, as shown in~\cite{Camp95b}; see~\cite{BrenCamp89,Gear88} and also~\cite{KunkMehr06} for a definition of the differentiation index. A higher index causes difficulties in the numerical solution and hence index reduction methods are frequently used, cf.~\cite{HairWann96}. Popular approaches are index reduction by projection onto the constrained and unconstrained directions~\cite{BlajKolo04} and index reduction by minimal extension~\cite{AltmBets16,BetsAltm16}.

The feedback control input~$u_{\rm fb}$ is generated by a dynamic state feedback of the form
\begin{equation}\label{eq:objcontr}
\begin{aligned}
\dot{z}(t)\,&=F\big(t,z(t),q(t),v(t),y_{\rm ref}(t)\big),\\
u_{\rm fb}(t)\,&=G\big(t,z(t),q(t),v(t),y_{\rm ref}(t)\big).
\end{aligned}
\end{equation}
We choose a feedback control design based on the funnel control methodology~\cite{BergLe18a,IlchRyan02b}. Recently, this method has been extended to linear non-minimum phase systems in~\cite{Berg20a}. Based on a linearization of the internal dynamics, the approach from~\cite{Berg20a} has been used for tracking control of a nonlinear non-minimum phase robotic manipulator in~\cite{BergLanz20}. {We stress that this approach to non-minimum phase systems, compared to classical funnel control, requires certain knowledge of system parameters and measurements of state variables.}

{As alternative suitable feedback control methods -- not pursued in this work -- we like to mention prescribed performance control~\cite{BechRovi14,DaiHe18} and control barrier functions~\cite{Jin19,TeeRen11}.}

The objective of funnel control is to design a feedback control law such that in the closed-loop system the tracking error~${e(t) = y(t) - y_{\rm ref}(t)}$ evolves within the boundaries of a prescribed performance funnel
\begin{equation}
\mathcal{F}_{\varphi} := \setdef{(t,e)\in\R_{\ge 0} \times\R^m}{\varphi(t) \|e\| < 1},\label{eq:perf_funnel}
\end{equation}
which is determined by a~function~$\varphi$ belonging to
\begin{equation}
\Phi_k \!:=\!
\left\{
\varphi\in  \cC^k(\R_{\ge 0}\to\R)
\left|\!\!\!
\begin{array}{l}
\text{ $\varphi, \dot \varphi,\ldots,\varphi^{(k)}$ are bounded,}\\
\text{ $\varphi (\tau)>0$ for all $\tau>0$,}\\
 \text{ and }  \liminf_{\tau\rightarrow \infty} \varphi(\tau) > 0
\end{array}
\right.\!\!\!
\right\},
\label{eq:Phir}\end{equation}
where $k = \max_{i=1,\ldots, m} r_i$ and $(r_1,\ldots,r_m)$ is the vector relative degree of the considered system. Furthermore, all involved signals should remain bounded.

The boundary of $\cF_\varphi$ is given by the reciprocal of $\varphi$, see Figure~\ref{Fig:funnel}. The case $\varphi(0)=0$ is explicitly allowed, which means that no restriction is put on the initial value and the funnel boundary $1/\varphi$ has a pole at $t=0$ in this case.

\begin{figure}[h!tp]
\begin{center}
\begin{tikzpicture}[scale=0.45]
\tikzset{>=latex}
  \filldraw[color=gray!25] plot[smooth] coordinates {(0.15,4.7)(0.7,2.9)(4,0.4)(6,1.5)(9.5,0.4)(10,0.333)(10.01,0.331)(10.041,0.3) (10.041,-0.3)(10.01,-0.331)(10,-0.333)(9.5,-0.4)(6,-1.5)(4,-0.4)(0.7,-2.9)(0.15,-4.7)};
  \draw[thick] plot[smooth] coordinates {(0.15,4.7)(0.7,2.9)(4,0.4)(6,1.5)(9.5,0.4)(10,0.333)(10.01,0.331)(10.041,0.3)};
  \draw[thick] plot[smooth] coordinates {(10.041,-0.3)(10.01,-0.331)(10,-0.333)(9.5,-0.4)(6,-1.5)(4,-0.4)(0.7,-2.9)(0.15,-4.7)};
  \draw[thick,fill=lightgray] (0,0) ellipse (0.4 and 5);
  \draw[thick] (0,0) ellipse (0.1 and 0.333);
  \draw[thick,fill=gray!25] (10.041,0) ellipse (0.1 and 0.333);
  \draw[thick] plot[smooth] coordinates {(0,2)(2,1.1)(4,-0.1)(6,-0.7)(9,0.25)(10,0.15)};
  \draw[thick,->] (-2,0)--(12,0) node[right,above]{\normalsize$t$};
  \draw[thick,dashed](0,0.333)--(10,0.333);
  \draw[thick,dashed](0,-0.333)--(10,-0.333);
  \node [black] at (0,2) {\textbullet};
  \draw[->,thick](4,-3)node[right]{\normalsize$\lambda$}--(2.5,-0.4);
  \draw[->,thick](3,3)node[right]{\normalsize$(0,e(0))$}--(0.07,2.07);
  \draw[->,thick](9,3)node[right]{\normalsize$\varphi(t)^{-1}$}--(7,1.4);
\end{tikzpicture}
\end{center}
\caption{Error evolution in a funnel $\cF_{\varphi}$ with boundary $\varphi(t)^{-1}$ for $t > 0$.}
\label{Fig:funnel}
\end{figure}
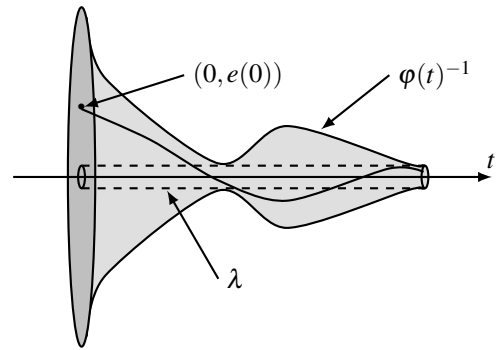

Each performance funnel given by~$\mathcal{F}_{\varphi}$ is bounded away from zero, since boundedness of $\varphi$ implies existence of $\lambda>0$ such that $1/\varphi(t)\geq\lambda$ for all $t > 0$. The funnel boundary is not necessarily monotonically decreasing and widening the funnel over some later time interval might be beneficial, e.g.\ in the presence of periodic disturbances. 


The detailed design of a feedback controller of the form~\eqref{eq:objcontr} yielding the feedback control input~$u_{\rm fb}$ is presented in Section~\ref{Sec:FunCon} and {extends} the recent approach~\cite{Berg20a,BergLanz20}. The funnel controller was already successfully applied e.g.\ in temperature control of chemical reactor models~\cite{IlchTren04}, DC-link power flow control~\cite{SenfPaug14}, voltage and current control of electrical circuits~\cite{BergReis14a}, control of peak inspiratory pressure~\cite{PompWeye15}, adaptive  cruise  control~\cite{BergRaue18,BergRaue20} and control of industrial servo-systems~\cite{Hack17} and underactuated multibody systems~\cite{BergOtto19}.

\subsection{Organization of the present paper}\label{Ssec:Org}
The present work is organized as follows. In Section~\ref{Sec:Computing-ID} we present a novel procedure to compute the internal dynamics of a multibody system. To this end, we introduce an auxiliary input and output to avoid the DAE formulation.
In Section~\ref{Sec:intdyn-coord} we present a feasible set of coordinates of the internal dynamics.
Section~\ref{Sec:Servo} provides a presentation of the open-loop control strategy based on servo-constraints and in particular details the case of unstable internal dynamics.
In Section~\ref{Sec:FunCon} we present a funnel-based feedback-controller design which is based on recently developed results for linear non-minimum phase systems.
Finally, in Section~\ref{Sec:Robot} we apply the findings of the present paper to achieve output tracking of a nonlinear, non-minimum phase, multi-input, multi-output multibody system with kinematic loop, namely a robotic manipulator which is described by a DAE that cannot be conveniently reformulated as an ordinary differential equation (ODE).

\section{Computing the internal dynamics} \label{Sec:Computing-ID}

In this section we present a novel approach to decouple the internal dynamics of the multibody system~\eqref{eq:MBS} avoiding the DAE formulation. In principle, it would be possible to obtain a Byrnes-Isidori form for a DAE directly, see~\cite{Berg17b}. Here we define auxiliary inputs and outputs by removing the constraints and adding the Lagrange multipliers to the input functions. For this auxiliary ODE system we derive the Byrnes-Isidori form as in~\eqref{eq:BIF} and then add the constraints which have been removed before, thus obtaining the internal dynamics.

\subsection{The general case}

Consider~\eqref{eq:MBS} and define an auxiliary input and output  as follows:
\begin{align*}
     u_{\rm aux}(t) &:= \begin{pmatrix} \mu(t)\\ \lambda(t)\\ u(t)\end{pmatrix},\\
     y_{\rm aux}(t) &:= O(q(t),v(t)) = \begin{pmatrix} O_1(q(t),v(t))\\ O_2(q(t),v(t))\\ O_3(q(t),v(t))\end{pmatrix}\\
      &= \begin{pmatrix} J(q(t)) v(t) + j(q(t))\\ g(q(t))\\ h(q(t),v(t))\end{pmatrix}.
\end{align*}
With the state variables $x_1 = q$, $x_2 = v$ and $x= (x_1^\top, x_2^\top)^\top$ we may now consider the auxiliary ODE system (omitting the argument~$t$ for brevity)
\begin{equation}\label{eq:auxODE}
\begin{aligned}
  \dot x  &= \underset{=:F(x)}{\underbrace{\begin{pmatrix} x_2\\ M(x_1)^{-1} f(x_1,x_2) \end{pmatrix}}} \\
  &\quad+ \underset{=:K(x)}{\underbrace{\diag\left(I_n, M(x_1)^{-1}\right) \begin{bmatrix} 0&0&0\\ J(x_1)^\top & G(x_1)^\top & B(x_1)\end{bmatrix}}} u_{\rm aux},\\
  y_{\rm aux} &= O(x),
\end{aligned}
\end{equation}
which is of the form~\eqref{eq:System}, and decouple its internal dynamics using the Byrnes-Isidori form. After that we add the constraints $y_{\rm aux,1} = O_1(x) = 0$ and $y_{\rm aux,2} = O_2(x) = 0$ to derive the internal dynamics of~\eqref{eq:MBS}.

First we calculate the vector relative degree of~\eqref{eq:auxODE}. Let $U\subseteq\R^{2n}$ be open such that the following Lie-derivatives exist for $x = (x_1^\top, x_2^\top)^\top \in U$:
\begin{align*}
  &(L_{K} O_1)(x) = O_1'(x) K(x)\\
  &= \begin{bmatrix} \frac{\partial}{\partial x_1} \big(J(x_1) x_2 + j(x_1)\big) & \frac{\partial}{\partial x_2} \big(J(x_1) x_2 + j(x_1)\big) \end{bmatrix} K(x)\\
  &= J(x_1)M(x_1)^{-1}\begin{bmatrix} J(x_1)^\top & G(x_1)^\top & B(x_1)\end{bmatrix}
\end{align*}
and
\[
  (L_{K} O_2)(x) = O_2'(x) K(x) = [g'(x_1),0] K(x) = 0,
\]
as well as
\begin{align*}
  &(L_{K} L_{F}  O_2)(x) = \big(O_2'(x) F(x)\big)' K(x)\\
  &= \left([ G(x_1) , 0] \begin{pmatrix} x_2\\ M(x_1)^{-1} f(x_1,x_2) \end{pmatrix}\right)' K(x)\\
  &= \begin{bmatrix} \frac{\partial}{\partial x_1} \big(G(x_1) x_2\big) & G(x_1)\end{bmatrix} K(x)\\
  &= G(x_1)M(x_1)^{-1}\begin{bmatrix} J(x_1)^\top & G(x_1)^\top & B(x_1)\end{bmatrix}.
\end{align*}
Since both $J(x_1)M(x_1)^{-1}J(x_1)^\top$ and $G(x_1)M(x_1)^{-1}G(x_1)^\top$ are invertible it follows that a vector relative degree of system~\eqref{eq:auxODE}, if it exists, is of the form $r = (r_1, r_2, r_3)$ with $r_1 = (1,\ldots,1) \in\N^{1\times p}$, $r_2 = (2,\ldots,2)\in\N^{1\times\ell}$ and some $r_3\in\N^{1\times m}$.

In the following, we restrict ourselves to the case where the multi-index~$r_3$ has equal entries, i.e.,
\[
    r_3 = (\hat r, \ldots, \hat r)\in \N^{1\times m}\quad \text{for some $\hat r\in \N$.}
\]
Then it is necessary that, for all $x\in U$,
\[
     \forall\, i=0,\ldots, \hat r-2:\ (L_{K} L_{F}^i  O_3)(x) = 0
\]
and
\[
    \Gamma(x) = \begin{bmatrix}  \Gamma_1(x)\\ \Gamma_2(x)\\ \Gamma_3(x)\end{bmatrix} = \begin{bmatrix}  (L_{K} O_1)(x)\\ (L_{K} L_{F}  O_2)(x)\\ (L_{K} L_{F}^{\hat r-1}  O_3)(x)\end{bmatrix} \in \Gl_{\ell + p + m}.
\]
The first two block rows are given above, the last row strongly depends on the shape of the function~$h$. In the special case that inputs and outputs are colocated we may explicitly calculate~$\Gamma_3(x)$ and show that the relative degree is well defined, see Section~\ref{Sec:colocate}.

Now we can compute the Byrnes-Isidori form of~\eqref{eq:auxODE} using the state space transformation as in~\eqref{eq:BIF-trafo}, which in this case reads
\begin{equation} \label{eq:Phi-aux}
    \Phi(x) = \begin{pmatrix} O_1(x)\\ O_2(x)\\ (L_{F}  O_2)(x) \\ O_3(x) \\ (L_{F}  O_3)(x)\\ \vdots \\ (L_{F}^{\hat r-1}  O_3)(x)\\ \phi_{\bar r + 1}(x)\\ \vdots\\ \phi_{2n}(x) \end{pmatrix},\quad x\in U,
\end{equation}
where $\bar r = p + 2\ell + \hat r m$ and $\phi_i: U \to\R$ for $i=\bar r +1,\ldots,2n$.  Then, in the new coordinates
\[
    \Phi(x) = \begin{pmatrix} \xi\\ \eta\end{pmatrix} = \begin{pmatrix} \xi_1\\ \xi_2\\ \dot \xi_2\\ \xi_3\\ \dot \xi_3\\ \vdots\\ \xi_3^{(\hat r-1)}\\ \eta\end{pmatrix},\quad \xi_1\in\R^p,\ \xi_2\in\R^\ell,\ \xi_3\in\R^m,
\]
the system has the form
\begin{equation}\label{eq:BIF-aux}
\begin{aligned}
  \dot \xi_1 &= b_1(\xi,\eta) + \Gamma_1(\Phi^{-1}(\xi,\eta)) u_{\rm aux},\\
  \ddot \xi_2 &= b_2(\xi,\eta) + \Gamma_2(\Phi^{-1}(\xi,\eta)) u_{\rm aux},\\
  \xi_3^{(\hat r)} &= b_3(\xi,\eta) + \Gamma_3(\Phi^{-1}(\xi,\eta)) u_{\rm aux},\\
  \dot \eta &= q(\xi,\eta) + p(\xi,\eta) u_{\rm aux}.
\end{aligned}
\end{equation}
If the distribution~$x\mapsto \cK(x) = \im K(x)$ is involutive, then the functions~$\phi_i$ can be chosen such that, for all $x\in U$,
\[
   (L_{K} \phi_i)(x) = \phi_i'(x) K(x) = 0\quad \text{and}\quad \Phi'(x)\in\Gl_{2n},
\]
by which $p(\cdot) = 0$, cf.\ Section~\ref{Ssec:RelDeg}. For the auxiliary ODE~\eqref{eq:auxODE} we have the following result.

\begin{lemma} \label{Lem:K-involutive-Holonomic}
For~$K:U\to \R^{2n\times (p+\ell+m)}$ as in \eqref{eq:auxODE}, the distribution~$x\mapsto \cK(x) = \im K(x)$ is involutive on~$U$.
\end{lemma}
\begin{proof}
First, recall that a distribution $\cD(\cdot) = \spn \{ d_1(\cdot),..., d_q(\cdot) \}$ is involutive if, and only if, for any $0 \le i,j \le q$ and $x\in U$, we have
\begin{equation} \label{eq:rank-condition-involutivity}
\rk \big[ d_1(x), ..., d_q(x) \big] = \rk \big[d_1(x), ..., d_q(x), [ d_i, d_j](x) \big],
\end{equation}
where~$[d_i, d_j](x)$ denotes the Lie bracket of~$d_i$ and~$d_j$ at~$x$. Observe that $K$ in~\eqref{eq:auxODE} has the following structure:
\begin{align*}
K(x) &= \big[\tilde k_1(x), ..., \tilde k_{p+\ell+m}(x)\big] \\
&= \begin{bmatrix} \begin{pmatrix} 0 \\ k_1(x_1) \end{pmatrix},...,\begin{pmatrix} 0 \\ k_{p+\ell+m}(x_1) \end{pmatrix}  \end{bmatrix}
\end{align*}
for all $x = (x_1^\top, x_2^\top)^\top \in U$, where $k_i: [I_n,0] U \to \R^n$ for $i=1,...,p+\ell+m$. Now, calculate for any pair~$\tilde k_i(x), \tilde k_j(x)$, $0\le i,j \le p+\ell+m$ the Lie bracket~$[\tilde k_i, \tilde k_j](x)$:
\begin{align*}
[\tilde k_i, \tilde k_j](x) &= \begin{bmatrix} 0 & 0 \\ \frac{\partial  k_i(x_1)}{\partial x_1} & 0 \end{bmatrix} \begin{smallpmatrix} 0 \\ k_j(x_1) \end{smallpmatrix}
-  \begin{bmatrix} 0 & 0 \\ \frac{\partial k_j(x_1)}{\partial x_1}  & 0 \end{bmatrix} \begin{smallpmatrix} 0 \\ k_i(x_1) \end{smallpmatrix} = 0
\end{align*}
for all $x = (x_1^\top, x_2^\top)^\top \in U$. Therefore, using~\eqref{eq:rank-condition-involutivity} the distribution~$\cK(\cdot)$ is involutive on~$U$.\qed
\end{proof}

Lemma~\ref{Lem:K-involutive-Holonomic} yields that it is always possible to choose the diffeomorphism~$\Phi$ such that~$p(\cdot)=0$ in~\eqref{eq:BIF-aux}, which we assume is true henceforth.

We are now in the position to invoke the original constraints $O_1(x) = 0$ and $O_2(x) = 0$, which in the new coordinates simply read $\xi_1 = 0$ and $\xi_2 = 0$. Therefore, we obtain the system
\begin{equation}\label{eq:MBS-BIF}
\begin{aligned}
  0 &= b_1(0,\xi_3,\dot \xi_3,\ldots,\xi_3^{(\hat r-1)},\eta)\\
   &\quad + \Gamma_1(\Phi^{-1}(0,\xi_3,\dot \xi_3,\ldots,\xi_3^{(\hat r-1)},\eta)) u_{\rm aux},\\
  0 &= b_2(0,\xi_3,\dot \xi_3,\ldots,\xi_3^{(\hat r-1)},\eta)\\
   &\quad + \Gamma_2(\Phi^{-1}(0,\xi_3,\dot \xi_3,\ldots,\xi_3^{(\hat r-1)},\eta)) u_{\rm aux},\\
  \xi_3^{(\hat r)} &= b_3(0,\xi_3,\dot \xi_3,\ldots,\xi_3^{(\hat r-1)},\eta) \\
  &\quad + \Gamma_3(\Phi^{-1}(0,\xi_3,\dot \xi_3,\ldots,\xi_3^{(\hat r-1)},\eta)) u_{\rm aux},\\
  \dot \eta &= q(0,\xi_3,\dot \xi_3,\ldots,\xi_3^{(\hat r-1)},\eta),\\
  y &= \xi_3,
\end{aligned}
\end{equation}
which is equivalent to the original multibody system~\eqref{eq:MBS}.

Considering the first two equations in~\eqref{eq:MBS-BIF} we obtain that
\begin{align*}
    &\begin{bmatrix} J(\bar x_1)M(\bar x_1)^{-1}J(\bar x_1)^\top & J(\bar x_1)M(\bar x_1)^{-1}G(\bar x_1)^\top \\
    G(\bar x_1)M(\bar x_1)^{-1}J(\bar x_1)^\top & G(\bar x_1)M(\bar x_1)^{-1}G(\bar x_1)^\top \end{bmatrix} \begin{pmatrix} \mu \\ \lambda\end{pmatrix}\\
     &= - \begin{pmatrix} b_1(0,y,\dot y,\ldots,y^{(\hat r-1)},\eta) \\ b_2(0,y,\dot y,\ldots,y^{(\hat r-1)},\eta)\end{pmatrix} - \begin{bmatrix} J(\bar x_1)\\ G(\bar x_1)\end{bmatrix} M(\bar x_1)^{-1} B(\bar x_1) u,
\end{align*}
where $\bar x_1 = [I_n,0] \Phi^{-1}(0,y,\dot y,\ldots,y^{(\hat r-1)},\eta)$. Since $\Gamma(x)$ can be written as
\[
    \Gamma(x) = \begin{bmatrix} O_1'(x)\\  \big(L_{F}  O_2\big)'(x) \\ \big(L_{F}^{\hat r-1}  O_3\big)'(x)\end{bmatrix} K(x)
\]
and is invertible, it follows that $K(x)$ has full column rank for all $x\in U$. Therefore,~$\begin{bmatrix} J(x_1)^\top & G(x_1)^\top & B(x_1)\end{bmatrix}$ has full column rank and the matrix~$\begin{bmatrix} J(x_1)^\top & G(x_1)^\top\end{bmatrix}$ has full column rank as well, hence, by using that $M$ is pointwise positive definite, we have that $A(x_1) := \begin{smallbmatrix} J(x_1)\\ G(x_1)\end{smallbmatrix} M(x_1)^{-1} \begin{bmatrix} J(x_1)^\top & G(x_1)^\top\end{bmatrix}$ is positive definite for all $x_1\in [I_n,0] U$. Therefore, we can solve for the Lagrange multipliers as follows:
\begin{align*}
    &\begin{pmatrix} \mu \\ \lambda\end{pmatrix} = -A(\bar x_1)^{-1} \left( \begin{pmatrix} b_1(0,y,\dot y,\ldots,y^{(\hat r-1)},\eta) \\ b_2(0,y,\dot y,\ldots,y^{(\hat r-1)},\eta)\end{pmatrix} \right.\\
    &\qquad\quad \left. + \begin{bmatrix} J(\bar x_1)\\ G(\bar x_1)\end{bmatrix} M(\bar x_1)^{-1} B(\bar x_1) u\right) \\
    &= b_4(y,\dot y,\ldots,y^{(\hat r-1)},\eta) - A(\bar x_1)^{-1} \begin{bmatrix} J(\bar x_1)\\ G(\bar x_1)\end{bmatrix} M(\bar x_1)^{-1} B(\bar x_1) u
\end{align*}
for some appropriately chosen function~$b_4$. Inserting this into the third equation in~\eqref{eq:MBS-BIF} gives
\begin{equation}\label{eq:MBS-solved-Lagrange}
\begin{aligned}
  y^{(\hat r)} &= b_3(0,y,\dot y,\ldots,y^{(\hat r-1)},\eta) + \Gamma_3(\bar x) \begin{pmatrix} \mu\\ \lambda \\ u\end{pmatrix}\\
  &= b_3(0,y,\dot y,\ldots,y^{(\hat r-1)},\eta) +\Gamma_3(\bar x) \begin{pmatrix} b_4(y,\dot y,\ldots,y^{(\hat r-1)},\eta)\\ 0\end{pmatrix} \\
  &\quad + \underset{=:S(\bar x)}{\underbrace{\Gamma_3(\bar x) \begin{bmatrix}  -A(\bar x_1)^{-1} \begin{bmatrix} J(\bar x_1)\\ G(\bar x_1)\end{bmatrix} M(\bar x_1)^{-1} B(\bar x_1)\\ I_m\end{bmatrix}}} u,
\end{aligned}
\end{equation}
where $\bar x = \Phi^{-1}(0,y,\dot y,\ldots,y^{(\hat r-1)},\eta)$ and $\bar x_1 = [I_n, 0] \bar x$. Note that it is not clear in general whether the matrix~$S(x)$ is invertible for some or all $x\in U$.
The multibody system~\eqref{eq:MBS} is then equivalent to
\begin{equation}\label{eq:MBS-decoupled}
\begin{aligned}
  y^{(\hat r)} &= b_5(y,\dot y,\ldots,y^{(\hat r-1)},\eta) + \tilde S(y, \dot y,\ldots,y^{(\hat r-1)}, \eta) u,\\
  \dot \eta &= q(y,\dot y,\ldots,y^{(\hat r-1)}, \eta),
\end{aligned}
\end{equation}
where~$b_5$ is some appropriate function, $\tilde S(\cdot) = S\big(\Phi^{-1}(0,\cdot)\big)$ and the internal dynamics are completely decoupled in~\eqref{eq:MBS-decoupled}.


\begin{remark}
  We note that an alternative to the above described approach would be to simply differentiate the holonomic and nonholonomic constraints in~\eqref{eq:MBS} and solve for the respective Lagrange multipliers. When these have been inserted in~\eqref{eq:MBS} (and the still present holonomic and nonholonomic constraints are neglected) an ODE of the form~\eqref{eq:System} is obtained, for which the Byrnes-Isidori form as in Section~\ref{Ssec:RelDeg} can be computed. This approach has some disadvantages compared to our approach presented above:
  \begin{enumerate}[1)]
    \item Since differentiating the constraints does not remove them, it must be guaranteed that any solution of the resulting ODE is also a solution of the overall DAE, i.e., it satisfies the constraints for all times. Nevertheless, this can be ensured, provided that the initial values are chosen appropriately.
    \item It cannot be avoided to solve for the Lagrange multipliers first, while in our approach the internal dynamics in~\eqref{eq:MBS-decoupled} are the same as obtained for the auxiliary system in~\eqref{eq:BIF-aux}, since~$p(\cdot)=0$ in the latter by Lemma~\ref{Lem:K-involutive-Holonomic}. The system~\eqref{eq:BIF-aux} can be computed without solving for the Lagrange multipliers and hence requires much less computational effort. Nevertheless, it is still important to choose the functions $\phi_{\bar r+1},\ldots,\phi_{2n}$ such that $L_K \phi_i = 0$ for $i=\bar r+ 1,\ldots,2n$; this problem is discussed in Section~\ref{Sec:intdyn-coord}.
  \end{enumerate}
\end{remark}

\subsection{Position-dependent output with relative degree two}\label{Sec:r=2}

In the special case $h(q,v) = h(q)$, i.e., $h:\R^n\to\R^m$ and $y(t) = h\big(q(t)\big)$, we obtain some additional structure. First of all, we may compute that
\[
  (L_{K} O_3)(x) = O_3'(x)K(x) = [h'(x_1), 0] K(x) = 0,
\]
whence $\hat r \ge 2$. Assuming that $\hat r = 2$ (a typical situation) we obtain, for all $x=(x_1^\top, x_2^\top)^\top\in U$,
\begin{align*}
    \Gamma(x) &= \begin{bmatrix}  (L_{K} O_1)(x)\\ (L_{K} L_{F}  O_2)(x)\\ (L_{K} L_{F}  O_3)(x)\end{bmatrix} \\
    &= \begin{bmatrix} J(x_1)\\
    G(x_1)\\ h'(x_1)\end{bmatrix}  M(x_1)^{-1} \begin{bmatrix} J(x_1)^\top & G(x_1)^\top & B(x_1)\end{bmatrix}.
\end{align*}
Furthermore,
\begin{align*}
   & S(x) = h'(x_1) M(x_1)^{-1} \bigg( M(x_1)  \\
   & - \begin{bmatrix} J(x_1)^\top & G(x_1)^\top\end{bmatrix} A(x_1)^{-1} \begin{bmatrix} J(x_1)\\ G(x_1)\end{bmatrix}  \bigg) M(x_1)^{-1} B(x_1),
\end{align*}
which is the Schur complement of~$A(x_1)$ in~$\Gamma(x)$.

\subsection{The colocated case}\label{Sec:colocate}

Now we consider the special case that the inputs and outputs of system~\eqref{eq:MBS} are colocated, which means that $h(q,v) = h(q)$ as in Section~\ref{Sec:r=2} and
\begin{equation}\label{eq:colocate}
     \forall\, q\in\R^n:\ h'(q) = B(q)^\top.
\end{equation}
We show that in this case the vector relative degree with respect to~$O_3$ is $\hat r = 2$. First we calculate, for all $x=(x_1^\top, x_2^\top)^\top\in U$,
\[
  (L_{K} O_3)(x) = O_3'(x)K(x) = [h'(x_1), 0] K(x) = 0
\]
and
\begin{align*}
  &(L_{K} L_{F} O_3)(x) = \big(O_3'(x) F(x)\big)' K(x)\\
  &= \left([B(x_1)^\top, 0] \begin{pmatrix} x_2\\ M(x_1)^{-1} f(x_1,x_2) \end{pmatrix}\right)' K(x)\\
  &= \begin{bmatrix} \frac{\partial}{\partial x_1} \big(B(x_1)^\top x_2\big) & B(x_1)^\top\end{bmatrix} K(x)\\
  &= B(x_1)^\top M(x_1)^{-1} \begin{bmatrix} J(x_1)^\top & G(x_1)^\top & B(x_1)\end{bmatrix}.
\end{align*}
Then we obtain the high-gain matrix
\begin{multline*}
    \Gamma(x) = \begin{bmatrix} \Gamma_1(x)\\ \Gamma_2(x)\\ \Gamma_3(x)\end{bmatrix} = \begin{bmatrix} J(x_1)\\ G(x_1)\\ B(x_1)^\top\end{bmatrix} M(x_1)^{-1} \begin{bmatrix} J(x_1)\\ G(x_1)\\ B(x_1)^\top\end{bmatrix}^\top \\ \in\R^{(p+\ell+m)\times (p+\ell+m)}.
\end{multline*}
Therefore, we find that, for some $x = (x_1^\top, x_2^\top)^\top \in U$,
\[
    \Gamma(x) > 0\ \ \iff \ \ \rk \begin{bmatrix} J(x_1)^\top & G(x_1)^\top & B(x_1)\end{bmatrix} = p+\ell+m,
\]
i.e., the latter matrix has full column rank. Roughly speaking, this condition means that (the components of) the Lagrange multipliers~$\mu$ and~$\lambda$ and the inputs~$u$ influence the system in a linearly independent (i.e., non-redundant) way, because in
\[
    M(q) \dot v= f(q,v) +  \begin{bmatrix} J(q)^\top & G(q)^\top & B(q)\end{bmatrix} \begin{pmatrix} \mu\\ \lambda\\ u\end{pmatrix}
\]
the matrix $\begin{bmatrix} J(q)^\top & G(q)^\top & B(q)\end{bmatrix}$ has full column rank, thus the auxiliary input~$u_{\rm aux}$ does not have any redundant components.

Then the transformation which puts~\eqref{eq:auxODE} into Byrnes-Isidori form is given by
\[
    \Phi(x) = \begin{pmatrix} O_1(x)\\ O_2(x)\\ O_2'(x) F(x) \\ O_3(x) \\ O_3'(x) F(x)\\ \phi_{\bar r + 1}(x)\\ \vdots\\ \phi_{2n}(x) \end{pmatrix}  = \begin{pmatrix}
        J(x_1) x_2 + j(x_1) \\ g(x_1) \\  G(x_1) x_2\\ h(x_1)\\ B(x_1)^\top x_2\\ \phi_{\bar r + 1}(x)\\ \vdots\\ \phi_{2n}(x) \end{pmatrix},
\]
where $\bar r = p + 2\ell + 2m$ and, since $\cK(\cdot)$ is involutive by Lemma~\ref{Lem:K-involutive-Holonomic}, it is possible to choose $\phi_i:U\to\R$, $i=\bar r +1,\ldots,2n$, such that $(L_{K} \phi_i)(x) = \phi_i'(x) K(x) = 0$ and $\Phi'(x)$ is invertible for all $x\in U$. Using the new coordinates
\[
    \Phi(x) = \begin{pmatrix} \xi\\ \eta\end{pmatrix} = \begin{pmatrix} \xi_1\\ \xi_2\\ \dot \xi_2\\ \xi_3\\ \dot \xi_3\\ \eta\end{pmatrix}
\]
and invoking the original constraints $O_1(x) = 0$ and $O_2(x) = 0$ we may rewrite the original multibody system~\eqref{eq:MBS} in the form~\eqref{eq:MBS-BIF}, which in the colocated case simplifies to
\begin{equation}\label{eq:MBS-BIF-col}
\begin{aligned}
  0 &= b_1(0,\xi_3,\dot \xi_3,\eta) + \Gamma_1(\Phi^{-1}(0,\xi_3,\dot \xi_3,\eta)) u_{\rm aux},\\
  0 &= b_2(0,\xi_3,\dot \xi_3,\eta) + \Gamma_2(\Phi^{-1}(0,\xi_3,\dot \xi_3,\eta)) u_{\rm aux},\\
  \ddot \xi_3 &= b_3(0,\xi_3,\dot \xi_3,\eta) + \Gamma_3(\Phi^{-1}(0,\xi_3,\dot \xi_3,\eta)) u_{\rm aux},\\
  \dot \eta &= q(0,\xi_3,\dot \xi_3,\eta),\\
  y &= \xi_3.
\end{aligned}
\end{equation}

As in the general case we may solve for the Lagrange multipliers $\mu$ and $\lambda$ and insert this in~\eqref{eq:MBS-BIF-col}, thus~\eqref{eq:MBS-solved-Lagrange} becomes
\begin{align*}
  \ddot y &= b_3(0,y,\dot y,\eta)\\
   &\quad + B(\bar x_1)^\top M(\bar x_1)^{-1} \begin{bmatrix} J(\bar x_1)^\top & G(\bar x_1)^\top & B(\bar x_1)\end{bmatrix} \begin{pmatrix} \mu\\ \lambda \\ u\end{pmatrix}\\
  &= b_3(0,y,\dot y,\eta) + B(\bar x_1)^\top M(\bar x_1)^{-1} \begin{bmatrix} J(\bar x_1) \\ G(\bar x_1)\end{bmatrix}^\top b_4(y,\dot y,\eta)\\
  &\quad + S(\bar x)\ u,
\end{align*}
where $\bar x = \Phi^{-1}(0,y,\dot y,\ldots,y^{(\hat r-1)},\eta)$, $\bar x_1 = [I_n, 0] \bar x$ and
\begin{align*}
   & S(x) = B(x_1)^\top M(x_1)^{-1} \bigg( M(x_1)  \\
   & - \begin{bmatrix} J(x_1)^\top & G(x_1)^\top\end{bmatrix} A(x_1)^{-1} \begin{bmatrix} J(x_1)\\ G(x_1)\end{bmatrix}  \bigg) M(x_1)^{-1} B(x_1)
\end{align*}
for  $x = (x_1^\top, x_2^\top)^\top \in U$. In the colocated case, the matrix~$S(x)$ is the Schur complement of~$A(x_1)$ in~$\Gamma(x)$. It is a well-known result that since~$\Gamma(x)$ is positive definite, then the Schur complement~$S(x_1)$ is positive definite as well. Therefore, the decoupled system~\eqref{eq:MBS-decoupled} is of the form
\begin{equation}\label{eq:MBS-decoupled-col}
\begin{aligned}
  \ddot y &= b_5(y,\dot y,\eta) + \tilde S(y, \dot y, \eta) u,\\
  \dot \eta &= q(y,\dot y, \eta),
\end{aligned}
\end{equation}
where~$\tilde S$ is pointwise positive definite.

\section{A feasible set of coordinates for the internal dynamics}\label{Sec:intdyn-coord}
In this section we derive a representation of the internal dynamics which depends on the auxiliary output~$y_{\rm aux}$ and its derivative~$\dot y_{\rm aux}$.
Consider a system~\eqref{eq:MBS} with holonomic and nonholonomic constraints and position-dependent output. Further set $H(x_1) =: h'(x_1)$ and assume that, for some open set $U_1 \subseteq \R^n$,
\begin{equation}\label{eq:ass-gamma}
   \begin{bmatrix} J(x_1) \\ G(x_1) \\ H(x_1) \end{bmatrix} M(x_1)^{-1} \begin{bmatrix} J(x_1)^\top & G(x_1)^\top & B(x_1) \end{bmatrix} \in \Gl_{p+\ell+m}
\end{equation}
for all $x_1\in U_1$, which means that the high-gain matrix $\Gamma(\cdot)$ is invertible on $U := U_1 \times \R^n$. Then, following Section~\ref{Sec:r=2}, the auxiliary ODE~\eqref{eq:auxODE} with (again omitting the argument~$t$ for brevity)  $u_{\rm aux}^\top = (\mu^\top, \lambda^\top, u^\top)$ and $y_{\rm aux}^\top = \big((J(x_1)x_2 + j(x_1))^\top, g(x_1)^\top, h(x_1)^\top\big)$ has vector relative degree $(1,\ldots,1,2,\ldots,2)$ on~$U$, thus~$\bar r = p+2\ell+2m$. The diffeomorphism~$\Phi:U\to W$ can be represented as
\begin{equation} \label{eq:Phi}
    \Phi(x) = \begin{pmatrix} \xi \\ \eta  \end{pmatrix}  =   \begin{pmatrix} O_1(x) \\ O_2(x) \\ O_2'(x) F(x)\\ O_3(x) \\  O_3'(x) F(x) \\ \phi_{\bar r + 1}(x)\\ \vdots\\ \phi_{2n}(x) \end{pmatrix}
 =   \begin{pmatrix}  J(x_1)x_2 + j(x_1) \\       g(x_1) \\  G(x_1) x_2  \\  h(x_1) \\ H(x_1) x_2\\ \phi_{\bar r + 1}(x)\\ \vdots\\ \phi_{2n}(x) \end{pmatrix},
\end{equation}
$x = (x_1^\top, x_2^\top)^\top \in U$. For the internal dynamics $\eta = \big(\phi_{\bar r+1}(x),\ldots,\phi_{2n}(x)\big)^\top$ we make, with some abuse of notation, the structural ansatz
\begin{align} \label{eq:eta-structure}
 \eta = \begin{pmatrix} \eta_1 \\ \eta_2  \end{pmatrix} = \begin{pmatrix} \phi_1(x_1) \\ \phi_2(x_1) x_2  \end{pmatrix},
\end{align}
where $\phi_1\in\cC^1(U_1 \to \R^{n-\ell-m})$ and $\phi_2\in\cC(U_1 \to \R^{(n-p-\ell-m) \times n})$. If we now rearrange~\eqref{eq:Phi}, i.e., we set
\begin{equation}\label{eq:xi-partition}
\begin{aligned}
 \xi_1 &= \begin{pmatrix} g(x_1) \\  h(x_1) \end{pmatrix},\quad \xi_2 = \begin{bmatrix} J(x_1) \\ G(x_1) \\ H(x_1) \end{bmatrix} x_2,\\
 \xi &= P \begin{pmatrix} \xi_1 \\ \xi_2 \end{pmatrix} + \begin{pmatrix} j(x_1) \\  0 \end{pmatrix},
\end{aligned}
\end{equation}
where $P\in\R^{\bar r\times \bar r}$ is some permutation matrix,
we see that $(\xi_1^\top, \xi_2^\top)^\top$ and $(\eta_1^\top, \eta_2^\top)^\top$ are of similar structure. Since~$\Phi$ is a diffeomorphism and system~\eqref{eq:auxODE} has vector relative degree $(1,\ldots,1,2,\ldots,2)$ on~$U$, its Jacobian is invertible on~$U$:
\begin{align} \label{eq:Phi-invertible}
&\forall\, x \!=\! (x_1^\top, x_2^\top)^\top \!\in\! U\!: \Phi'(x)\!=\! \begin{bmatrix} * + j'(x_1) & J(x_1) \\ G(x_1) & 0 \\ * & G(x_1) \\ H(x_1) & 0 \\ * & H(x_1) \\ \phi_1'(x_1) & 0 \\ * & \phi_2(x_1) \end{bmatrix} \!\in\! \Gl_{2n}   \nonumber \\
&\iff\ \ \forall\, x_1\in U_1 \subseteq \R^n: \nonumber \\
&\qquad\quad \begin{bmatrix} G(x_1) \\ H(x_1) \\ \phi'_1(x_1) \end{bmatrix} \in \Gl_{n} \
\wedge \  \begin{bmatrix} J(x_1) \\ G(x_1) \\ H(x_1) \\ \phi_2(x_1) \end{bmatrix} \in \Gl_{n},
\end{align}
where $*$ is of the form $\tfrac{\partial}{\partial x_1} \big(\zeta(x_1)  x_2\big)$ with appropriate $\zeta: U_1 \to \R^{q \times n}$ and $q \in \N$.
Since we are interested in the case~${p(\cdot) = 0}$ in~\eqref{eq:BIF-aux} (which can be achieved by Lemma~\ref{Lem:K-involutive-Holonomic}) we aim to find~$\phi_2: U_1 \to \R^{(n-p-\ell-m) \times n}$ such that
\begin{align} \label{eq:K-involutive}
\forall x_1 \in U_1: \phi_2(x_1) M(x_1)^{-1} \begin{bmatrix} J(x_1)^\top & G(x_1)^\top & B(x_1) \end{bmatrix} = 0.
\end{align}
We summarize this in the following result.

%

\begin{lemma}\label{Lem:existence-phi1/2}
Consider the 
multi-body system~\eqref{eq:MBS} and assume that~\eqref{eq:ass-gamma} holds on an open set $U_1\subseteq\R^n$. For any~$x_1^0 \in U_1$ there exists an open neighbourhood $U_1^0 \subseteq U_1$ of~$x_1^0$ and $\phi_1\in\cC^1(U_1^0 \to \R^{n-\ell-m})$, $\phi_2\in\cC(U_1^0 \to \R^{(n-p-\ell-m) \times n})$
 such that~\eqref{eq:Phi-invertible} and~\eqref{eq:K-involutive} hold locally on~$U_1^0$, resp.
\end{lemma}

\begin{proof}
We exploit~\cite[Lem. 4.1.5]{Sont98a} which states the following: Consider~$W \in \cC( U_1 \to \R^{w \times n})$ with  $\rk W(x_1) = w$ for all~$x_1 \in U_1$. Then for each~$x_1^0 \in U_1$ there exist an open neighbourhood $U_1^0 \subseteq U_1$ and~$T \in \cC(U_1^0 \to \Gl_n)$ such that
\[
\forall\, x_1 \in U_1^0: \ W(x_1) T(x_1) = \begin{bmatrix} I_w & 0 \end{bmatrix}.
\]
Now, we use this to show the existence of~$\phi_1$. Since by~\eqref{eq:ass-gamma} we have~$\rk [G(x_1)^\top, H(x_1)^\top] = \ell + m$ for all~${x_1 \in U_1}$ there exist for each~$x_1^0 \in U_1$ an open neighbourhood~$U_1^0 \subseteq U_1$ and~$T = [T_1, T_2] \in \cC(U_1^0 \to \Gl_n)$ such that
\begin{equation*}
\forall\, x_1 \in U_1^0: \
\begin{bmatrix}
G(x_1) \\ H(x_1)
\end{bmatrix}
\begin{bmatrix}
T_1(x_1) & T_2(x_1)
\end{bmatrix}
= \begin{bmatrix}
I_{\ell + m} & 0
\end{bmatrix},
\end{equation*}
i.e., $\im T_2(x_1) = \ker [G(x_1)^\top , H(x_1)^\top]^\top$ and $\rk T_2(x_1) = n-\ell-m$ for all~$x_1 \in U_1^0$.
{Define $E :=T_2(x_1^0)^\top$, then
\begin{equation*}
\begin{bmatrix}
G(x_1) \\ H(x_1) \\ E
\end{bmatrix}
\begin{bmatrix}
T_1(x_1) & T_2(x_1)
\end{bmatrix}
= \begin{bmatrix}
I_{\ell + m} & 0 \\ * & E T_2(x_1)
\end{bmatrix}.
\end{equation*}
Since $\rk T_2(x_1^0) = n-\ell-m$ it follows that for $x_1=x_1^0$ we have $E T_2(x_1^0) = T_2(x_1^0)^\top T_2(x_1^0) \in \Gl_{n-\ell-m}$.} Furthermore, by~$T_2 \in \cC(U_1^0 \to \R^{n \times (n-\ell-m)})$, the mapping~$x_1 \mapsto \det(ET_2(x_1))$ is continuous on~$U_1^0$, hence there is an open neighbourhood~$\bar U_1^0 \subseteq U_1^0$ of~$x_1^0$ such that~$\det(ET_2(\bar x_1)) \neq 0$ for all~$\bar x_1 \in \bar U_1^0$. Thus,
\begin{equation*}
\forall\, x_1 \in \bar U_1^0: \
\rk \begin{bmatrix}
G(x_1) \\ H(x_1) \\ E
\end{bmatrix} = n
\end{equation*}
and, with $\phi_1 :  \bar U_1^0 \to \R^{n-\ell-m},\ x_1 \mapsto E x_1$, we have $\phi_1 \in \cC^1(\bar U_1^0 \to \R^{n-\ell-m})$ and the first condition in~\eqref{eq:Phi-invertible} is satisfied on~$\bar U_1^0$ since $\phi_1'(x_1) = E$.
\ \\
Now, we show the existence of~$\phi_2$. Observe that by~\eqref{eq:ass-gamma} we have~$ \rk [J(x_1)^\top, G(x_1)^\top, B(x_1)] = p+\ell + m$ and~$\rk M(x_1) = n$ for all~$x_1 \in U_1$. Therefore, again via~\cite[Lem.~4.1.5]{Sont98a},
there exist for each~$x_1^0 \in U_1$ an open neighbourhood~$U_1^0 \subseteq U_1$ and~$T = [T_1, T_2] \in \cC(U_1^0 \to \Gl_n)$ such that
\begin{equation*}
\forall \, x_1 \in U_1^0: \, \begin{bmatrix}J(x_1) \\ G(x_1) \\ B(x_1)^\top \end{bmatrix} \begin{bmatrix} T_1(x_1) & T_2(x_1)\end{bmatrix}
= \begin{bmatrix} I_{p+\ell + m} & 0 \end{bmatrix},
\end{equation*}
i.e., $\im T_2(x_1) = \ker [J(x_1)^\top , G(x_1)^\top , B(x_1)]^\top$ and $T_2 \in \cC(U_1^0 \to \R^{n \times (n-p-\ell-m)})$. Now, choosing
${\phi_2(\cdot) = T_2(\cdot)^\top M(\cdot)}$ we obtain
\begin{equation*}
\begin{aligned}
\begin{bmatrix} J(x_1) \\ G(x_1) \\ H(x_1) \\ \phi_2(x_1) \end{bmatrix}
\begin{bmatrix} M(x_1)^{-1} \begin{bmatrix} J(x_1)^\top & G(x_1)^\top & B(x_1) \end{bmatrix} & \phi_2(x_1)^\top \end{bmatrix} \\
=
\begin{bmatrix} \Gamma(x_1) & * \\ 0 & \phi_2(x_1) \phi_2(x_1)^\top  \end{bmatrix} \in \R^{n \times n}, \  x_1 \in U_1^0,
\end{aligned}
\end{equation*}
which is invertible on~$U_1^0$ since by~\eqref{eq:ass-gamma} the high-gain~$\Gamma(x_1)$ is invertible, and $\rk T_2(x_1) = n-p-\ell-m$ for all~$x_1 \in U_1^0$. Hence the second condition in~\eqref{eq:Phi-invertible} is satisfied on~$U_1^0$. Furthermore, equation~\eqref{eq:K-involutive} holds on~$U_1^0$ by construction of~$\phi_2$.
\qed
\end{proof}

Lemma~\ref{Lem:existence-phi1/2} justifies the structural ansatz for the internal state~$\eta$ in~\eqref{eq:eta-structure}. Note that the construction of~$\phi_1$ in the proof of Lemma~\ref{Lem:existence-phi1/2} is not unique. There is a lot of freedom in the choice of~$\phi_1$, {which is basically free up to~\eqref{eq:Phi-invertible}. However, a straightforward choice is $\phi_1(x_1) = Ex_1$ with $E = T_2(x_1^0)^\top$ as in the proof of Lemma~\ref{Lem:existence-phi1/2}, which can be computed via $\im T_2(x_1^0) = \ker [G(x_1^0)^\top, H(x_1^0)^\top]^\top$.}

On the other hand,~$\phi_2$ is uniquely determined up to an invertible left transformation. To find all possible representations, assume that $\phi_2\in\cC(U_1 \to \R^{(n-p-\ell-m) \times n})$ is such that~\eqref{eq:Phi-invertible} and~\eqref{eq:K-involutive} hold. Then there exist $P: U_1\to \R^{n \times (p+\ell+m)}$ and $V: U_1\to \R^{n \times (n-p-\ell-m)}$ such that
\begin{equation}\label{eq:PV-inverse}
\forall\, x_1\in U_1:\ [P(x_1), V(x_1)] \begin{bmatrix} \begin{bmatrix} J(x_1) \\ G(x_1) \\ H(x_1) \end{bmatrix} \\ \phi_2(x_1) \end{bmatrix} = I_n.
\end{equation}
Furthermore, $P, V$ have pointwise full column rank, by which the pseudoinverse of~$V$ is given by $V^\dagger(x_1) = (V(x_1)^\top V(x_1))^{-1} V(x_1)^\top$ for $x_1\in U_1$. {Note that it follows from~\eqref{eq:PV-inverse} that~$V$ is the pointwise basis matrix of a kernel, i.e.,
\[
    \forall\, x_1\in U_1:\ \im V(x_1) =  \ker  \begin{bmatrix} J(x_1) \\ G(x_1) \\ H(x_1) \end{bmatrix}.
\]
}

\begin{lemma}\label{Lem:phi_2}
Let $\phi_2\in\cC(U_1 \to \R^{(n-p-\ell-m) \times n})$ be such that~\eqref{eq:Phi-invertible} and~\eqref{eq:K-involutive} hold. Then~$\phi_2$ is uniquely determined up to an invertible left transformation. All possible functions are given by
\begin{equation}\label{eq:phi2-def}
\begin{aligned}
\phi_2(x_1) &=  V^\dagger(x_1) \bigg(I_n -  \\
& M(x_1)^{-1} \begin{bmatrix} J(x_1)^\top & G(x_1)^\top & B(x_1) \end{bmatrix} \Gamma(x_1)^{-1} \begin{smallbmatrix} J(x_1) \\ G(x_1) \\ H(x_1) \end{smallbmatrix}  \bigg),
\end{aligned}
\end{equation}
$x_1\in U_1$, for feasible choices of~$V$ satisfying~\eqref{eq:PV-inverse}.
\end{lemma}
\begin{proof} Assume that~\eqref{eq:Phi-invertible} and~\eqref{eq:K-involutive} hold, and hence we have~\eqref{eq:PV-inverse} for some corresponding~$P$ and~$V$. Observe that multiplying~\eqref{eq:PV-inverse} from the left by~$V(x_1)^\dagger$ gives
\[
\phi_2(x_1) = V^\dagger(x_1) \left(I_n - P(x_1) \begin{bmatrix} J(x_1) \\ G(x_1) \\ H(x_1) \end{bmatrix}  \right),\quad x_1\in U_1.
\]
Invoking~\eqref{eq:K-involutive}, we further obtain from~\eqref{eq:PV-inverse} that
\begin{align*}
P(x_1) & =  M(x_1)^{-1} \begin{bmatrix} J(x_1)^\top & G(x_1)^\top & B(x_1) \end{bmatrix} \\
&\quad \cdot \underset{= \Gamma(x_1)^{-1}}{\underbrace{\left( \begin{bmatrix} J(x_1) \\ G(x_1) \\ H(x_1) \end{bmatrix} M(x_1)^{-1} \begin{bmatrix} J(x_1)^\top & G(x_1)^\top & B(x_1) \end{bmatrix}  \right)^{-1}}},
\end{align*}
and hence~$P$ is uniquely determined by~$M, J, G, H, B$. Therefore,~$\phi_2$ is given as in~\eqref{eq:phi2-def}. Furthermore, it follows from~\eqref{eq:PV-inverse} that
\[
 \begin{bmatrix} \begin{bmatrix} J(x_1) \\ G(x_1) \\ H(x_1) \end{bmatrix} \\ \phi_2(x_1) \end{bmatrix} [P(x_1), V(x_1)] = \begin{bmatrix} I_{p+\ell+m} &0\\0& I_{n-p-\ell-m} \end{bmatrix},
\]
from which we may deduce $\phi_2(x_1) V(x_1) = I_{n-p-\ell-m}$ and $\im V(x_1) = \ker [J(x_1)^\top, G(x_1)^\top, H(x_1)^\top]^\top$. Therefore, the representation of~$\phi_2$ in~\eqref{eq:phi2-def} only depends on the choice of the basis of $\ker [J(x_1)^\top, G(x_1)^\top, H(x_1)^\top]^\top$. Now, let $\tilde V(x_1) = V(x_1) R(x_1)$, $x_1\in U_1$, for some $R: U_1\to \Gl_{n-p-\ell-m}$ and set
\begin{align*}
\tilde \phi_2(x_1) &=  \tilde V^\dagger(x_1) \bigg(I_n -  \\
& M(x_1)^{-1} \begin{bmatrix} J(x_1)^\top & G(x_1)^\top & B(x_1) \end{bmatrix} \Gamma(x_1)^{-1} \begin{smallbmatrix} J(x_1) \\ G(x_1) \\ H(x_1) \end{smallbmatrix}  \bigg).
\end{align*}
Then a short calculation shows the claim $\tilde \phi_2(x_1) = R(x_1)^{-1} \phi_2(x_1)$ for all $x_1\in U_1$.
\qed
\end{proof}

Now, let $\phi_1\in\cC^1(U_1 \to \R^{n-\ell-m})$ and $\phi_2\in\cC(U_1 \to \R^{(n-p-\ell-m) \times n})$ be such that~\eqref{eq:Phi-invertible} and~\eqref{eq:K-involutive} are satisfied and $V: U_1\to \R^{n \times (n-p-\ell-m)}$ as in~\eqref{eq:PV-inverse}. We continue deriving a representation of the internal dynamics.
We define $\bar y_{\rm aux} := [0,I_{\ell+m}] y_{\rm aux}$ and observe that, using the first equation in~\eqref{eq:auxODE},
\[
\ddt {\bar y}_{\rm aux} = \begin{bmatrix} G(x_1) \\ H(x_1) \end{bmatrix} x_2.
\]
Then, with $y_{\rm aux, 1} = O_1(x)$ we find
\begin{equation}
\begin{aligned}
\begin{pmatrix} y_{\rm aux,1} \\ \ddt{\bar y}_{\rm aux} \\ \eta_2 \end{pmatrix} &= \begin{bmatrix} J(x_1) \\ \begin{bmatrix}  G(x_1) \\ H(x_1) \end{bmatrix} \\ \phi_2(x_1) \end{bmatrix} x_2
+ \begin{pmatrix} j(x_1) \\ \begin{pmatrix} 0 \\ 0\end{pmatrix}  \\ 0 \end{pmatrix}
\\
\Longrightarrow\quad
x_2 &= \begin{bmatrix}  J(x_1) \\ \begin{bmatrix} G(x_1) \\ H(x_1) \end{bmatrix} \\ \phi_2(x_1) \end{bmatrix}^{-1} \begin{pmatrix} y_{\rm aux,1} - j(x_1) \\ \ddt{ \bar y}_{\rm aux} \\ \eta_2 \end{pmatrix}  \\
&\overset{\eqref{eq:PV-inverse}}{=}  M(x_1)^{-1} \begin{bmatrix} J(x_1)\\ G(x_1)\\  B(x_1)^\top\end{bmatrix}^\top \Gamma(x_1)^{-1} \begin{pmatrix} y_{\rm aux,1} - j(x_1) \\ \ddt{ \bar y}_{\rm aux} \end{pmatrix} \\
& \quad + V(x_1)\eta_2.  \label{eq:x2}
\end{aligned}
\end{equation}
Thus, using~\eqref{eq:eta-structure} and~\eqref{eq:x2} the dynamics of~$\eta_1$ are given by
\begin{equation}
\begin{aligned}
\dot{\eta}_1 &= \phi'_1(x_1) M(x_1)^{-1} \begin{bmatrix} J(x_1)\\ G(x_1)\\  B(x_1)^\top\end{bmatrix}^\top \Gamma(x_1)^{-1} \begin{pmatrix} y_{\rm aux,1} - j(x_1) \\ \ddt{ \bar y}_{\rm aux} \end{pmatrix} \\
&\quad + \phi'_1(x_1) V(x_1) \eta_2 \label{eq:eta1-dynamik_x1_dependent}.
\end{aligned}
\end{equation}
{Furthermore, we have that}
\begin{align*}
\begin{pmatrix} \bar y_{\rm aux} \\ \eta_1 \end{pmatrix} = \begin{pmatrix} \begin{pmatrix} g(x_1) \\ h(x_1) \end{pmatrix} \\ \phi_1(x_1) \end{pmatrix} =: {\vt(x_1)}
\end{align*}
for a continuously differentiable {$\vt: U_1 \to \R^n$}, the Jacobian of which is invertible on~$U_1$ by~\eqref{eq:Phi-invertible}. In order to ensure that~${\vt}$ is a diffeomorphism on~$U_1$ we need to additionally require that there exist a diffeomorphism $\psi: U_1\to \R^n$ and $\omega:[0,\infty)\to (0,\infty)$ which is continuous, nondecreasing and satisfies $\int_0^\infty \frac{1}{\omega(t)} \, {\rm d}t = \infty$ such that
\[
    \forall\, x\in U_1:\ \|\psi'(x) \cdot \big({\vt}'(x)\big)^{-1}\| \le \omega(\|x\|).
\]
Then~\cite[Thm.~2.1]{BergHall20pp} yields that ${\vt}: U_1 \to W_1 := {\vt}(U_1)$ is a diffeomorphism. Then we have
\begin{equation} \label{eq:x1}
x_1 = {\vt}^{-1}\Bigg( \begin{pmatrix} \bar y_{\rm aux} \\ \eta_1 \end{pmatrix} \Bigg).
\end{equation}
To combine~\eqref{eq:x1} with~\eqref{eq:eta1-dynamik_x1_dependent} we define the following functions as concatenations on $W_1$:
{
\begin{align*}
  \phi'_{1,\vt}(\cdot) &:=  \big(\phi'_1 \circ \vt^{-1} \big)(\cdot), &M_{\vt}(\cdot)^{-1} &:= \big(M \circ \vt^{-1}\big)(\cdot)^{-1},\\
J_{\vt}^\top(\cdot) &:= \big(J \circ \vt^{-1} \big)(\cdot)^\top, &j_{\vt}(\cdot) &:= \big(j \circ \vt^{-1} \big)(\cdot),\\
G_{\vt}^\top(\cdot) &:= \big(G \circ \vt^{-1} \big)(\cdot)^\top, &B_{\vt}(\cdot) &:= \big(B \circ \vt^{-1} \big)(\cdot),\\
\Gamma_{\vt}(\cdot)^{-1} &:= \big(\Gamma \circ \vt^{-1} \big)(\cdot)^{-1} &V_{\vt}(\cdot) &:= \big(V \circ \vt^{-1}\big)(\cdot).
\end{align*}
Therefore, we obtain the following representation of~\eqref{eq:eta1-dynamik_x1_dependent}:
\begin{equation*}
\begin{aligned}
\dot \eta_1
&= \phi'_{1, \vt}\left(\begin{smallpmatrix} \bar y_{\rm aux} \\ \eta_1 \end{smallpmatrix} \right)
M_{\vt}\left(\begin{smallpmatrix} \bar y_{\rm aux} \\ \eta_1 \end{smallpmatrix} \right)^{-1}  \\
&\quad \cdot \begin{bmatrix} J_{\vt}\left(\begin{smallpmatrix} \bar y_{\rm aux} \\ \eta_1 \end{smallpmatrix} \right)^\top & G_{\vt}\left(\begin{smallpmatrix} \bar y_{\rm aux} \\ \eta_1 \end{smallpmatrix} \right)^\top & B_{\vt}\left(\begin{smallpmatrix} \bar y_{\rm aux} \\ \eta_1 \end{smallpmatrix} \right)\end{bmatrix} \\
&\quad \cdot \Gamma_{\vt}\left(\begin{smallpmatrix} \bar y_{\rm aux} \\ \eta_1 \end{smallpmatrix} \right)^{-1} \begin{pmatrix} y_{\rm aux,1} - j_\vt \left(\begin{smallpmatrix} \bar y_{\rm aux} \\ \eta_1 \end{smallpmatrix} \right)\\ \ddt{\bar y}_{\rm aux}\end{pmatrix}  \\
&\quad+ \phi'_{1,\vt}\left(\begin{smallpmatrix} \bar y_{\rm aux} \\ \eta_1 \end{smallpmatrix} \right) V_{\vt}\left(\begin{smallpmatrix} \bar y_{\rm aux} \\ \eta_1 \end{smallpmatrix} \right) \eta_2.
\end{aligned}
\end{equation*}
}
Now, we investigate the dynamics of~$\eta_2$. Define $\phi'_2[x_1,x_2] := \tfrac{\partial}{\partial x_1} \big(\phi_2(x_1) x_2\big) \in \R^{(n-p-\ell-m) \times n }$, then from~\eqref{eq:auxODE},~\eqref{eq:eta-structure} and~\eqref{eq:K-involutive} we obtain
\[
\dot \eta_2 = \phi'_2[x_1,x_2]x_2 + \phi_2(x_1) M(x_1)^{-1} f(x_1,x_2).
\]
Let {$\phi_{2,\vt}(\cdot) :=  \big(\phi_2 \circ \vt^{-1} \big)(\cdot)$, $H_{\vt}(\cdot) :=  \big(H \circ \vt^{-1} \big)(\cdot)$} on $W_1$ and define for $w \in W_1$, $v \in \R^n$
{
\begin{align*}
\phi'_{2,\vt}[w,v] &:= \phi'_2 \left[\vt^{-1}(w), \begin{smallbmatrix} J_\vt(w) \\ \begin{smallbmatrix} G_{\vt}(w)\\ H_{\vt}(w) \end{smallbmatrix} \\ \phi_{2,\vt}(w) \end{smallbmatrix}^{-1}  v \right],\\
f_{\vt}(w,v) &:= f \left(\vt^{-1}(w), \begin{smallbmatrix} J_\vt(w) \\ \begin{smallbmatrix} G_{\vt}(w)\\ H_{\vt}(w) \end{smallbmatrix} \\ \phi_{2,\vt}(w) \end{smallbmatrix}^{-1}  v \right).
\end{align*}
Then the dynamics of~$\eta_2$ are given by
\begin{equation*}
\begin{aligned}
&\dot \eta_2 =
\phi'_{2,\vt} \Bigg[ \begin{pmatrix} \bar y_{\rm aux} \\ \eta_1 \end{pmatrix}, \begin{smallpmatrix} y_{\rm aux,1} - j_\vt \\ \ddt{ \bar y}_{\rm aux} \\ \eta_2 \end{smallpmatrix} \Bigg] \\
& \cdot \Big[M_\vt^{-1} \big[J_\vt^\top, G_\vt^\top, B_\vt \big]\Gamma_\vt^{-1}, V_\vt  \Big]\left(\begin{smallpmatrix} \bar y_{\rm aux} \\ \eta_1 \end{smallpmatrix} \right)\cdot  \begin{smallpmatrix} y_{\rm aux,1} - j_\vt \\ \ddt{ \bar y}_{\rm aux} \\ \eta_2 \end{smallpmatrix}   \\
&+ \phi_{2,\vt}\left(\begin{smallpmatrix} \bar y_{\rm aux} \\ \eta_1 \end{smallpmatrix} \right) M_\vt\left(\begin{smallpmatrix} \bar y_{\rm aux} \\ \eta_1 \end{smallpmatrix} \right)^{-1} f_\vt \left( \begin{smallpmatrix} \bar y_{\rm aux} \\ \eta_1 \end{smallpmatrix}, \begin{smallpmatrix} y_{\rm aux,1} - j_\vt \\ \ddt{ \bar y}_{\rm aux} \\ \eta_2 \end{smallpmatrix} \right).
\end{aligned}
\end{equation*}
}
Finally, we invoke the original constraints, which means that $y_{\rm aux,1} = 0$ and ${ \bar y}_{\rm aux} = \begin{smallpmatrix} 0 \\ y\end{smallpmatrix}$. With this, and omitting the argument $(\bar y_{\rm aux}^\top, \eta_1^\top)^\top = (0, y^\top, \eta_1^\top)^\top$ where it is obvious, the internal dynamics read
{
\begin{equation}
\begin{aligned} \label{eq:eta-dynamics}
\dot \eta_1
&= \phi'_{1, \vt}M_{\vt}^{-1} \begin{bmatrix} J_{\vt}\\ G_{\vt}\\ B_{\vt}^\top\end{bmatrix}^\top \Gamma_{\vt}^{-1} \begin{pmatrix}  - j_\vt \\ 0\\ \dot y\end{pmatrix} + \phi'_{1,\vt} V_{\vt} \eta_2,\\
\dot \eta_2 &=
\phi'_{2,\vt} \Bigg[ \begin{smallpmatrix} 0\\ y \\ \eta_1 \end{smallpmatrix}, \begin{smallpmatrix}  - j_\vt \\ 0\\ \dot y \\ \eta_2 \end{smallpmatrix} \Bigg] \cdot \left[M_\vt^{-1} \begin{bmatrix} J_{\vt}\\ G_{\vt}\\ B_{\vt}^\top\end{bmatrix}^\top, V_\vt  \right]\cdot  \begin{smallpmatrix}  - j_\vt \\ 0\\ \dot y \\ \eta_2 \end{smallpmatrix}   \\
&\quad + \phi_{2,\vt}  M_\vt^{-1} f_\vt \left( \begin{smallpmatrix} 0\\ y \\ \eta_1 \end{smallpmatrix}, \begin{smallpmatrix} - j_\vt \\ 0\\ \dot y \\ \eta_2 \end{smallpmatrix} \right).
\end{aligned}
\end{equation}
}

\begin{remark}
If no nonholonomic constraints are present we may obtain further structure for the internal dynamics in~\eqref{eq:eta-dynamics}. Consider the functions~$\phi_1$ and ~$\phi_2$ as above and recall the concept of a conservative vector field. For~$U \subseteq \R^n$ open, a vector field~$\zeta: U \to \R^n$ is said to be conservative, if it is the gradient of a~scalar field. More precisely, if
there exists a scalar field {$\chi: U \to \R$ such that~$\chi'(x) = \zeta(x)^\top$.}

Now, if every row~$\phi_{2,i}(\cdot)$ of~$\phi_2(\cdot)$ defines a conservative vector field, then the components of~$\phi_1(\cdot)$ can be chosen to be multiples of the associated scalar fields. More precisely, if there exist {$\chi_i\in \cC^1(U_1\to\R)$ such that $\chi_i'(x_1) = \phi_{2,i}(x_1)$} for all $x_1\in U_1$ and all $i=1,\ldots,n-\ell-m$, then we may choose {$\phi_1 = (\lambda_1 \chi_1,\ldots,\lambda_{n-\ell-m}\chi_{n-\ell-m})^\top$} for some $\lambda_i\in\R$, $i=1,\ldots,n-\ell-m$. This gives $\phi'_{1} = \Lambda \phi_{2}$ for $\Lambda = \diag(\lambda_1,\ldots,\lambda_{n-\ell-m})$
and, likewise, {$\phi'_{1,\vt} = \Lambda \phi_{2,\vt}$.} Hence via~\eqref{eq:K-involutive} and Lemma~\ref{Lem:phi_2} the dynamics of~$\eta_1$ in~\eqref{eq:eta-dynamics} simplify to
\[
    \dot \eta_1(t) = \Lambda \eta_2(t),
\]
where the diagonal entries of $\Lambda$ can be chosen as desired.

Hence, in this case the variables $\eta_1 = \phi_1(q)$ can be interpreted as (transformed) positions and $\eta_2 = \Lambda^{-1} \dot \eta_1$ as velocities. Therefore, the states $\eta = (\eta_1^\top, \eta_2^\top)^\top$ admit a partition which is structurally similar to that of the states $(q^\top, v^\top)^\top$ of the multi-body system~\eqref{eq:MBS}. In the presence of nonholonomic constraints this simplification is not possible, because of the different dimensions of~$\phi_1'$ and~$\phi_2$.
\end{remark}

\begin{remark}
  The computation of the Byrnes-Isidori form leading to the decoupling of the internal dynamics as in~\eqref{eq:MBS-decoupled} often requires a lot of effort. The choice of variables for the internal dynamics presented in this section offers an alternative, leading to the internal dynamics given in~\eqref{eq:eta-dynamics} directly in terms of the original system parameters; a computation of the Byrnes-Isidori form is not necessary.
\end{remark}

\section{Servo-constraints approach}\label{Sec:Servo}

We seek to find a feedforward control~$u_{\rm ff}$ regarding the control methodology described in Section~\ref{Ssec:ContrMeth} by inverting the system model. For minimum phase systems, the transformed system~(\ref{eq:MBS-decoupled}) can be used as an inverse model for feedforward control. Solving the respective first equation for the input~$u$ with $y=y_{\rm ref}$ yields the feedforward control signal~$u_{\rm ff}$, provided the internal dynamics are integrated simultaneously. However, a direct integration of the internal dynamics is not possible for non-minimum phase systems, due to unstable dynamics. Thus, stable inversion is introduced in \cite{ChenPade96,DevaChen96} for finding a bounded solution based on a boundary value problem (BVP). This yields a non-causal solution, in the sense that a system input $u_{\rm ff}(t)\neq0$ is present in a preactuation phase~$t<t_0$ before the start of the trajectory at time~$t_0$. Thus, the BVP is solved on a time interval $\left[T_0\,, T_f\right]$ with $T_0<t_0$ and $T_f>t_f$, with $t_f$ denoting the end of the desired trajectory~$y_{\rm ref}$. Provided that the equilibrium of the internal dynamics is hyperbolic, there are $n^s$ eigenvalues in the left half plane and $n^u$ eigenvalues in the right half plane, such that $n^s+n^u=2n-\bar{r}$. The boundary conditions of the BVP are chosen such that the states~$\eta$ of the internal dynamics start on the unstable manifold~$W_0^u\in\mathbb{R}^{n^u}$ of the equilibrium point~$\eta_{\rm eq,0}$ at time~$T_0$. At the end~$T_f$, the boundary conditions constrain the states~$\eta$ to lie on the stable manifold~$W_f^s\in\mathbb{R}^{n^s}$ of the equilibrium point~$\eta_{\rm eq,f}$. This is summarized as
\begin{equation}\label{eq:origbcs}
\eta(T_0) \in W_0^u \quad\wedge\quad \eta(T_f) \in W_f^s.
\end{equation}
The case of non-hyperbolic equilibria is discussed in~\cite{Deva99}. The BVP consisting of the dynamics~(\ref{eq:MBS-decoupled}) subject to the constraints~(\ref{eq:origbcs}) yields a bounded solution to the internal dynamics. It was recently demonstrated
in {the brief note}~\cite{DrueSeif20} that the boundary conditions~(\ref{eq:origbcs}) can be simplified as
\begin{equation}\label{eq:simpbcs}
L_0 \, \eta\left(T_0\right) = L_0 \, \eta_{\rm eq,0}\quad \wedge\quad
L_f \, \eta\left(T_f\right) = L_f \,  \eta_{\rm eq,f},
\end{equation}
where the {binary} matrices~$L_0\in\mathbb{R}^{n^0 \times (2n-r)}$ and $L_f\in\mathbb{R}^{n^f \times (2n-r)}$ select in total $n^0+n^f=2n-\bar{r}$ conditions to constrain~$n^0$ states of the internal dynamics to the equilibrium point at time~$T_0$ and $n^f$ states of the internal dynamics to the equilibrium point at time~$T_f$. Thus, the number of constraints equals the number of unknowns.

For higher relative degree and multi-input, multi-output systems, the symbolic derivation of the internal dynamics, and especially of the stable and unstable manifolds, becomes tedious and complex. Thus, the method of servo-constraints introduced in~\cite{Blaj92} is applied to demonstrate an alternative approach. Motivated from modelling of classical holonomic mechanical constraints, such as bearings, the equations of motion~(\ref{eq:MBS}) are appended by $m$ servo-constraints
\begin{equation}\label{eq:servo-constraints}
    h(q(t)) - y_{\rm ref}(t) = 0,
\end{equation}
which enforce the output to stay on the prescribed trajectory $y_{\rm ref}\in\cW^{r,\infty}(\R_{\ge 0}\to\R^m)$. This yields a set of DAEs
\begin{equation}\label{eq:servo-constraintsDAE}
\begin{aligned}
    \dot q(t) &= v(t),\\
    M(q(t)) \dot v(t) &= f\big(q(t),v(t)\big) + J(q(t))^\top \mu(t)\\
    &\quad  + G(q(t))^\top \lambda(t) + B(q(t))\, u(t),\\
    0 &= J(q(t)) v(t) + j(q(t)),\\
    0 &= g(q(t)),\\
    0 &=   h(q(t)) - y_{\rm ref}(t)
\end{aligned}
\end{equation}
to be solved numerically for the coordinates $q$ and $v$, the Lagrange multipliers~$\lambda$ and $\mu$ and the input $u$. We stress that the initial values $q^0, v^0, \lambda^0, \mu^0, u^0$ must be chosen so that they are consistent and the desired trajectory must be compatible with the possible motion of the system, i.e., it is required that a solution of~\eqref{eq:servo-constraintsDAE} exists on $\R_{\ge 0}$. Note that~\eqref{eq:servo-constraintsDAE} might have multiple solutions in the case that the desired output trajectory can be achieved by choosing different state trajectories. This is for example already the case for a 2-arm manipulator and most serial manipulators. However, from a control perspective, any solution that generates the desired output trajectory can be chosen. The freedom of choosing any solution can be utilized to optimize e.g.\ input energy. The chosen solution $(q,v,\lambda,\mu,u):\R_{\ge 0}\to\R^n\times\R^n\times\R^\ell\times\R^p\times\R^m$ of~\eqref{eq:servo-constraintsDAE} is used to define the feedforward control signal~$u_{\rm ff} := u$. As a consequence, if all the system parameters and initial values of the multibody system~\eqref{eq:MBS} are known exactly and the control $u= u_{\rm ff}$ is applied to it, then the closed-loop system has a solution which satisfies $y=y_{\rm ref}$. However, this result may not be true in the presence of disturbances, uncertainties and/or modelling errors. Then, a feedback controller can be added to reduce the tracking errors. This is demonstrated in the application example.

While classical mechanical constraints are enforced by reaction forces orthogonal to the tangent of the constraint manifold, the servo-constraints~\eqref{eq:servo-constraints} are enforced by the generalized input~$ B\big(q(t)\big) u(t)$ which is not necessarily perpendicular to the tangent of the constraint manifold. Different configurations are distinguished depending on 
the vector relative degree of the associated auxiliary ODE~\eqref{eq:auxODE}, see~\cite{BlajKolo04}. 
For example, if inputs and outputs are colocated as discussed in Section~\ref{Sec:colocate}, then the system input can directly actuate the output and the inverse model is in so-called ideal orthogonal realization. 
Recall that this corresponds to a relative degree two framework. If less than~$m$ components of the system input directly influence the system output, then the inverse model is in so-called mixed tangential-orthogonal or purely tangential realization.

\begin{remark}
It is interesting to note that in the colocated case, i.e., $h'(q) = B(q)^\top$ for all $q\in\R^n$, and for $y_{\rm ref} = 0$ the servo-constraints act like holonomic constraints in~\eqref{eq:servo-constraintsDAE}, where the corresponding Lagrange multipliers are given by the inputs~$u(t)$ of the system. Therefore, the system~\eqref{eq:servo-constraintsDAE} can again be interpreted as a multibody system (without any inputs or outputs) in this case.
\end{remark}

While DAEs describing mechanical systems dynamics with classical constraints are of differentiation index~3, the set of DAEs~(\ref{eq:servo-constraintsDAE}) might have a larger differentiation index. As a rule of thumb for single-input single-output systems, the differentiation index is larger by one than the relative degree of the respective system, if the internal dynamics is modeled as an ODE and not affected by a constraint~\cite{Camp95b}.

It is shown in~\cite{BrueBast13} that the original stable inversion problem, solving the boundary value problem of (\ref{eq:MBS-decoupled}) and~(\ref{eq:origbcs}) can be formulated and solved directly for the inverse model described by the DAEs~(\ref{eq:servo-constraintsDAE}). However, the derivation of the boundary conditions~\eqref{eq:origbcs} is still tedious. In order to avoid the boundary conditions, it is proposed in~\cite{BastSeif13} to reformulate the stable inversion problem as an optimization problem. Alternatively, the simplified boundary conditions~(\ref{eq:simpbcs}) can also be formulated for the inverse model described by the DAEs~(\ref{eq:servo-constraintsDAE}). They then read
\begin{equation}\label{eq:simpbcs_dae}
\begin{aligned}
&L_0 \begin{bmatrix} q(T_0)^\top & v(T_0)^\top &\lambda(T_0)^\top&\mu(T_0)^\top&u(T_0)^\top \end{bmatrix}^\top  \\
&= L_0  \begin{bmatrix} q_{\rm eq,0}^\top & v_{\rm eq,0}^\top &\lambda_{\rm eq,0}^\top&\mu_{\rm eq,0}^\top&u_{\rm eq,0}^\top \end{bmatrix}^\top\qquad \text{and}  \\
& L_f \begin{bmatrix} q(T_f)^\top & v(T_f)^\top &\lambda(T_f)^\top&\mu(T_f)^\top&u(T_f)^\top \end{bmatrix}^\top\\
&= L_f   \begin{bmatrix} q_{\rm eq,f}^\top & v_{\rm eq,f}^\top &\lambda_{\rm eq,f}^\top&\mu_{\rm eq,f}^\top&u_{\rm eq,f}^\top \end{bmatrix}^\top
\end{aligned}
\end{equation}
with the {binary selection} matrices~$L_0\in\mathbb{R}^{n^0 \times (2n+\ell+p+m)}$ and $L_f\in\mathbb{R}^{n^f \times (2n+\ell+p+m)}$. {These describe} in total $n^0+n^f=2n+\ell+p+m$ conditions, which constrain $n^0$ states to the initial equilibrium and $n^f$ states to the final equilibrium of the internal dynamics. Compared to the Byrnes-Isidori form, the resulting boundary value problem of the dynamics~(\ref{eq:servo-constraintsDAE}) subject to the constraints~(\ref{eq:simpbcs_dae}) greatly simplifies the problem setup, since less analytical derivations are required.

The boundary value problem needs to be solved numerically. Since high index DAEs are difficult to solve numerically, see e.g.~\cite{HairWann96}, different index reduction strategies can be applied. In the context of servo-constraints, index reduction by projection is proposed in~\cite{BlajKolo04}, minimal extension is proposed in~\cite{AltmBets16,BetsAltm16} and Baumgarte stabilization is applied in~\cite{BencLasz17}. Available methods for solving BVPs are for example single shooting, multiple shooting or finite differences~\cite{AschMatt95}. In the present paper, we utilize finite differences with Simpson discretization~\cite{ShamGlad03}.

{The computational effort for the solution of the BVP depends on the number of unactuated states. For systems with few unactuated states this might be done within seconds or minutes. However, the solution of the boundary value problem has to be computed offline. This is due to the fact that the BVP requires the specification of the output trajectory over the complete time horizon. If real time trajectory adaption is desired, one might have to use a different approach which allows real time capable forward time integration of the servo-constraint problem, see e.g.~\cite{OttoSeif18a,OttoSeif18}. However, this requires a differentially flat or minimum phase system. This could be achieved by changing the mechanical design of the system or by output redefinition, see e.g~\cite{MorlMeye21}. While output redefinition is quite simple, the close tracking of the original end-effector output is deteriorated.}

{Finally, it should be noted that the feedforward control of non-minimum phase systems, such as obtained by solving the BVP, is non-causal. This means that there is a so-called pre-actuation phase where some control action occurs before the trajectory tracking starts. During this pre-actuated phase the system is driven into the required initial condition for exact trajectory tracking. If this pre-actuation phase is neglected, then there is a mismatch between the actual initial conditions of the system and the ones provided by the inverse model. However, this mismatch is often small and can be compensated quickly by the feedback controller.}

\section{Funnel-based feedback controller}~\label{Sec:FunCon}

In this section we present a feedback-controller design for the multibody system~\eqref{eq:MBS}, which is based on the funnel controller for systems with arbitrary vector relative degree recently developed in~\cite{BergLe20}. Furthermore, the design invokes a recently developed methodology for {linear} non-minimum phase systems from~\cite{Berg20a}, which has been successfully applied to a nonlinear non-minimum phase robotic manipulator in~\cite{BergLanz20}~-- however, here we extend this methodology to the case of vector relative degree. The method is based on the introduction of a new output, which is differentially flat for the unstable part of the internal dynamics, cf.~\cite{FlieLevi95} for differential flatness. With respect to the new output, the overall system has a higher vector relative degree (in some of the components), but the unstable part of the internal dynamics is eliminated. This allows to apply the controller from~\cite{BergLe20} to the system~\eqref{eq:MBS} with new output and appropriate new reference signal. Since the required differentially flat output does often not exist for real-world nonlinear multibody systems (cf.~\cite{BergLanz20}), we first linearize the internal dynamics, i.e., the second equation in~\eqref{eq:MBS-decoupled}, around an equilibrium $(\eta_1^0,\eta_2^0)\in\R^{n-\ell-m}\times\R^{n-p-\ell-m}$ and $y = \dot y = \ldots = y^{(\hat r-1)} = 0$. With
\begin{align*}
    Q &= \frac{\partial \tilde q}{\partial \eta}(0,\eta_1^0,\eta_2^0)\in\R^{(2n-\bar r)\times(2n-\bar r)},\\
    P_i&=\frac{\partial \tilde p}{\partial y^{(i-1)}}(0,\eta_1^0,\eta_2^0) \in\R^{(2n-\bar r)\times m},\quad i=1,\ldots,\hat r,
\end{align*}
we obtain the linearized internal dynamics (with some abuse of notation, again using~$\eta$ as state variable)
\begin{equation}\label{eq:intdyn-lin-aux}
    \dot \eta(t) = Q\eta(t) + \sum_{i=1}^{\hat r} P_i y^{(i-1)}(t).
\end{equation}
Note that~$\tilde q$ and~$\tilde p$ may also be given in the special form~\eqref{eq:eta-dynamics}, obtained using the approach presented in Section~\ref{Sec:intdyn-coord}. In the next step, as in~\cite{BergLanz20}, we need to transform~\eqref{eq:intdyn-lin-aux} so that the derivatives of~$y$ are removed from the right hand side. To this end,
{{%
we define
\begin{equation*}
\zeta^0 := \eta, \quad
\zeta^j := \zeta^{j-1} - \sum_{i=j+1}^{\hat r} P_i y^{(i-j-1)},\quad j=1,\ldots,\hat r-1,
\end{equation*}
then a straightforward induction shows that
\begin{equation*}
\ddt \zeta^j(t) = Q \zeta^j(t) + \sum_{i=1}^{j}Q^{i-1} P_i y(t) + Q \sum_{i=j+1}^{\hat r} P_i y^{(i-j-1)}.
\end{equation*}
With some abuse of notation, again using~$\eta$ as state variable, we set $\eta := \zeta^{\hat r -1}$ and obtain
\begin{equation}\label{eq:intdyn-lin}
    \dot \eta(t) = Q\eta(t) + P y(t),
\end{equation}
where $P = \sum_{i=1}^{\hat r} Q^{i-1} P_i$. In the following we assume that $\sigma(Q) \cap \overline{\C_+} \neq \emptyset$, so that the internal dynamics are not minimum phase.
}}%

%
%
%
%
%
%
In order to derive the controller design, we first require some assumptions which are extensions of those stated in~\cite{Berg20a}. Compared to the latter, here we increase each component of the vector relative degree separately (and possibly differently), instead of uniformly increasing the strict relative degree. To this end, choose $T\in\Gl_{2n-\bar r}$ and $k\in\N$  such that
\begin{equation}\label{eq:decompQ}
    TQT^{-1} = \begin{bmatrix} \hat Q_1 & \hat Q_2 \\ 0&\tilde Q\end{bmatrix},\quad TP = \begin{bmatrix} \hat P\\ \tilde P\end{bmatrix},
\end{equation}
where $\hat Q_1\in\R^{(2n-\bar r-k)\times (2n-\bar r-k)}$, $\hat Q_2\in\R^{(2n-\bar r-k)\times k}$, $\tilde Q\in\R^{k\times k}$, $\hat P\in\R^{(2n-\bar r-k)\times m}$, $\tilde P\in\R^{k\times m}$, with $\sigma(\hat Q_1)\subseteq\C_-$ and $\sigma(\tilde Q)\subseteq \overline{\C_+}$. We assume:
\begin{enumerate}[\hspace{2pt}\textbf{(A1)}]
  \item[\textbf{(A1)}] There exist $q\in \{1,\ldots,m\}$ and $k_i\in\N$, $i=1,\ldots,q$ such that $k_1+\ldots +k_q = k$. Furthermore, there exist
  ${\rho} \in\{-1,1\}$,
  $K_1, \ldots, K_q\in\R^{1\times k}$ and $K_{q+1}, \ldots, K_m\in\R^{1\times m}$ such that, for the matrix-valued function~$\tilde S:W\to\R^{m\times m}$, $W = [0,I_{2n-p-2\ell}] \Phi(U)$, that appears in~\eqref{eq:MBS-decoupled} (which is basically the high-gain matrix of this equation), we have that
      \[
        \forall\, i=1,\ldots,q\ \forall\, j=0,\ldots,k_i-2:\ K_i \tilde Q^j \tilde P = 0
       \]
      and the matrix
      \[
        {\rho} \cdot \begin{bmatrix} K_1 \tilde Q^{k_1 - 1} \tilde P\\ \vdots \\ K_q \tilde Q^{k_q-1} \tilde P\\ K_{q+1} \\ \vdots \\ K_m\end{bmatrix} \tilde S(w)
      \]
      is positive definite for all $w\in W$.
  \item[\textbf{(A2)}] Let $y_{\rm ref}\in\cW^{\hat r,\infty}(\R_{\ge 0}\to\R^m)$ be a given reference signal and $W\in\Gl_{k}$ be such that
  \[
    W\tilde Q W^{-1} = \begin{bmatrix} Q_1&0\\ 0&Q_2\end{bmatrix},\quad W\tilde P = \begin{bmatrix} P_1\\ P_2\end{bmatrix},
  \]
  where $Q_i\in\R^{l_i\times l_i}$, $i=1,2$, with $\sigma(Q_1)\subseteq\C_+$ and $\sigma(Q_2)\subseteq {\rm i}\R$. Then the equation
  \[
    \dot z_2(t) = Q_2 z_2(t) + P_2 y_{\rm ref}(t),\quad z_2(0) = 0
  \]
  has a bounded solution $z_2:\R_{\ge 0}\to\R^{l_2}$.
\end{enumerate}
Note that in the colocated case discussed in Section~\ref{Sec:colocate}, the matrix~$S(\cdot)$ (and hence also $\tilde S(\cdot)$) is always pointwise positive definite, independent of the specific system parameters. Hence, the last condition in~\textbf{(A1)} is satisfied, if, for instance, it is possible to choose $K_1,\ldots K_m$ such that $\begin{smallbmatrix} K_1 \tilde Q^{k_1 - 1} \tilde P\\ \vdots \\ K_q \tilde Q^{k_q-1} \tilde P\\ K_{q+1} \\ \vdots \\ K_m\end{smallbmatrix} = I_m$.

We may now define the new output. Let, for~$\eta$ as in~\eqref{eq:intdyn-lin} and~$T$ as in~\eqref{eq:decompQ}, $(\eta_1^\top, \eta_2^\top)^\top = T \eta$ (not to be confused with the variables introduced in Section~\ref{Sec:intdyn-coord}). Then
\begin{equation}\label{eq:intdyn-lin-unst}
    \dot \eta_2(t) = \tilde Q \eta_2(t) + \tilde P y(t).
\end{equation}
With this we define $y_{\rm new}:\R_{\ge 0}\to\R^m$ by
\begin{equation}\label{eq:new-output}
\begin{aligned}
    y_{{\rm new}, i}(t) &:= K_i \eta_2(t),\quad &i&=1,\ldots,q,\\
    y_{{\rm new}, i}(t) &:= K_i y(t),\quad &i&=q+1,\ldots,m.
\end{aligned}
\end{equation}
Similar to~\cite{Berg20a} it can be computed that the vector relative degree (of the auxiliary ODE~\eqref{eq:auxODE}) increases when the output~$y$ is replaced by~$y_{\rm new}$. More precisely, in the context of Section~\ref{Sec:Computing-ID} we find that $r_3$ changes to $(\hat r + k_1,\ldots, \hat r +k_q, \hat r, \ldots, \hat r)\in\N^{1\times m}$.

In order to track the original reference signal~$y_{\rm ref}$ with the original output~$y$, we need to introduce a new reference signal for system~\eqref{eq:MBS} with new output~\eqref{eq:new-output}. The new reference signal is generated by the subsystem~\eqref{eq:intdyn-lin-unst} corresponding to the unstable part~$\eta_2$ when the output~$y$ is substituted by the reference~$y_{\rm ref}$, i.e.,
\begin{equation}\label{eq:new-ref}
\begin{aligned}
    \dot \eta_{2,\rm ref}(t) &= \tilde Q \eta_{2,\rm ref}(t) + \tilde P y_{\rm ref}(t),\quad &\eta_{2,\rm ref}(0)=\eta_{2,\rm ref}^0,\\
    \hat y_{{\rm ref},i}(t) &= K_i \eta_{2,\rm ref}(t),\quad &i=1,\ldots,q,\\
    \hat y_{{\rm ref},i}(t) &= K_i y_{\rm ref}(t),\quad &i=q+1,\ldots,m.
\end{aligned}
\end{equation}
Equation~\eqref{eq:new-ref} will be part of the controller design and constitutes a dynamic part of the overall controller. In order for the controller from~\cite{BergLe20} to be applicable we require that $\hat y_{{\rm ref},i}\in \cW^{\hat r+k_i,\infty}(\R_{\ge 0}\to\R^m)$ for $i=1,\ldots,q$. As shown in~\cite{Berg20a} this is the case, if, using the notation from~\textbf{(A2)},
\begin{equation}\label{eq:eta2-ref-0}
  \eta_{2,\rm ref}^0  = W^{-1} \begin{bmatrix} -I_{l_1}\\ 0_{l_2\times l_1}\end{bmatrix} \int_0^\infty e^{-Q_1 s} P_1 y_{\rm ref}(s) \ds{s}.
\end{equation}
Furthermore, if the original reference $y_{\rm ref}$ is generated by a linear exosystem (as in linear regulator problems, cf.~\cite{Fran77}) of the form
\begin{equation}\label{eq:exo}
    \dot w(t) = A_e w(t),\quad y_{\rm ref}(t) = C_e w(t),\quad w(0)=w^0,
\end{equation}
where the parameters~$A_e\in\R^{n_e\times n_e}, C_e\in\R^{m\times n_e}$ and $w^0\in\R^{n_e}$ are known, and $\sigma(A_e)\subseteq \overline{\C_-}$ so that any eigenvalue $\lambda\in\sigma(A_e)\cap {\rm i}\R$ is semisimple (note that this guarantees $y_{\rm ref}\in\cW^{\hat r,\infty}(\R_{\ge 0}\to\R^m)$), then $\eta^0_{2,\rm ref}$ can be calculated via the solution $X\in\R^{l_1\times n_e}$ of the Sylvester equation
\begin{equation*}\label{eq:Sylv}
    Q_1 X - X A_e = P_1 C_e
\end{equation*}
as
\begin{equation*}\label{eq:eta20-exo}
    \eta^0_{2,\rm ref} = W^{-1} \begin{bmatrix} -I_{l_1}\\ 0_{l_2\times l_1}\end{bmatrix} X w^0.
\end{equation*}

Now,  using the new output~\eqref{eq:new-output} and the new reference signal~\eqref{eq:new-ref}, the auxiliary tracking error is defined by
\begin{align}
\label{eq:e-new}
e_0(t) = y_{\rm new}(t) - \hat{y}_{\rm ref}(t).
\end{align}
The application of the funnel control law from~\cite{BergLe20} requires the derivatives of~$y_{\rm new}$, and for implementation we need to express them in terms of the original variables~$q$ and~$v$ of the multibody system~\eqref{eq:MBS}. We use the linearization~\eqref{eq:intdyn-lin-unst} to obtain these derivatives. We replace~$\eta_2$ in~$y_{\rm new}$ and its derivatives by the original coordinates from~\eqref{eq:MBS}. To this end, we pretend that the coordinates of the linearization of the internal dynamics in~\eqref{eq:intdyn-lin-aux} coincide with the original coordinates in the second equation of~\eqref{eq:MBS-decoupled}, which can be obtained via
\[
    \Phi(q,v) = (0,y,\dot y, \ldots, y^{(\hat r-1)}, \eta).
\]
Assume that the transformation from~\eqref{eq:intdyn-lin-aux} to~\eqref{eq:intdyn-lin} is given by
\[
    (y,\dot y, \ldots, y^{(\hat r-1)}, \eta) \mapsto F_0 \eta + \sum_{i=1}^{\hat r} F_i y^{(i-1)},
\]
then we consider
\[
\Psi: U \to \R^{k},\ (q,v) \mapsto [0, I_{k}] T [0, F_1, \ldots, F_{\hat r}, F_0] \Phi(q,v).
\]
Now we replace $\eta_2(t) = \Psi(q(t),v(t))$ in~$y_{{\rm new}, i}(t) = K_i \eta_2(t)$ and its derivatives for $i=1,\ldots,q$. To this end, observe that, similar as in~\cite{Berg20a}, we may calculate that
\[
    \begin{pmatrix}
      y_{{\rm new}, i}(t) \\ \dot y_{{\rm new}, i}(t) \\ \vdots \\ y_{{\rm new}, i}^{(k_i-1)}(t) \\ y_{{\rm new}, i}^{(k_i)}(t)
    \end{pmatrix} = \begin{bmatrix} K_i\\ K_i \tilde Q \\ \vdots \\ K_i \tilde Q^{k_i-1} \\ K_i \tilde Q^{k_i}\end{bmatrix} \eta_2(t) + \begin{bmatrix} 0\\ \vdots\\ 0\\ K_i \tilde Q^{k_i-1} \tilde P\end{bmatrix} y(t).
\]
Hence, we obtain the replacement rule
\begin{align}
\label{eq:derivatives-ynew-approx}
y_{{\rm new}, i}(t) &= K_i \Psi(q(t),v(t)),\nonumber \\
\dot y_{{\rm new}, i}(t) &=  K_i \tilde Q \Psi(q(t),v(t)),\nonumber \\
 &\ \vdots  \nonumber\\
 y_{{\rm new}, i}^{(k_i-1)}(t) &= K_i \tilde Q^{k_i-1} \Psi(q(t),v(t)),\nonumber\\
y_{{\rm new}, i}^{(k_i)}(t) &= K_i \tilde Q^{k_i} \Psi(q(t),v(t)) + K_i \tilde Q^{k_i-1} \tilde P y(t),\nonumber\\
y_{{\rm new}, i}^{(k_i+1)}(t) &= K_i \tilde Q^{k_i+1} \Psi(q(t),v(t)) + K_i\tilde Q^{k_i} \tilde P y(t) \nonumber \\
&\quad + K_i \tilde Q^{k_i-1} \tilde P \dot y(t),\nonumber\\
&\ \vdots \nonumber\\
y_{{\rm new}, i}^{(\hat r+k_i-1)}(t) &= K_i \tilde Q^{\hat r + k_i-1} \Psi(q(t),v(t)) + \sum_{j=0}^{\hat r-2} K_i\tilde Q^{k_i+j} \tilde P y^{(j)}(t) \nonumber \\
&\quad +K_i\tilde Q^{k_i-1} \tilde P y^{(\hat r -1)}(t),
\end{align}
for $i=1,\ldots,q$, where the derivatives~$y^{(j)}(t)$ are to be expressed via Lie derivatives as in~\eqref{eq:Phi-aux} and read
\begin{equation} \label{eq:y-derivatives-via-Lie}
y^{(j)}(t) = \big(L_{F}^{j} O_3\big)(q(t),v(t)), \quad 0 \le j \le \hat r -1.
\end{equation}
Similarly, we obtain the replacement rule
\[
    y_{{\rm new}, i}^{(j)}(t) = K_i \big(L_{F}^{j} O_3\big)(q(t),v(t))
\]
for $i=q+1,\ldots,m$ and $j=0,\ldots,\hat r-1$.

{If the transformation $\Psi:U\to\R^k$ is not available or hard to compute, then $\eta_2$ may instead be computed as the solution of the linear differential equation~\eqref{eq:intdyn-lin-unst} with $y(t) = h(q(t),v(t))$, which is inserted at the respective places above.}

With this, the application of the controller from~\cite{BergLe20} to a multibody system~\eqref{eq:MBS} with new output~\eqref{eq:new-output} and reference signal as in~\eqref{eq:new-ref} leads to the following overall control law:
\fbox{\parbox{0.47\textwidth}{%
\begin{align}
\label{eq:Funnel-Control}
\dot \eta_{2,\rm ref}(t) &= \tilde Q \eta_{2,\rm ref}(t) + \tilde P y_{\rm ref}(t),\quad \eta_{2,\rm ref}(0)=\eta_{2,\rm ref}^0,\nonumber \\
i=1,&\ldots,q: \nonumber\\
\hat y_{{\rm ref},i}(t) &= K_i \eta_{2,\rm ref}(t),\nonumber\\
    y_{{\rm new},i}(t) &= K_i \Psi(q(t),v(t)),\nonumber \\
e_{i,0}(t) &= y_{{\rm new},i}(t) - \hat{y}_{{\rm ref},i}(t), \nonumber \\
e_{i,1}(t) &= \dot e_{i,0}(t) + k_{i,0}(t)e_{i,0}(t)\quad \text{using~\eqref{eq:derivatives-ynew-approx}+\eqref{eq:y-derivatives-via-Lie}}, \nonumber \\
e_{i,2}(t) &= \dot e_{i,1}(t) + k_{i,1}(t)e_{i,1}(t)\quad \text{using~\eqref{eq:derivatives-ynew-approx}+\eqref{eq:y-derivatives-via-Lie}}, \nonumber \\
&\, \ \vdots \nonumber\\
e_{i,\hat r+k_i-1}(t)&=\dot{e}_{i,\hat r+k_i-2}(t)+k_{i,\hat r+k_i-2}(t)\,e_{i,\hat r+k_i-2}(t)\nonumber\\
&\quad \text{using~\eqref{eq:derivatives-ynew-approx}+\eqref{eq:y-derivatives-via-Lie}},\nonumber\\
k_{i,j}(t)&=\frac{\kappa_{i,j}}{1-\varphi_{i,j}(t)^2e_{i,j}(t)^2},\quad j=0,\ldots,\hat r+k_i-2,\nonumber \\
i=q+1,&\ldots,m: \nonumber\\
e_{i,0}(t) &= K_i \big(y(t) - y_{{\rm ref}}(t)\big), \nonumber \\
e_{i,1}(t) &= \dot e_{i,0}(t) + k_{i,0}(t)e_{i,0}(t)\quad \text{using~\eqref{eq:y-derivatives-via-Lie}}, \nonumber \\
e_{i,2}(t) &= \dot e_{i,1}(t) + k_{i,1}(t)e_{i,1}(t)\quad \text{using~\eqref{eq:y-derivatives-via-Lie}}, \nonumber \\
&\, \ \vdots \nonumber\\
e_{i,\hat r-1}(t)&=\dot{e}_{i,\hat r-2}(t)+k_{i,\hat r-2}(t)\,e_{i,\hat r-2}(t)\quad \text{using~\eqref{eq:y-derivatives-via-Lie}},\nonumber\\
k_{i,j}(t)&=\frac{\kappa_{i,j}}{1-\varphi_{i,j}(t)^2e_{i,j}(t)^2},\quad j=0,\ldots,\hat r-2,\nonumber \\
\bar e(t) &= \big(e_{1,\hat r+k_1-1},\ldots,e_{q,\hat r+k_q-1}, \nonumber \\
&\qquad e_{q+1,\hat r-1},\ldots, e_{m,\hat r-1}\big)^\top,\nonumber \\
\bar k(t) &= \frac{\bar \kappa}{1-\varphi(t)^2\|\bar e(t)\|^2},\nonumber \\
u(t) &= - {\rho} \,\bar k(t)\,\bar e(t),
\end{align}
}}
with~${\eta}_{2,\rm ref}^0$ as in~\eqref{eq:eta2-ref-0}, $ \rho$ as in~\textbf{(A1)}, tuning parameters $\bar \kappa, \kappa_{i,j} > 0$, reference signal $y_{\rm ref}\in \cW^{\hat r,\infty}(\R_{\ge 0}\to\R)$ and funnel functions $\varphi_{i,j} \in\Phi_{\hat r + k_i-j}$ for $i=1,\ldots,q$, $j=0,\ldots,\hat r + k_i -2$, $\varphi_{i,j} \in\Phi_{\hat r -j}$ for $i=q+1,\ldots,m$, $j=0,\ldots,\hat r  -2$ and $\varphi\in \Phi_1$.

We note that in~\eqref{eq:Funnel-Control} the variables $e_{i,j}$ are only short-hand notations and they can be expressed in terms of~$y_{\rm new}, \hat y_{\rm ref}, \varphi_{i,j}$ and the derivatives of these. Therefore, the expressions in~\eqref{eq:derivatives-ynew-approx} must be used for the replacement of the derivatives of~$y_{\rm new}$ in terms of the original coordinates~$q$ and~$v$. Further note that the controller \eqref{eq:Funnel-Control} implicitly assumes that the physical dimensions of all input and output variables coincide. For systems in which these variables have different physical dimensions a straightforward renormalization approach as presented in~\cite{BergOtto19} can be used.
{{%
\begin{remark}
Since the controller design developed in this section involves numerous transformations, a summarizing remark is warranted. For a given multibody system~\eqref{eq:MBS} with unstable internal dynamics we perform the following steps.
\begin{enumerate}
\item Decouple the internal dynamics as in~\eqref{eq:MBS-decoupled} using the method of Section~\ref{Sec:Computing-ID} or, alternatively, that of Section~\ref{Sec:intdyn-coord}.
\item Linearize the internal dynamics around an equilibrium to obtain~\eqref{eq:intdyn-lin-aux}.
\item Transform the linearization such that the derivatives of~$y$ are removed from the right hand side and~\eqref{eq:intdyn-lin} is obtained.
\item Decompose the internal dynamics via~\eqref{eq:decompQ} such that the unstable part of the linearized internal dynamics is decoupled and given by~\eqref{eq:intdyn-lin-unst}.
\item Define the new output $y_{\rm new}$ as in~\eqref{eq:new-output} and express it as a function of the original variables~$q$ and~$v$ of~\eqref{eq:MBS}.
\item Calculate the new reference signal $\hat y_{\rm ref}$ via~\eqref{eq:new-ref}.
\item Calculate the required derivatives of $y_{\rm new}$ and $\hat y_{\rm ref}$ using~\eqref{eq:derivatives-ynew-approx},~\eqref{eq:y-derivatives-via-Lie} and~\eqref{eq:new-ref}.
\item Calculate the control input~$u$ via~\eqref{eq:Funnel-Control}.
\end{enumerate}
\end{remark}
}}
%
%

\section{A robotic manipulator with kinematic loop}~\label{Sec:Robot}

In this section, we illustrate our findings by a robotic manipulator which is described by a DAE that cannot be reformulated as an ODE and, at the same time, has unstable internal dynamics. The model is shown in Figure~\ref{fig:flexarm}. This example is motivated by a similar flexible manipulator in~\cite{MorlMeye18}. It consists of three bodies with mass~$m_i$, length~$L_i$ and moments of inertia~$I_i$. The first body is fixed to a translationally moving actuator on one end. The other end is connected to the end of body~2. The center of gravity of body~2 is mounted on a translationally moving actuator. Body~3 is attached to the end of body~2 by a linear spring-damper combination with coefficients~$c$ and $D$. Due to the described configuration, there exists a kinematic loop in the model and the resulting equations of motion in DAE-form cannot be easily reformulated as an ODE.

With the variables $q = (s_1, s_2, \alpha, \beta, \gamma)^\top$, $\lambda = (\lambda_1, \lambda_2)^\top$ and the system input $u=(u_1, u_2)^\top$ the equations of motion are given by
\begin{align*}
    \dot q(t) &= v(t),\\
  M(q(t)) \dot v(t) &= f\big(q(t),v(t)\big) + G(q(t))^\top \lambda(t) + B u(t),\\
  0&= g(q(t)),
\end{align*}
which are of the form~\eqref{eq:MBS} with~$f, g, M, G, B$ defined in~\eqref{eq:robotMgB} in a brief way. The model parameters are listed in Table~\ref{tab:flexarm-param}. Note that a homogeneous mass distribution is assumed for bodies 1 and 2, while the center of gravity of body 3 can be varied by~$X_3$. The parameters of the reference model differ from the parameters for the forward simulation of the actual model in order to show robustness of the presented approach. For this example, the mass~$m_3$ of the third arm is chosen 20\% larger compared to the reference model. This might be the case if the robot is hoisting an unknown mass with its third arm.

The reference model is inverted according to Section~\ref{Sec:Servo} to obtain the feedforward control input~$u_{\rm ff}$. Using identical parameters for both models would yield exact tracking of the feedforward controller. Therefore, the parameters for the simulated model differ from the ones of the inverted model. The funnel controller introduced in Section~\ref{Sec:FunCon} decreases the tracking errors due to parameter missmatch of the pure feedforward control. However, note that the funnel controller~\eqref{eq:Funnel-Control} also requires some system parameters, for instance to compute the matrices~$\tilde Q$ and~$\tilde P$, and hence relies on the reference model. Nevertheless, we will demonstrate that it is able to compensate the tracking error induced by feedforward control.

\begin{figure*}[h!tb]
\begin{equation}\label{eq:robotMgB}
\begin{aligned}
M&=\begin{bmatrix}
 m_{1} & 0 & \frac{L_{1}m_{1}\cos\left(\alpha\right)}{2} & 0 & 0\\
0 & m_{2}+m_{3} & 0 & -X_{3}m_{3}\sin\left(\beta+\gamma\right)-\frac{L_{2}m_{3}\sin\left(\beta\right)}{2} & -X_{3}m_{3}\sin\left(\beta+\gamma\right)\\
\frac{L_{1}m_{1}\cos\left(\alpha\right)}{2} & 0 & \frac{m_{1}L_{1}^2}{4}+I_1 & 0 & 0\\
0 & -X_{3}m_{3}\sin\left(\beta+\gamma\right)-\frac{L_{2}m_{3}\sin\left(\beta\right)}{2} & 0 & \frac{m_{3}L_{2}^2}{4}+m_{3}\,\cos\left(\gamma\right)L_{2}X_{3}+m_{3}X_{3}^2+I_{2}+I_{3} & m_{3}X_{3}^2+\frac{L_{2}m_{3}\,\cos\left(\gamma\right)X_{3}}{2}+I_{3}\\
0 & -X_{3}m_{3}\,\sin\left(\beta+\gamma\right) & 0 & m_{3}X_{3}^2+\frac{L_{2}m_{3}\cos\left(\gamma\right)X_{3}}{2}+I_{3} & m_{3}X_{3}^2+I_{3}
\end{bmatrix},
\\
f &= \begin{pmatrix}
 \tfrac12 {L_{1}\,{\dot{\alpha}}^2\,m_{1}\,\sin\left(\alpha\right)}\\
 \tfrac12 {m_{3}\,\left(2\,X_{3}\,{\dot{\beta}}^2\,\cos\left(\beta+\gamma\right)+2\,X_{3}\,{\dot{\gamma}}^2\,\cos\left(\beta+\gamma\right) +L_{2}\,{\dot{\beta}}^2\,\cos\left(\beta\right)+4\,X_{3}\,\dot{\beta}\,\dot{\gamma}\,\cos\left(\beta+\gamma\right)\right)}\\
 0\\
 \tfrac12 {L_{2}\,X_{3}\,\dot{\gamma}\,m_{3}\,\sin\left(\gamma\right)\,\left(2\,\dot{\beta}+\dot{\gamma}\right)}\\
 -\tfrac12 {L_{2}\,X_{3}\,m_{3}\,\sin\left(\gamma\right)\,{\dot{\beta}}^2}-D\,\dot{\gamma}-c\,\gamma
 \end{pmatrix},\qquad B = \begin{bmatrix}
     1  &   0 \\
     0  &  1 \\
     0  &  0 \\
     0  &   0 \\
     0  &   0
\end{bmatrix}, \\
 \\
g &= \begin{bmatrix}
 L_1\,\cos(\alpha) - s_2 - d + \frac{1}{2}(L_2\,\cos(\beta)) \\
     s_1 + L_1\,\sin(\alpha) - \frac{1}{2}(L_2\,\sin(\beta))
\end{bmatrix},\qquad
G = \begin{bmatrix}
0 & -1 & -L_1\,\sin(\alpha) & -\frac{1}{2}(L_2\,\sin(\beta)) & 0 \\
1 &  0 &  L_1 \,\cos(\alpha) &  -\frac{1}{2}(L_2\cos(\beta)) & 0
\end{bmatrix}. \\
\end{aligned}
\end{equation}
\end{figure*}

\begin{figure}[tbp]
\begin{center}
\includegraphics[width=7cm]{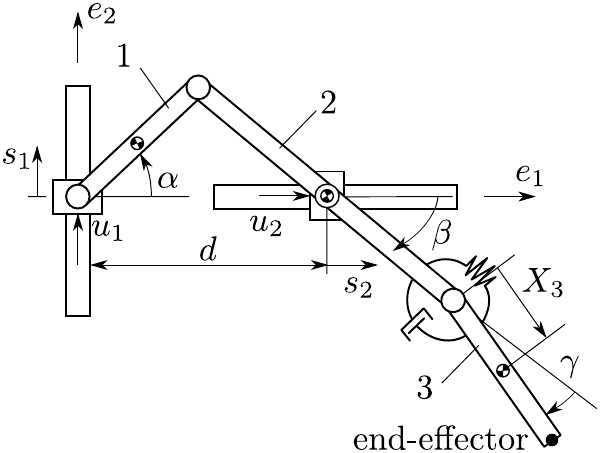}
\caption{Model of the flexible robot.}
\label{fig:flexarm}
\end{center}
\end{figure}

\begin{table}[tbp]
\begin{center}
\begin{tabular}{ccc}
\hline
& Inv. model & Sim. model \\
\hline
$m_1$ (kg) & 3.4 & 3.4 \\
$L_1$ (m) & 0.5 & 0.5 \\
$I_1$ (kgm$^2$) & 0.071 & 0.071\\
$m_2$ (kg) & 6.8 & 6.8 \\
$L_2$ (m) & 1 & 1 \\
$I_2$ (kgm$^2$) & 0.567 & 0.567 \\
$m_3$ (kg) & 3.4 & \textbf{4.1} \\
$L_3$ (m) & 0.5 & 0.5 \\
$X_3$ (m) & 0.25 & 0.25\\
$I_3$ (kgm$^2$) & 0.071 &  \textbf{ 0.085}\\
$c$ (\si{\newton\per{\radian}}) & 50 & 50 \\
$D$ (\si{\newton\second\per{\radian}}) & 0.25 & 0.25 \\
$d$ (m) & 0.8 & 0.8 \\
\hline
\end{tabular}
\caption{Parameters of the inverted reference model and the simulated model.}
\label{tab:flexarm-param}
\end{center}
\end{table}


The output is chosen as
\[
    y(t) = \begin{pmatrix} y_1(t)\\ y_2(t)\end{pmatrix} = h(q(t)) = \begin{pmatrix} s_2(t)\\ \beta(t) + \frac{2L_3}{L_2+ 2L_3} \gamma(t)\end{pmatrix},
\]
where the first component is the horizontal position of the second body and the second component corresponds to an auxiliary angle describing the end-effector position for small angles~$\gamma$. Both components may be used to approximate the end-effector position by
\[
    r_{\rm app}(t) = \begin{pmatrix} d + y_1(t) + (\frac{1}{2}L_2 + L_3) \cos (y_2(t) )\\ -(\frac{1}{2}L_2 + L_3) \sin (y_2(t))\end{pmatrix} \,.
\]
{Note that in practice the output $y(t)$ can be easily measured. The first component of~$y(t)$ is directly obtained by measuring the position $s_2(t)$ of the second translational actuator. The second component of the output is computed from the measured joint angles $\beta(t)$ and $\gamma(t)$ and the known lengths $L_2$ and $L_3$, see Figure~\ref{fig:flexarm}. A detailed explanation and derivation of such an output to approximate the end-effector position can be found in~\cite{Seif14}.}

To demonstrate the tracking capability of the proposed controller design we choose the reference trajectory~$r_{\rm app,ref}$ in end-effector coordinates as the path
\begin{align*}
r_{\rm app,ref}(t) = \begin{cases}
r_{\rm app,0}  , \quad &t < t_0 \\
r_{\rm app,0} + {{r(t)}} \left( r_{\rm app,f}-r_{\rm app,0} \right) , \quad &t_0\leq t\leq t_f \\
r_{\rm app,f}  , \quad &t> t_f
\end{cases}
 \end{align*}
parametrized by~${{r(t)}}$. The timing law of the scalar parameter~${r(t)}$ is chosen as the polynomial
\begin{align*}
{r(t)} = \,&70\left(\dfrac{t}{t_f}\right)^9 -315\left(\dfrac{t}{t_f}\right)^8 + 540\left(\dfrac{t}{t_f}\right)^7 \\ &-420\left(\dfrac{t}{t_f}\right)^6 + 126\left(\dfrac{t}{t_f}\right)^5 \,
\end{align*}
with initial time $t_0=\SI{0}{\second}$ and final transition time~$t_f=\SI{1}{\second}$. The initial position is chosen as $r_{\rm app,0}=\begin{bmatrix}1.6 &-0.6\end{bmatrix}^\top\,$m and the final position is $r_{\rm app,f}=\begin{bmatrix}0.9 &-0.9\end{bmatrix}^\top\,$m. The reference trajectory~$r_{\rm app,ref}$ is then transformed to the system output trajectory~$y_{\rm ref}$.

In the following we aim to calculate the relative degree in a preferably large open set around the equilibrium point
\[
    (s_1^0, s_2^0, \alpha^0, \beta^0, \gamma^0) = \left(0, 0, \frac{\pi}{4}, \frac{\pi}{4}, 0\right),\quad \lambda^0 = 0,\quad u^0 = 0.
\]
To this end, we choose (and this will be justified later)
\begin{align*}
    U_q &:= \R^2 \times  \left(0, \frac{\pi}{2}\right)^2 \times \setdef{\gamma\in [-\pi,\pi) }{ \cos \gamma > \tfrac23},\\
     U &:= U_q \times \R^5,
\end{align*}
which accordingly restricts the operating range of the system.  Henceforth, we identify~$M$, $G$, $g$ and~$h$ with their restrictions to~$U_q$ and~$f$ with its restriction to~$U$.
It is straightforward to check that the conditions in~\eqref{eq:assumptions} are satisfied on~$U_q$. However, the inputs and outputs are not colocated here, i.e.,~\eqref{eq:colocate} does \emph{not} hold. It is then easy to see that for all $x = (q^\top, v^\top)^\top\in U$ we have $L_{K(x)} O_3(x) = 0$ and
\begin{align*}
    \Gamma(x) &= \begin{bmatrix}  \Gamma_2(x)\\ \Gamma_3(x)\end{bmatrix} = \begin{bmatrix}  (L_{K} L_{F} O_2)(x)\\  (L_{K} L_{F}  O_3)(x)\end{bmatrix}\\
     &= \begin{bmatrix} G(q) M(q)^{-1} G(q)^\top & G(q) M(q)^{-1} B\\ h'(q) M(q)^{-1} G(q)^\top & h'(q) M(q)^{-1} B\end{bmatrix},
\end{align*}
where
\[
    h'(q) = \begin{bmatrix} 0&1&0&0&0\\ 0&0&0&1&\frac{2L_3}{L_2+ 2L_3}\end{bmatrix},
\]
and a Matlab calculation yields that
\[
    \det \Gamma(x) =  -L_1^2 L_2^2 \sin(\alpha) \sin(\alpha+\beta) \frac{ I_3 + m_3 X_3^3 - m_3 L_3 X_3 \cos\gamma}{(2L_2 + 4L_3) \det(M(q))}.
\]
Let us assume in the following that the third body has a homogeneous mass distribution, and hence
\begin{equation}\label{eq:body3-hom}
    X_3 = \tfrac12 L_3,\quad I_3 = \tfrac{1}{12} m_3 L_3^2.
\end{equation}
Now, for $x = (q^\top, v^\top)^\top\in U$ we find that $\sin(\alpha)>0$, $\sin(\alpha+\beta)>0$ and $\cos \gamma > \tfrac23$, where the latter gives
\begin{align*}
    I_3 + m_3 X_3^3 - m_3 L_3 X_3 \cos\gamma &\stackrel{\eqref{eq:body3-hom}}{=} \tfrac{1}{12} m_3 L_3^2 + \tfrac{1}{4} m_3 L_3^2 - \tfrac{1}{2} m_3 L_3^2 \cos\gamma \\
    &  <   m_3 L_3^2 \left(\tfrac{1}{12} +  \tfrac{1}{4} - \tfrac{1}{2} \cdot\tfrac23\right) = 0,
\end{align*}
thus $\det \Gamma(x) < 0$ and hence~$\Gamma(x)$ is invertible and the vector relative degree is well defined on~$U$ with $\hat r = 2$. As a consequence we find that $\bar r = 4 + 4 = 8 < 10 = 2n$, so the system has nontrivial internal dynamics.

For later use set $\kappa := m_3 X_3^2 + I_3  = \frac{m_3 L_3^2}{3}$. In order to decouple the internal dynamics we invoke the transformation~$\Phi$ with
\[
   \xi = \begin{pmatrix} \phi_1(x)\\ \vdots\\ \phi_8(x)\end{pmatrix} = \begin{pmatrix} g_1(q)\\ g_2(q)\\ G_1(q) v\\ G_2(q) v\\ h_1(q)\\ h_2(q)\\ h_1'(q) v \\ h_2'(q) v\end{pmatrix}
\]
and following the ansatz~\eqref{eq:eta-structure} in Section~\ref{Sec:intdyn-coord} we set
\[
    \eta = \begin{pmatrix} \eta_1\\ \eta_2\end{pmatrix} = \begin{pmatrix} \phi_9(x)\\ \phi_{10}(x) \end{pmatrix},
\]
where~$\phi_{10}(x) = \tilde \phi_{10}(q) v$. As in the proof of Lemma~\ref{Lem:existence-phi1/2} we choose
\begin{align*}
\tilde \phi_{10}: U_q &\to \R^{1 \times 5} \\
q &\mapsto [0, -\frac{m_3 L_3}{2} \sin(\beta + \gamma), 0, \kappa + \frac{m_3 L_2 L_3}{4} \cos(\gamma), \kappa],
\end{align*}
and observe that~$\tilde \phi_{10}(q) = [0,0,0,0,1]M(q)$. We check that
\begin{align*}
    \begin{bmatrix}  \frac{\partial}{\partial q} \big( \tilde \phi_{10}(q)v \big) & \tilde \phi_{10}(q)   \end{bmatrix} K(q) &= \tilde \phi_{10}(q) M(q)^{-1} [G(q)^\top, B] \\
    &= [0,0,0,0,1] [G(q)^\top, B] = 0.
\end{align*}
Further, we choose~$\phi_9: U_q \to \R,\ q \mapsto \gamma$
and thus
\[
    \phi'_9(q) K(q) = [0,0,0,0,1,0,0,0,0,0] K(q) = 0.
\]
Furthermore, the Jacobian of~$\Phi$ is given by~\eqref{eq:dotPhi}.
\begin{figure*}[h!tb]
\begin{equation}\label{eq:dotPhi}
    \Phi'(x) = \begin{bmatrix} 0 & -1 & -L_1\sin\alpha & -\frac{L_2}{2}\sin\beta &0&0&0&0&0&0\\
    1 & 0& L_1\cos\alpha & -\frac{L_2}{2}\cos\beta &0&0&0&0&0&0 \\
    0&0& -L_1\cos(\alpha)\dot \alpha & -\frac{L_2}{2} \cos(\beta)\dot\beta &0&0&-1&-L_1\sin\alpha&-\frac{L_2}{2}\sin\beta&0\\
    0&0& -L_1\sin(\alpha)\dot\alpha&\frac{L_2}{2}\sin(\beta)\dot\beta&0&1&0&L_1\cos\alpha&-\frac{L_2}{2}\cos\beta&0\\
    0&1&0&0&0&0&0&0&0&0\\
    0&0&0&1&\frac{2L_3}{L_2+ 2L_3}&0&0&0&0&0\\
    0&0&0&0&0&0&1&0&0&0\\
    0&0&0&0&0&0&0&0&1&\frac{2L_3}{L_2+ 2L_3}\\
        0&0&0&0&1&0&0&0&0&0\\
    0&0&0&-\frac{L_3 m_3}{2}\cos(\beta+\gamma)\dot s_2&-\frac{L_3 m_3}{2}\cos(\beta+\gamma)\dot s_2 - \frac{m_3 L_2 L_3}{4} \sin(\gamma)\dot\beta & 0 & -\frac{L_3 m_3}{2} \sin(\beta+\gamma) & 0 & \kappa  + \tfrac{m_3 L_2 L_3 }{4} \cos \gamma & \kappa
    \end{bmatrix}
\end{equation}
\end{figure*}
Invertibility of the matrix in~\eqref{eq:dotPhi} is equivalent to invertibility of the submatrix
\[
    \begin{bmatrix} -L_1\sin\alpha &0&0&0\\
    -L_1\cos(\alpha)\dot \alpha&-L_1\sin\alpha&-\frac{L_2}{2}\sin\beta&0\\
    0&0&1& \frac{2L_3}{L_2+ 2L_3}\\
    0&0 & \kappa  + \tfrac{m_3 L_2 L_3}{4}  \cos \gamma & \kappa\end{bmatrix}.
\]
Since $x\in U$ implies that $\sin \alpha >0$ this matrix is invertible if, and only if, the determinant of the lower right $2\times 2$ submatrix is nonzero, which is given by
\begin{align*}
    &\det \begin{bmatrix} 1&\frac{2L_3}{L_2+ 2L_3}\\
    \kappa + \tfrac{m_3  L_2 L_3}{4} \cos \gamma & \kappa \end{bmatrix}\\
    &= \kappa  - \frac{2L_3}{L_2+ 2L_3} \big(\kappa  + \tfrac{m_3 L_2 L_3}{4}  \cos \gamma\big) \\
    &= \frac{2L_3}{L_2+ 2L_3} \left( \frac{\kappa L_2}{2L_3}   - \frac{m_3  L_2 L_3 }{4}\cos \gamma \right) \\
    &=     \frac{L_2}{L_2+ 2L_3} \left( \tfrac{m_3 L_3^2}{3} - \tfrac{m_3  L_3^2}{2} \cos \gamma \right)
     <  \frac{ m_3 L_2 L_3^2}{L_2+ 2L_3} \left( \frac13  -  \frac13 \right) = 0,
\end{align*}
since $\cos \gamma > \tfrac23$. Therefore, $\Phi'(x)$ is invertible everywhere for all $x\in U$. In the next step we aim to obtain the internal dynamics, i.e., the second equation in the decoupled system~\eqref{eq:MBS-decoupled}. First we calculate that
\begin{align*}
  \dot \eta_1(t) &= \dot \gamma(t) ,\\
  \dot \eta_2(t) &= \big(\ddt \tilde \phi_{10}(q(t))\big) \dot q(t) + \tilde \phi_{10}(q(t)) \ddot q(t)\\
  &= -\tfrac{L_3 m_3}{2} \cos(\beta(t) + \gamma(t)) (\dot \beta(t) + \dot \gamma(t)) \dot s_2(t)\\
  &\quad -\! \tfrac{m_3}{4} L_2 L_3 \sin(\gamma(t)) \dot \beta(t) (\dot \beta(t) \!+\! \dot \gamma(t)) \!-\! D \dot \gamma(t) \!-\! c \gamma(t).
\end{align*}
In order to resolve the right hand sides in the above equation, we calculate the inverse of the diffeomorphism~$\Phi$. First observe that the new coordinates admit the representation in~\eqref{eq:robot-coord}.
\begin{figure*}[h!tb]
\begin{equation}\label{eq:robot-coord}
\begin{aligned}
  \begin{pmatrix} \phi_1\\ \phi_2\\ \phi_5\\ \phi_6\\ \phi_9\end{pmatrix}
  &=
  \begin{bmatrix}
  0&-1&0&0&0\\ 1&0&0&0&0\\ 0&1&0&0&0\\ 0&0&0&1&\frac{2L_3}{L_2+ 2L_3}\\ 0&0&0&0&1
  \end{bmatrix}
  \begin{pmatrix} s_1\\ s_2\\ \alpha\\ \beta\\\gamma\end{pmatrix}
  + \begin{pmatrix}
  L_1\cos\alpha + \tfrac{L_2}{2} \cos\beta - d\\ L_1\sin\alpha - \tfrac{L_2}{2} \sin\beta \\ 0\\ 0\\ 0
  \end{pmatrix},\\
  \begin{pmatrix} \phi_3\\ \phi_4\\ \phi_7\\ \phi_8\\ \phi_{10}\end{pmatrix}
  &=\begin{bmatrix}
  0&-1&-L_1\sin\alpha&-\tfrac{L_2}{2}\sin\beta&0\\ 1&0&L_1\cos\alpha&-\tfrac{L_2}{2}\cos\beta&0\\ 0&1&0&0&0\\ 0&0&0&1&\frac{2L_3}{L_2+ 2L_3}\\ 0&-\tfrac{L_3 m_3}{2} \sin(\beta+\gamma) &0& \kappa + \tfrac{m_3}{4} L_2 L_3 \cos \gamma & \kappa \end{bmatrix} \begin{pmatrix} \dot s_1\\ \dot s_2\\ \dot\alpha\\ \dot\beta\\\dot\gamma
  \end{pmatrix}.
\end{aligned}
\end{equation}
\end{figure*}
Upon solving, and invoking the original constraints $\phi_1 = \ldots = \phi_4 = 0$, with $\delta := \frac{2L_3}{L_2+ 2L_3}$ we obtain for the required coordinates $s_2, \beta, \gamma$ and~$\dot s_2,\dot \beta$ and $\dot\gamma$ that
\begin{align*}
   s_2 &= \phi_5,\quad
   \gamma = \phi_{9},\quad
   \beta =  \phi_6 - \delta \phi_{9},\quad
   \dot s_2 = \phi_7,\\
 \dot \beta &=
 \tfrac{2(L_2 + 2L_3)}{\kappa L_2 \big(2 - 3\cos(\phi_9 ) \big)}
 \bigg( \kappa \phi_8 - \delta \phi_{10} \\
 &\quad - \tfrac{\delta m_3 L_3 }{2} \sin\big(\phi_6  + (1- \delta) \phi_9\big) \phi_7 \bigg)  \\
\dot \gamma &= \tfrac{2(L_2 + 2L_3)}{\kappa L_2 \big(2 - 3\cos(\phi_9)\big)}
\bigg( \phi_{10} - \big(\kappa + \tfrac{m_3 L_2 L_3}{4} \cos(\phi_9) \big) \phi_8  \\
&\quad + \tfrac{m_3 L_3 }{2} \sin\big(\phi_6 + (1- \delta)\phi_9 \big) \phi_7 \bigg).
%
%
%
%
\end{align*}
%
With this and $\phi_5 = y_1$, $\phi_6 = y_2$, $\phi_7 = \dot y_1$, $\phi_8 = \dot y_2$, $\phi_9 = \eta_1$ and $\phi_{10} = \eta_2$ the internal dynamics are given by~\eqref{eq:ID-robot}.

\begin{figure*}[h!tb]
\begin{equation} \label{eq:ID-robot}
\begin{aligned}
\dot \eta_1 &=
\frac{2(L_2 + 2L_3)}{\kappa L_2 \big(2 - 3\cos(\eta_1)\big)}
\bigg( \eta_2  - \big(\kappa + \tfrac{m_3 L_2 L_3}{4} \cos(\eta_1) \big) \dot y_2  +
\tfrac{m_3 L_3 }{2} \sin \big( y_2 + (1 - \delta) \eta_1 \big) \dot y_1 \bigg) =: \cF_1(\eta_1,\eta_2,y,\dot y)\\
\dot \eta_2 &= -\frac{2(L_2 + 2L_3)}{\kappa L_2 \big(2-3 \cos(\eta_1)\big)}
\Bigg( (1-\delta) \eta_2 - (1-\delta) \tfrac{L_3 m_3}{2}\sin \big(y_2 + (1-\delta) \eta_1 \big) \dot y_1 - \tfrac{L_2 L_3 m_3}{4}\cos(\eta_1) \dot y_2 \Bigg) \cdot \\
&\qquad \Bigg( \tfrac{L_3 m_3}{2} \cos\big(y_2 + (1-\delta) \eta_1 \big) \dot y_1 +
\tfrac{L_2 L_3 m_3}{4} \sin(\eta_1) 
\tfrac{2(L_2 + 2L_3)}{\kappa L_2 \big(2 - 3\cos(\eta_1 ) \big)}
 \bigg( \kappa \dot y_2 - \delta \eta_2 - \tfrac{\delta m_3 L_3 }{2} \sin\big(y_2  + (1- \delta) \eta_1 \big) \dot y_1 \bigg) \Bigg) \\
 & \qquad - D\frac{2(L_2 + 2L_3)}{\kappa L_2 \big(2 - 3\cos(\eta_1)\big)}
 \bigg( \eta_2  - \big(\kappa + \tfrac{m_3 L_2 L_3}{4} \cos(\eta_1) \big) \dot y_2  +
\tfrac{m_3 L_3 }{2} \sin \big( y_2 + (1 - \delta) \eta_1 \big) \dot y_1 \bigg)
- c \eta_1 := \cF_2(\eta_1,\eta_2,y,\dot y)
\end{aligned}
\end{equation}
\end{figure*}
Denote with $J_{\cF,x}(x^0,y^0) \in \R^{p \times q}$ the Jacobian of a function~$\cF \in \cC^1(\R^q \times \R^k \to \R^p)$ with respect to~$x$ at a point $(x^0,y^0) \in \R^q \times \R^k$. Instead of linearizing~\eqref{eq:ID-robot} around the equilibrium point, in order to increase performance we linearize it around the starting point of the reference trajectory at~$t_0$ and the end point at~$t_f$ given by
\[
y^0 = y_{\rm ref}(t_0),\  y^f  = y_{\rm ref}(t_f),\ \eta^0 = \eta^f = \dot y^0 = \dot y^f = 0 \in \R^2,
\]
and combine both points linearly to obtain $Q, P_1$ and $P_2$ as in Section~\ref{Sec:FunCon} from the function~$\cF = (\cF_1, \cF_2)^\top$ defined in~\eqref{eq:ID-robot} by
\begin{equation*} 
\begin{aligned}
Q &=\frac{1}{2} \big( J_{\cF,\eta}(\eta^0,y^0,\dot y^0) + J_{\cF,\eta}(\eta^f,y^f,\dot y^f) \big), \\
 P_1&= \frac{1}{2} \big( J_{\cF,y}(\eta^0,y^0,\dot y^0) + J_{\cF,y}(\eta^f,y^f,\dot y^f) \big), \\
 P_2 &= \frac{1}{2} \big( J_{\cF,\dot y}(\eta^0,y^0,\dot y^0) + J_{\cF,\dot y}(\eta^f,y^f,\dot y^f) \big).
\end{aligned}
\end{equation*}
With the coefficients
\begin{align*}
C_0 &=  \frac{2 L_{2} + 4 L_{3}}{ \kappa L_{2}  }, \\
C_1 &= -\frac{L_3 m_3(L_2 + 2L_3)}{2 \kappa L_2}\big( \sin(\beta^0) + \sin(\beta^f) \big), \\
C_2 &= \frac{4 \kappa + L_2 L_3 m_3}{10 \kappa L_2}(L_2 + 2 L_3),
\end{align*}
where $\beta^0, \beta^f$ are the values of $\beta$ at $t_0$ and $t_f$, resp., induced by the reference trajectory,
these matrices are given by
\begin{align*}
  Q&= \begin{bmatrix} 0 & -C_0 \\ -c & D C_0  \end{bmatrix},\quad P_1 = \begin{bmatrix} 0&0\\0&0 \end{bmatrix} ,\quad  P_2 =\begin{bmatrix} C_1 & C_2 \\ -DC_1 & -DC_2  \end{bmatrix}
\end{align*}
and we consider the corresponding system in the form as in~\eqref{eq:intdyn-lin-aux}, that is
\begin{equation} \label{eq:ID-robot-linearized}
\begin{pmatrix} \dot \eta_1(t) \\ \dot \eta_2(t) \end{pmatrix} = Q \begin{pmatrix} \eta_1(t) \\ \eta_2(t) \end{pmatrix} + P_2 \begin{pmatrix} \dot y_1(t) \\ \dot y_2(t) \end{pmatrix}.
\end{equation}
We may observe that~$Q$ has a positive and a negative eigenvalue, hence the system is not minimum phase, not even locally.
After the transformation
$\bar \eta(t) = \eta(t) - P_2 y(t) $ we obtain
\begin{equation}
\ddt \bar \eta(t) = Q \bar \eta(t) + P y(t),
\end{equation}
where $P =  Q P_2$.
Now, we calculate that~$Q$ has eigenvalues $\mu_{1,2} = \tfrac{D C_0}{2} \mp \sqrt{\big(\tfrac{D C_0}{2}\big)^2 + C_0 c}$.
Note, that for~$c >$ 0 we have $\mu_1 < 0 < \mu_2$, thus~\eqref{eq:ID-robot-linearized} has a hyperbolic equilibrium, whence the linearized internal dynamics have an unstable part. Next we seek a transformation which diagonalizes~$Q$ and hence separates the stable and the unstable part of the internal dynamics:
\begin{align*}
T = \begin{bmatrix} \tfrac{\mu_2}{c} & \tfrac{\mu_1}{c} \\ 1 & 1 \end{bmatrix} \in \Gl_2 \ \
\text{s.t.} \ \
T^{-1} Q T &= \begin{bmatrix} \mu_1 & 0 \\ 0& \mu_2 \end{bmatrix} =: \begin{bmatrix} \hat Q_1 & 0 \\ 0& \tilde Q \end{bmatrix}
\end{align*}
Using the transformation~$\hat \eta(t) = T^{-1} \bar \eta(t)$ we obtain the linearized internal dynamics
\begin{equation*}
\begin{aligned}
\ddt \begin{pmatrix} {\hat \eta}_1(t) \\ {\hat \eta}_2(t) \end{pmatrix} =
\begin{bmatrix} \mu_1 & 0\\ 0 & \mu_2 \end{bmatrix} \begin{pmatrix} \hat \eta_1(t) \\ \hat \eta_2(t) \end{pmatrix} + T^{-1}P y(t).
\end{aligned}
\end{equation*}
According to~\textbf{(A1)} we set~$\tilde P = [0,1] T^{-1} P$, then this has exactly the form of~\eqref{eq:intdyn-lin-unst} in Section~\ref{Sec:FunCon}:
\begin{equation} \label{eq:ID-robot-linearized-diagonalized}
\ddt \hat \eta_2(t) = \tilde Q \hat \eta_2(t) + \tilde P y(t),
\end{equation}
where~$\tilde Q, \tilde P$ satisfy~\textbf{(A1)} with~$q = 1$, $k_1=1$,
\[
    K_1 = -0.1\quad \text{and}\quad K_2 = [1, 0.01];
\]
however, strictly speaking the last condition in~\textbf{(A1)} is only satisfied in a neighborhood of the equilibrium point.
\textbf{(A2)} is satisfed with~$W = 1$.

Hereinafter~$\hat \eta_2$ plays the role of~$\eta_2$ in Section~\ref{Sec:FunCon} which is not to be confused with the expressions in~\eqref{eq:ID-robot}.
We have
\begin{equation} \label{eq:ynew}
y_{{\rm new},1} (t) := K_1 \hat \eta_2(t),\quad y_{{\rm new},2} (t) = K_2 y(t).
\end{equation}
To replace~$\hat \eta_2$ and its derivatives in~$y_{\rm new,1}$ consider
\begin{small}
\begin{align} \label{eq:Psi-roboter}
\Psi: U & \to \R,  \nonumber \\
x \mapsto
& [0,1] T^{-1}\begin{pmatrix} x_5  \\ \tfrac{m_3 L_3}{2} \sin( x_4 + x_5) x_7 + \kappa (x_9 + x_{10}) + \tfrac{m_3 L_2 L_3}{4} \cos(x_5) x_9 \end{pmatrix} \nonumber \\
 + &[0,1] T^{-1} P_2 \begin{pmatrix} x_2 \\ x_4 + \delta x_5\end{pmatrix}.
\end{align}
\end{small}
Then~$\hat \eta_2(t) = \Psi(q(t),v(t))$ and~$y_{\rm new,1}(t) = K_1 \Psi(q(t),v(t))$. In order to implement the controller~\eqref{eq:Funnel-Control} we calculate the required derivatives of~$y_{\rm new,1}$ in terms of the system's original variables.
Recall $y(t) = h(q(t))$ and $\dot y(t) = h'(q(t))v(t)$. Then, using~\eqref{eq:ID-robot-linearized-diagonalized}:
\begin{equation} \label{eq:ynew-derivatives}
\begin{aligned}
y_{\rm new,1}^{[1]}(t) &\overset{\eqref{eq:ID-robot-linearized-diagonalized}}{=} K_1 \tilde Q \Psi(q(t),v(t)) + K_1 \tilde P  h(q(t)) \\
y_{\rm new,1}^{[2]}(t) &= K_1 \tilde Q^2 \Psi(q(t),v(t)) + K_1 \tilde Q \tilde P h(q(t))\\
&\quad + K_1 \tilde P h'(q(t))v(t).
\end{aligned}
\end{equation}
Further, we need to express~$e_{i,j}(t)$ in~\eqref{eq:Funnel-Control} in terms of the system's original variables. Therefore, we need to calculate~$\hat y_{\rm ref,1}$ according to~\eqref{eq:new-ref}:
\begin{equation} \label{eq:robot-hat-y-ref}
\begin{aligned}
\ddt{\hat{\eta}}_{2,\rm ref}(t) &= \tilde Q \hat \eta_{2,\rm ref}(t) + \tilde P y_{\rm ref}(t), && \hat \eta_{2,\rm ref}(0)= \hat \eta_{2,\rm ref}^0, \\
    \hat y_{\rm ref,1}(t) &= K_1 \hat \eta_{2,\rm ref}(t),
\end{aligned}
\end{equation}
where~$\hat \eta_{2,\rm ref}^0$ is computed accoding to~\eqref{eq:eta2-ref-0}:
\begin{equation} \label{eq:robot-eta-ref}
\hat \eta_{2,\rm ref}^0 = - \int_0^\infty e^{- \mu_2 s} \tilde P \, y_{\rm ref}(s)\, {\rm d}s \ \in \R.
\end{equation}
Now, for $e_{1,0}(t) := y_{\rm new,1}(t) - \hat y_{\rm ref,1}(t)$ we define for $\vp_0 \in \Phi_3$ and~$\kappa_0 > 0$ the expressions
\begin{equation} \label{eq:error-derivative-approx}
\begin{aligned}
e_{1,0}^{[1]}(t) &= y_{\rm new,1}^{[1]}(t) - \ddt \hat y_{\rm ref,1}(t), \\
e_{1,0}^{[2]}(t) &= y_{\rm new,1}^{[2]}(t) -  \hat y_{\rm ref,1}^{(2)}(t), \\
k_{1,0}^{[1]}(t) &= \frac{2 \kappa_0 }{1-\vp_0(t)^2  e_{1,0}(t)^2} \cdot \\
&\quad \cdot \big(\vp_0(t) \dot \vp_0(t) e_{1,0}(t)^2 + \vp_0(t)^2 e_{1,0}(t) e_{1,0}^{[1]}(t) \big), \\
e_{1,1}^{[1]}(t) &= e_{1,0}^{[2]}(t) + k_{1,0}(t) e_{1,0}^{[1]}(t) + k_0^{[1]}(t) e_{1,0}(t),
%
\end{aligned}
\end{equation}
where~$\hat y_{\rm ref,1}$ is from~\eqref{eq:robot-hat-y-ref} and~$y_{\rm new,1}^{[i]}$ are from~\eqref{eq:ynew-derivatives}, resp.
For both errors we choose the same funnel function, i.e., $\vp_{i,j} = \vp_j \in \Phi_{3-j}$ and the same amplification factor, i.e., $\kappa_{i,j} = \kappa_j$ for $i=1,2$ and $j=0,\ldots,2-i$.
With this the controller reads
\ \\
\fbox{\parbox{0.47\textwidth}{%
\begin{equation*}
\begin{aligned}
\ddt{\hat{\eta}}_{2,\rm ref}(t)
&= \tilde Q \hat \eta_{2,\rm ref}(t) + \tilde P y_{\rm ref}(t), \ \hat \eta_{2,\rm ref}(0)= \hat \eta_{2,\rm ref}^0 \, \rm via~\eqref{eq:robot-eta-ref} \\
e_{1,0}(t) &= K_1 \big( \Psi(q(t),v(t)) - \hat \eta_{2,\rm ref}(t)\big),  \\
e_{1,1}(t) &= e_{1,0}^{[1]}(t) + k_{1,0}(t)e_{1,0}(t) \quad \text{using~\eqref{eq:error-derivative-approx}+\eqref{eq:Psi-roboter}},  \\
e_{1,2}(t) &= e_{1,1}^{[1]}(t) + k_{1,1}(t)e_{1,1}(t) \quad \text{using~\eqref{eq:error-derivative-approx}+\eqref{eq:Psi-roboter}},  \\
e_{2,0}(t) &= K_2 \big( y(t) - y_{\rm ref}(t) \big), \\
e_{2,1}(t) &= K_2 \big( \dot y(t) - \dot y_{\rm ref}(t) \big) + k_{2,0}(t) e_{2,0}(t) , \\
k_{i,j}(t)&=\frac{\kappa_j}{1-\varphi_j(t)^2 e_{i,j}(t)^2}, \quad i=1,2, \ j=0,\ldots,2\!-\!i,  \\
\bar e(t) &= \big(e_{1,2}(t), e_{2,1}(t) \big)^\top, \\
\bar k(t) &= \frac{\bar \kappa}{1-\varphi(t)^2\|\bar e(t)\|^2},  \\
u(t) &= - \bar k(t) \, \bar e(t),
\end{aligned}
\end{equation*}
}}
for $\bar \kappa = \kappa_2$ and $\varphi = \vp_2$.
For the simulation we have chosen the funnel functions $\vp_j(t) = (p_j e^{-q_j t} + r_j)^{-1}$ for $j=0,1,2$. The simulation parameters are given in Table~\ref{Tab:Parameters}.
\begin{table}[h!]
\begin{center}
\begin{tabular}{ccccc}
$j$ & $p_j$ & $q_j$ & $r_j$ & $\kappa_j$ \\ \hline
 0 & 0.5 & 2 & 0.001 & 1\\
  1 & 1 & 2 & 0.001  & 1\\
   2 & 1 & 2 & 0.001 & 50
\end{tabular}
\caption{Funnel control design parameters.} 
\label{Tab:Parameters}
\end{center}
\end{table}
Based on the parameters for the reference model from Table~\ref{tab:flexarm-param} the matrices $T^{-1} Q T$, and $T^{-1} P$ are given by
\begin{align*}
T^{-1}Q T &= \begin{bmatrix} -24.8623 & 0 \\ 0 & 28.3917 \end{bmatrix}, \\
T^{-1} P &= \begin{bmatrix} -211.5623 & 235.0693 \\ -305.9704 & 339.9671\end{bmatrix}
\end{align*}
and thus $\tilde Q = 28.3917$, $\tilde P = [-305.9704, 339.9671]$.
\\

The feedforward control~$u_{\rm ff}$ is obtained from solving the boundary value problem described in Section~\ref{Sec:Servo}, based on the inverse model in the DAE formulation~(\ref{eq:servo-constraintsDAE}). {In order to apply the simplified boundary conditions~\eqref{eq:simpbcs_dae} to this example, $n^0+n^f=2n+\ell+p+m=14$ initial or final conditions have to be selected from the equilibrium points at time $T_0$ and $T_f$, while the remaining $14$ entries in~\eqref{eq:simpbcs_dae} remain as free bounds. There are many possible combinations to select these $14$ conditions. Here, in an heuristic way one may choose e.g.}
\begin{align*}
\begin{aligned}
&\begin{bmatrix} q(T_0)^\top & v(T_0)^\top &\lambda(T_0)^\top&\mu(T_0)^\top&u(T_0)^\top \end{bmatrix}^\top = \\
&\qquad \left[ \setlength\arraycolsep{4pt} \begin{array}{cccccccccccccc} 0& 0&\alpha_0&\beta_0&0&0&0&\star&\star&\star &0 &\star&0 &0 \end{array} \right]^\top \\
&\begin{bmatrix} q(T_f)^\top & v(T_f)^\top &\lambda(T_f)^\top&\mu(T_f)^\top&u(T_f)^\top \end{bmatrix}^\top =\\
&\qquad\left[ \setlength\arraycolsep{4pt} \begin{array}{cccccccccccccc}   \star & \star &\star &\star &0&\star &\star &\star &\star &\star & \star & 0 & 0 & 0 \end{array} \right]^\top \,, \label{eqn:twoarm_simp_bcs}
\end{aligned}
\end{align*}
where free bounds are denoted by~$\star$. The initial angles~$\alpha_0$ and $\beta_0$ are given by the geometry of the model as~$\alpha_0=\arccos\left(\frac{L_1^2+d^2-0.25L_2^2}{2L_1\,d}\right)$ and $\beta_0=\arcsin\left(\frac{2\,L_1\sin(\alpha_0)}{L_2}\right)$. The BVP solution includes a pre-actuation phase before the start of the trajectory at time $t_0=\SI{0}{\second}$. However, here only the causal part for $t\geq t_0$ is considered, such that
\[
u_{\rm ff}(t) =
\begin{cases}
0, & t<t_0, \\
u_{\rm bvp}(t), & t\geq t_0,
\end{cases}
\]
where $u_{\rm bvp}$ denotes the solution of the BVP. Note that the clipping of the input signal will also yield small tracking errors, which have to be reduced by the feedback controller. Furthermore, the solution of the BVP does not only yield the feedforward control~$u_{\rm ff}$, but also the desired state trajectories~$q_{\rm ref},v_{\rm ref},\lambda_{\rm ref}$, which may be used for reference in a feedback control law.

In the simulations, three controller configurations are compared. The controller~C$_1$ is the combination $u(t)=u_{\rm fb}(t) + u_{\rm ff}(t)$ of the feedback and feedforward strategy discussed above. The controller~C$_2$ is pure feedback control $u(t)=u_{\rm fb}(t)$ and the controller~C$_3$ is pure feedforward control $u(t)=u_{\rm ff}(t)$. The simulations over the interval~$0-2\,\rm s$ are performed on \textsc{Matlab} with the solver ode15s (AbsTol and RelTol at default values).

Snapshots of an animation of the simulation with controller~C$_1$ are shown in Figure~\ref{fig:animation_snaps} for different time instances. The motion of the complete robot is visualized and the motion of the actuators attached to the moving bases with~$s_1(t)$ and $s_2(t)$ can be seen. Moreover, the deflection $\gamma(t)$ of the third, passive arm is visible at time~$t=0.4\;\text{s}$.

\begin{figure}[h]
\begin{center}
\includegraphics[width=\linewidth]{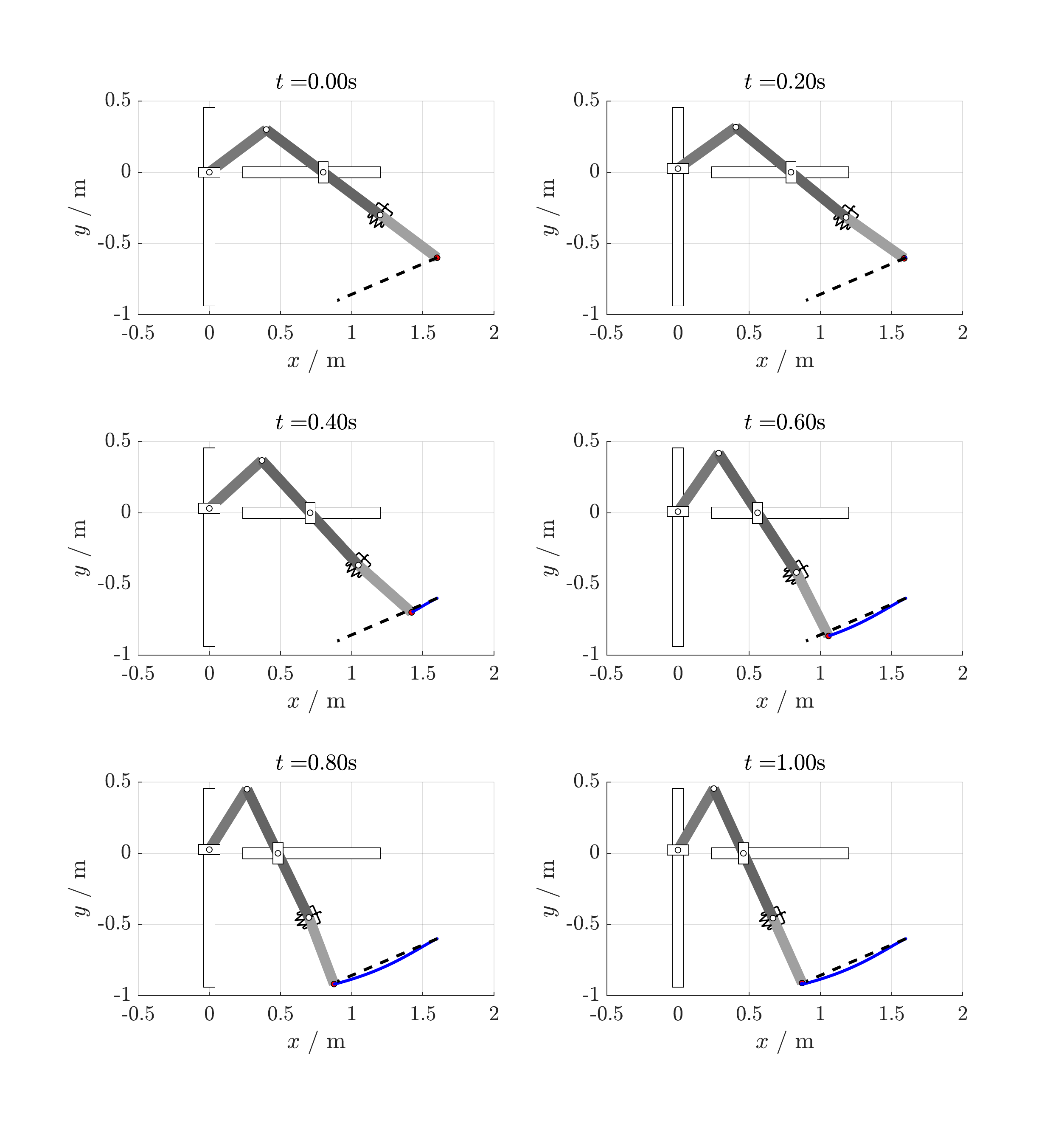}
\caption{Snapshots of an animation of the simulation under controller C$_1$.}
\label{fig:animation_snaps}
\end{center}
\end{figure}

Figure~\ref{fig:States} shows the associated states $s_1(t)$, $s_2(t)$, $\alpha(t)$, $\beta(t)$ and~$\gamma(t)$ during the motion, applying the combined controller~$\rm C_1$. The oscillation of the third, passive body due to the motion is visible as an oscillation of~$\gamma(t)$. This motion cannot be actuated directly and is therefore responsible for the involved controller design.
\begin{figure}[h]
\includegraphics[width=0.95\linewidth]{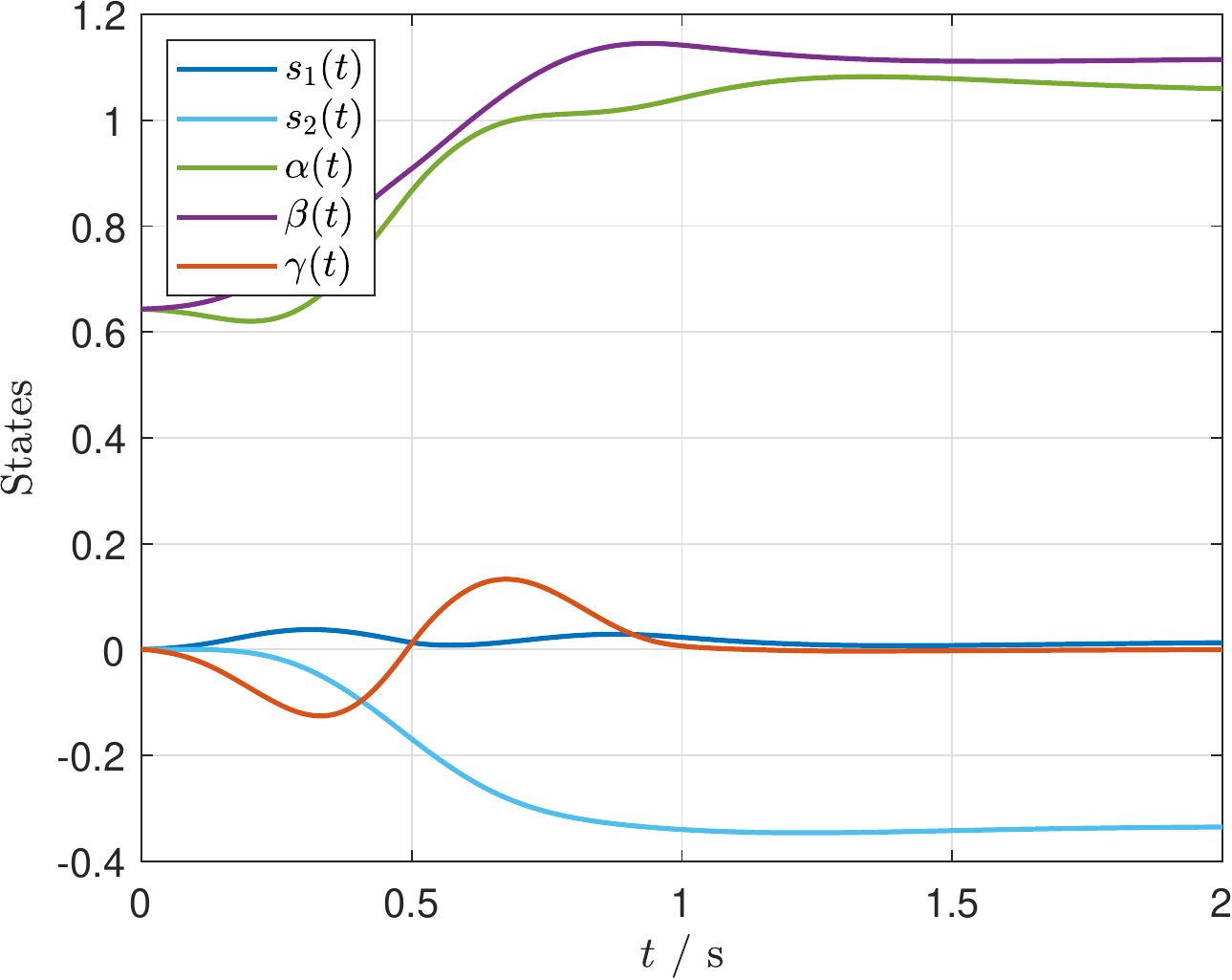}
\caption{State trajectories using the combined controller C$_1$ ($s_1$, $s_2$ in m and $\alpha$, $\beta$, $\gamma$ in rad).}
\label{fig:States}
\end{figure}

Figure~\ref{fig:Output_in_x_y} shows the end-effector position in the two-dimensional space under all three controller configurations C$_1$--C$_3$. While at the beginning a good tracking is achieved, the pure feedforward controller C$_3$ does not reach the desired final position at $(e_1,e_2)=(0.9,-0.9)\;\text{m}$. While both controllers C$_1$ and C$_2$ reach the correct desired final position, the controller C$_1$ shows smaller tracking errors over the course of the trajectory. The tracking errors are evaluated in the following.
\begin{figure}[h]
\includegraphics[width=1\linewidth]{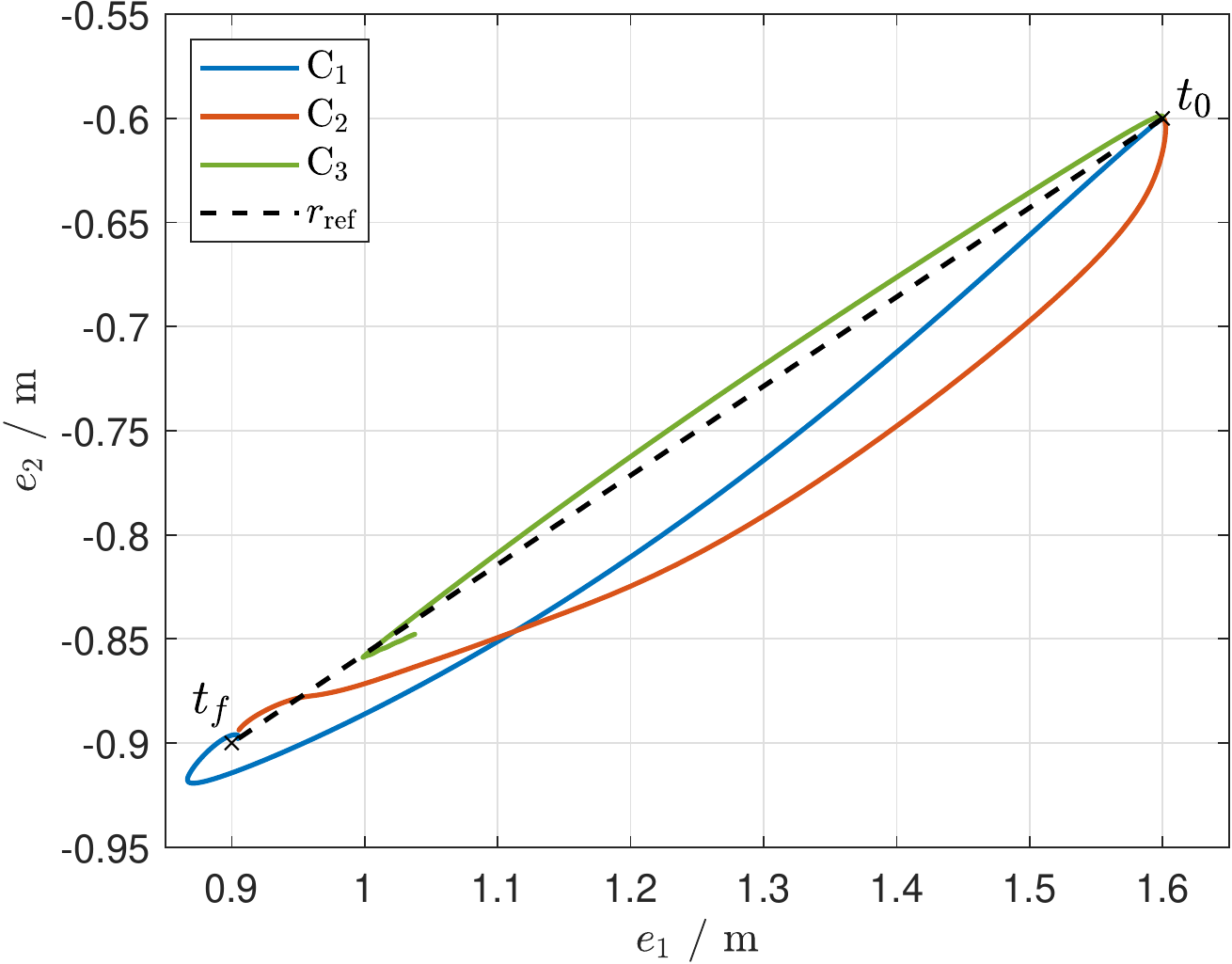}
\caption{End-effector position and reference.}
\label{fig:Output_in_x_y}
\end{figure}

Figure~\ref{fig:y-yref} shows the errors between the original output~$y$ and the reference signal~$y_{\rm ref}$, and Figure~\ref{fig_r-rapp} the errors in end-effector coordinates $\| r_{\rm app} - r_{\rm ref} \|$ under the controllers $\rm C_1, C_2$ and C$_3$, resp. The pure feedforward controller C$_3$ cannot compensate the increasing tracking errors due to the mass~$m_3$, which is larger in the simulated model compared to the reference model used for controller design. In contrast, both controllers C$_1$ and C$_2$, which include a feedback channel, reduce the tracking error after an initial growth period. Applying the combined controller~$\rm C_1$ results in the smallest tracking error at the end time~$t_f = 2\,\rm s$. Moreover, the cumulative error of the controller C$_1$ is considerably smaller compared to the controller C$_2$. The controller C$_1$ reduces the cumulative error of C$_2$ by 40\% and 50\%, resp., for the coordinates of~$y$ and~$r_{\rm app}$.

\begin{figure}[h]
\begin{center}
\subfigure[Errors between the original output $y(t)$ and reference signal $y_{\rm ref}(t)$. \label{fig:y-yref}]{\includegraphics[width=0.95\linewidth]{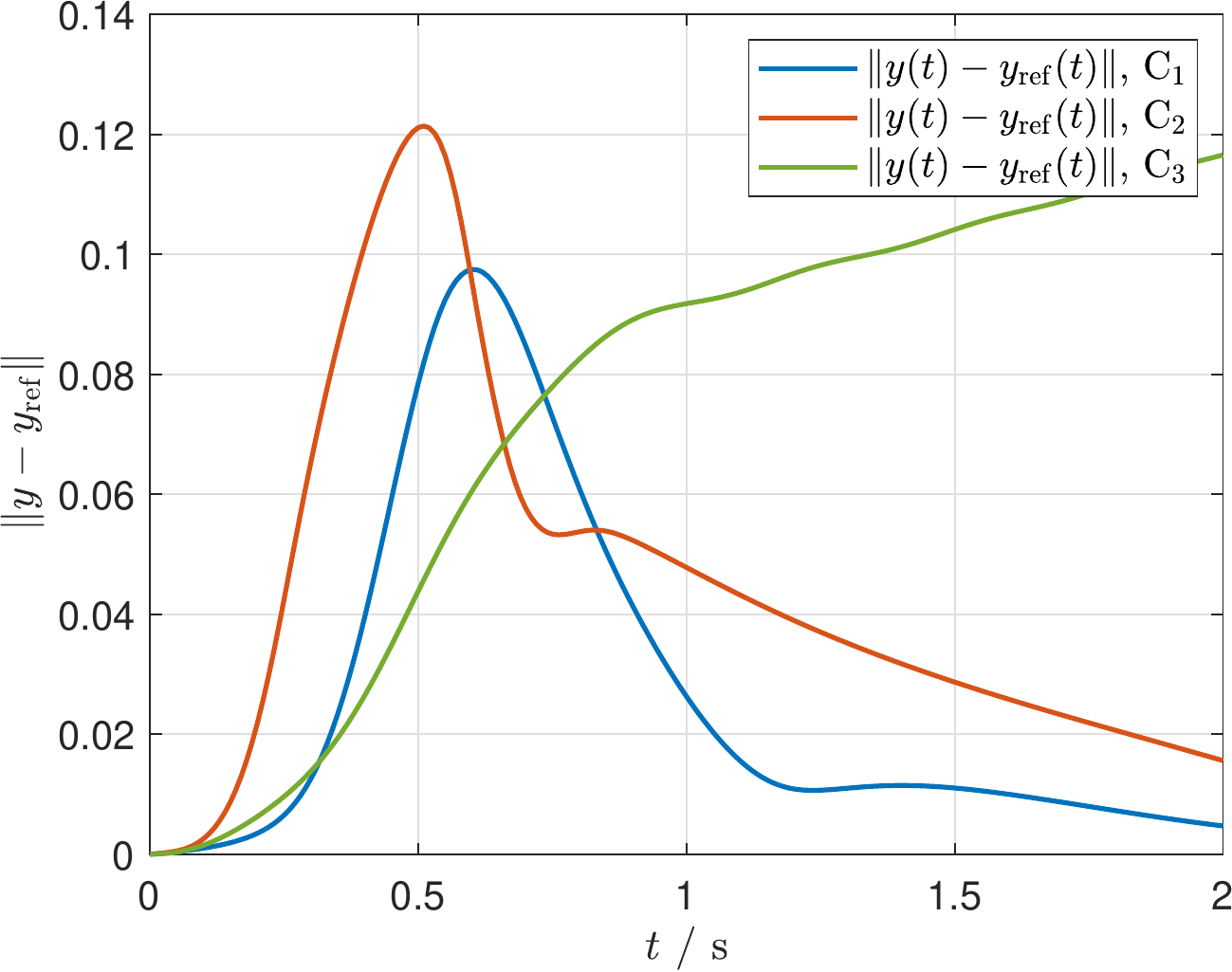}}
\subfigure[Errors between the approximated end-effector $r_{\rm  app}(t)$ and the refence $r_{\rm ref}(t)$ in end-effector coordinates.
\label{fig_r-rapp}]{\includegraphics[width=0.95\linewidth]{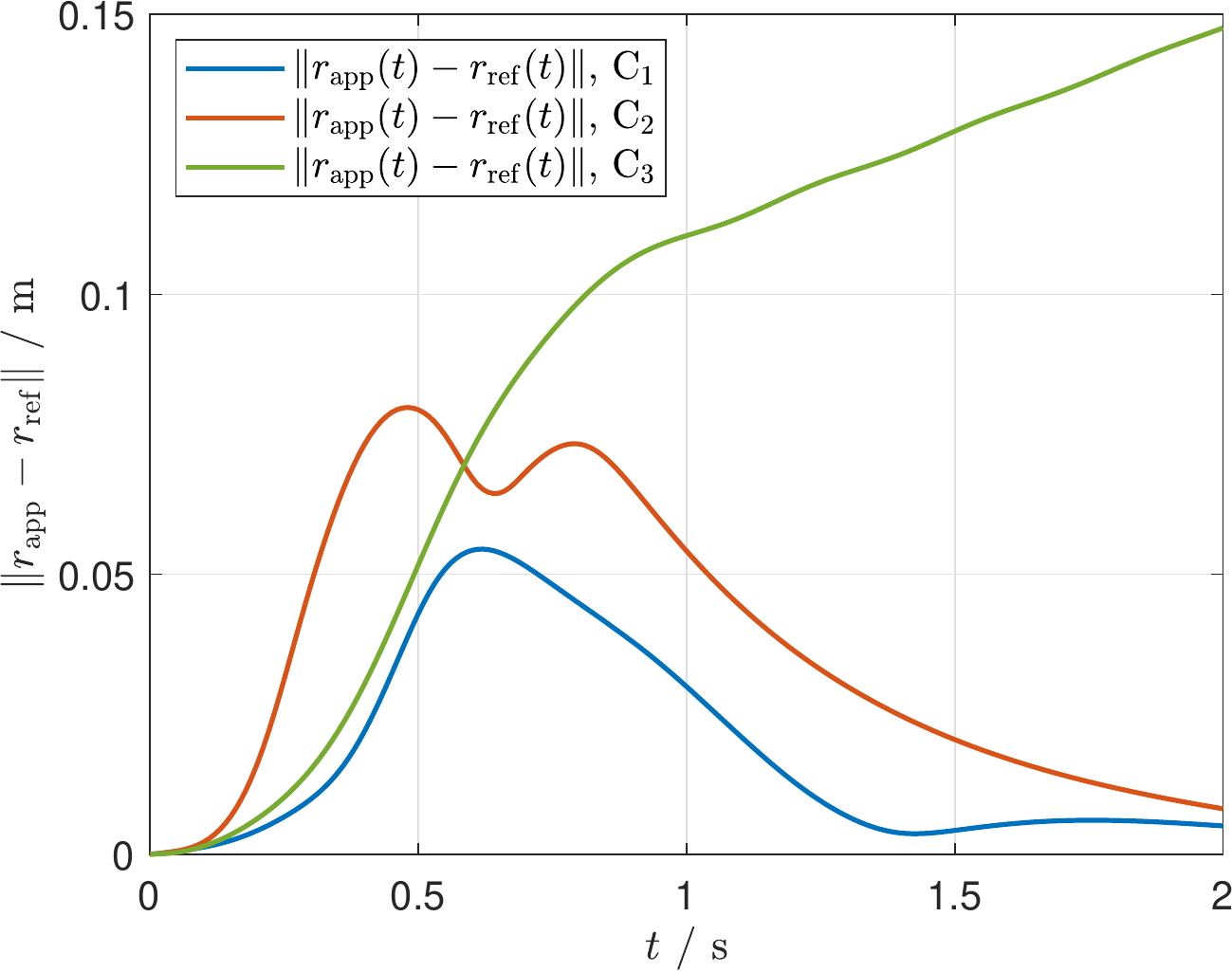}}
\caption{Output errors and end-effector errors.}
\label{fig:Output_and_endeffector_errors}
\end{center}
\end{figure}

In Figure~\ref{fig:Inputs} the input functions~$u_1$ and $u_2$ are depicted. The controller C$_3$ shows the smooth feedforward control signal from model inversion of the reference model. The funnel controller C$_2$ varies strongly, especially in the first component~$u_1$. These large derivatives result in large accelerations of the physical system and should be avoided in order to reduce the loads on the mechanical parts. This can be achieved by the controller C$_1$. Since it includes the smooth feedforward signal, which roughly moves the system on the desired path, the work load of the funnel feedback controller is reduced and the peaks and strong variations of the input signal can be completely avoided here. The input signals of C$_1$ are smooth and therefore reduce the loads on the mechanical parts of the system.

\begin{figure}[h]
\subfigure[First input components $u_1$.]{\includegraphics[width=0.95\linewidth]{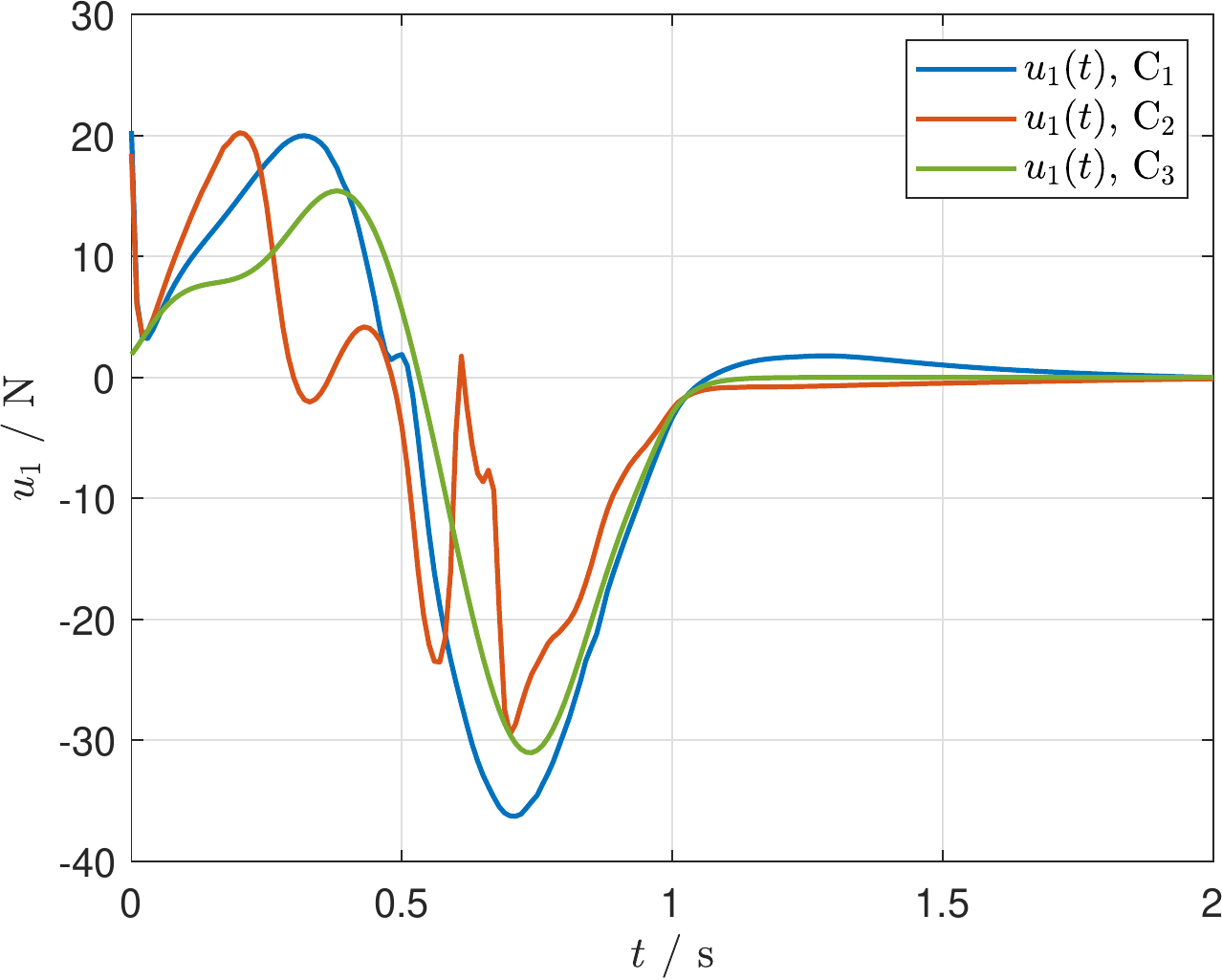}} 
\subfigure[Second input components $u_2$.]{\includegraphics[width=0.95\linewidth]{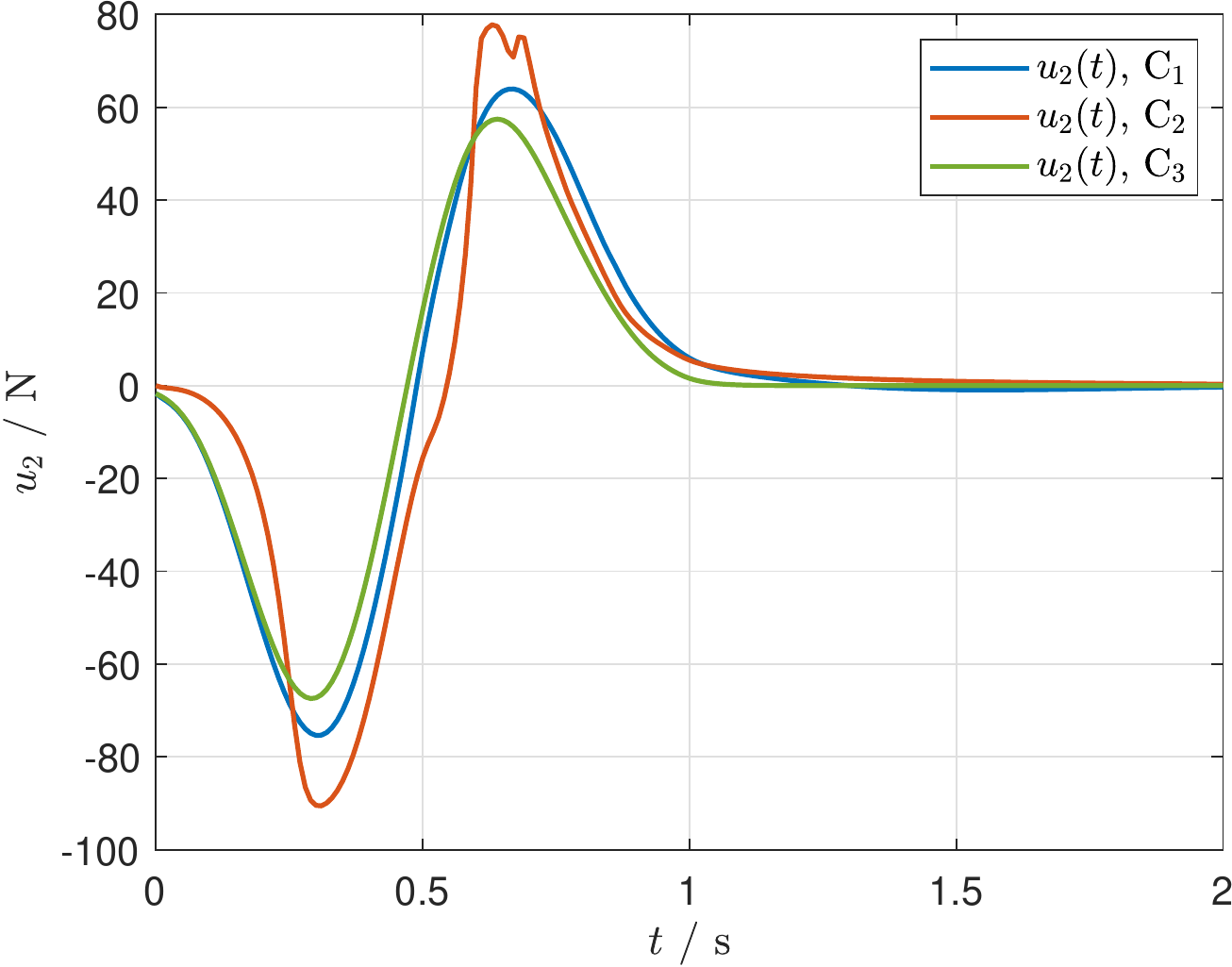}} 
\caption{Components of the control input.}
\label{fig:Inputs}
\end{figure}

In Figure~\ref{fig:Funnel_and_error} the error norm~$\|\bar e(\cdot)\|$ and the funnel boundary $\vp(\cdot)^{-1}$  under the controllers~$\rm C_1$ and~$\rm C_2$ are depicted. In both controller configurations the error lies within the funnel boundary.
However, in the scenario~$\rm C_2$ when only the feedback controller is applied, the error is closer to the funnel boundary which results in a larger input, cf. Figure~\ref{fig:Inputs}.

\begin{figure}[h]
\includegraphics[width=0.95\linewidth]{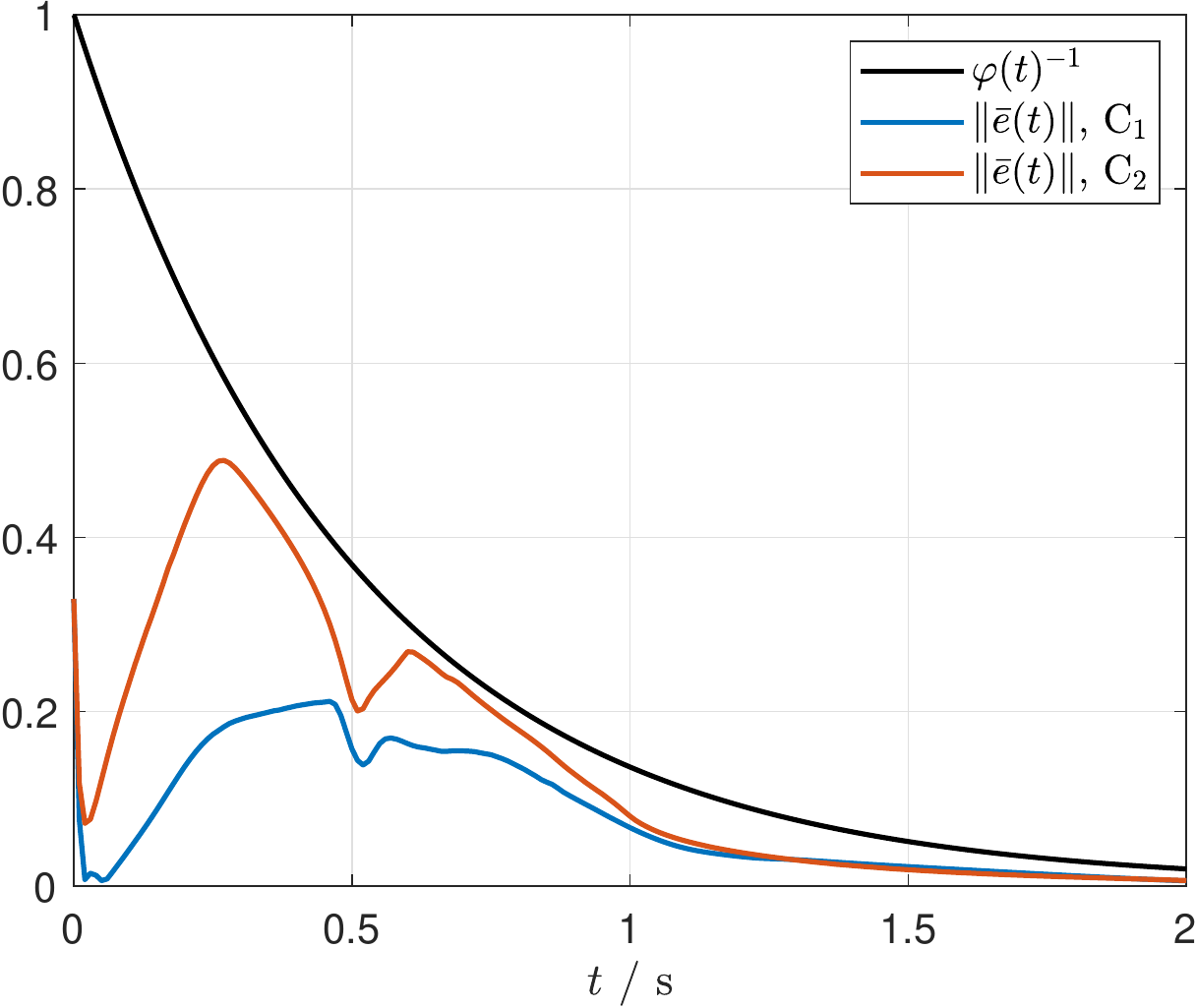}
\caption{Funnel boundary and error norm.}
\label{fig:Funnel_and_error}
\end{figure}

\section{Summary and Conclusion}

In the present work we successfully combined a feedforward control strategy based on servo-constraints with a novel high-gain feedback controller design in order to achieve output tracking for a multibody system with kinematic loop and unstable internal dynamics. First, we introduced auxiliary inputs and outputs in order to avoid the DAE formulation and decouple the internal dynamics. Then we derived a feasible set of coordinates for the internal dynamics which allows a decoupling in terms of the original system parameters and avoids an explicit calculation of the Byrnes-Isidori form.
In order to apply the recent result from~\cite{Berg20a} on funnel control of linear non-minimum phase systems, we considered a linearization of the internal dynamics of the multibody system. Furthermore, a system inversion based on servo-constraints is applied as a feedforward control. Due to the unstable internal dynamics, this includes solving a boundary value problem beforehand.

As a proof of concept we applied the combination of both control strategies to a nonlinear multi-input, multi-output multibody system with unstable internal dynamics in a simulation and compared the results. We found the combination of both controllers to have the best tracking performance in terms of the error between the system output and the reference signal. Moreover, the combination reduces peaks and strong variations in the input signal compared to the pure feedback strategy. The results motivate further research on the combination of open-loop and closed-loop control strategies for non-minimum phase systems.

\section*{Compliance with ethical standards}
Disclosure of potential conflicts of interest: The authors de-clare that they have no conflict of interest.

\bibliographystyle{spmpsci}

\begin{thebibliography}{10}
\providecommand{\url}[1]{{#1}}
\providecommand{\urlprefix}{URL }
\expandafter\ifx\csname urlstyle\endcsname\relax
  \providecommand{\doi}[1]{DOI~\discretionary{}{}{}#1}\else
  \providecommand{\doi}{DOI~\discretionary{}{}{}\begingroup
  \urlstyle{rm}\Url}\fi

\bibitem{AltmBets16}
Altmann, R., Betsch, P., Yang, Y.: Index reduction by minimal extension for the
  inverse dynamics simulation of cranes.
\newblock Multibody Sys. Dyn. \textbf{36}(3), 295--321 (2016)

\bibitem{AltmHeil17}
Altmann, R., Heiland, J.: Simulation of multibody systems with servo
  constraints through optimal control.
\newblock Multibody Sys. Dyn. \textbf{40}(1), 75--98 (2017)

\bibitem{AschMatt95}
Ascher, U., Mattheij, R., Russell, R.: Numerical Solution of Boundary Value
  Problems for Ordinary Differential Equations.
\newblock Society for Industrial and Applied Mathematics (1995).

\bibitem{BastSeif13}
Bastos, G., Seifried, R., Br{\"u}ls, O.: Inverse dynamics of serial and
  parallel underactuated multibody systems using a {DAE} optimal control
  approach.
\newblock Multibody System Dynamics \textbf{30}(3), 359--376 (2013).

\bibitem{BechRovi14}
Bechlioulis, C.P., Rovithakis, G.A.:  A low-complexity global
approximation-free control scheme with prescribed performance for unknown pure feedback systems.
\newblock Automatica \textbf{50}, 1217--1226 (2014)

\bibitem{BencLasz17}
Bencsik, L., Kov{\'a}cs, L.L., Zelei, A.: Stabilization of internal dynamics of
  underactuated systems by periodic servo-constraints.
\newblock International Journal of Structural Stability and Dynamics
  \textbf{17}(05), 1740004 (2017)

\bibitem{Berg14a}
Berger, T.: On differential-algebraic control systems.
\newblock Ph.D. thesis, Institut f{\"u}r Mathematik, Technische Universit{\"a}t
  Ilmenau, Universit{\"a}tsverlag Ilmenau, Germany (2014)

\bibitem{Berg17b}
Berger, T.: The zero dynamics form for nonlinear differential-algebraic
  systems.
\newblock {IEEE} Trans. Autom. Control \textbf{62}(8), 4131--4137 (2017)

\bibitem{Berg20a}
Berger, T.: Tracking with prescribed performance for linear non-minimum phase
  systems.
\newblock Automatica \textbf{115}, Article 108909 (2020)

\bibitem{BergHall20pp}
Berger, T., Haller, F.E.: On an extension of a global implicit function theorem
  (2020).
\newblock Submitted for publication, preprint available at arXiv:
  \url{https://arxiv.org/abs/2004.04427}

\bibitem{BergLanz20}
Berger, T., Lanza, L.: Output tracking for a non-minimum phase robotic
  manipulator (2020).
\newblock To appear in Proc. MTNS 2020. Preprint available at arXiv:
  \url{https://arxiv.org/abs/2001.07535}

\bibitem{BergLe18a}
Berger, T., L{\^e}, H.H., Reis, T.: Funnel control for nonlinear systems with
  known strict relative degree.
\newblock Automatica \textbf{87}, 345--357 (2018).

\bibitem{BergLe20}
Berger, T., L\^{e}, H.H., Reis, T.: Vector relative degree and funnel control
  for differential-algebraic systems.
\newblock In: S.~Grundel, T.~Reis, S.~Sch\"ops (eds.) {Progress in
  Differential-Algebraic Equations II}, Differential-Algebraic Equations Forum.
  Springer-Verlag, Berlin-Heidelberg (2020).
\newblock To appear. Preprint available at arXiv:
  \url{https://arxiv.org/abs/2001.05391}

\bibitem{BergOtto19}
Berger, T., Otto, S., Reis, T., Seifried, R.: Combined open-loop and funnel
  control for underactuated multibody systems.
\newblock Nonlinear Dynamics \textbf{95}, 1977--1998 (2019).

\bibitem{BergRaue18}
Berger, T., Rauert, A.L.: A universal model-free and safe adaptive cruise
  control mechanism.
\newblock In: Proceedings of the MTNS 2018, pp. 925--932. Hong Kong (2018)

\bibitem{BergRaue20}
Berger, T., Rauert, A.L.: Funnel cruise control.
\newblock Automatica \textbf{119}, Article 109061 (2020)

\bibitem{BergReis14a}
Berger, T., Reis, T.: Zero dynamics and funnel control for linear electrical
  circuits.
\newblock J. Franklin Inst. \textbf{351}(11), 5099--5132 (2014)

\bibitem{BetsAltm16}
Betsch, P., Altmann, R., Yang, Y.: Numerical integration of underactuated
  mechanical systems subjected to mixed holonomic and servo constraints.
\newblock In: J.M. Font-Llagunes (ed.) Multibody Dynamics, \emph{Computational
  Methods in Applied Sciences}, vol.~42, pp. 1--18. Springer-Verlag, Cham
  (2016).

\bibitem{BetsQuas09}
Betsch, P., Quasem, M., Uhlar, S.: Numerical integration of discrete mechanical
  systems with mixed holonomic and control constraints.
\newblock J. Mech. Sci. Tech. \textbf{23}(4), 1012--1018 (2009)

\bibitem{Blaj92}
Blajer, W.: Index of differential-algebraic equations governing the dynamics of
  constrained mechanical systems.
\newblock Appl. Math. Mod. \textbf{16}(2), 70--77 (1992)

\bibitem{BlajKolo04}
Blajer, W., Ko{\l}odziejczyk, K.: A geometric approach to solving problems of
  control constraints: {T}heory and a {DAE} framework.
\newblock Multibody Sys. Dyn. \textbf{11}(4), 343--364 (2004)

\bibitem{BlajKolo11}
Blajer, W., Ko{\l}odziejczyk, K.: Improved {DAE} formulation for inverse
  dynamics simulation of cranes.
\newblock Multibody Sys. Dyn. \textbf{25}(2), 131--143 (2011)

\bibitem{BrenCamp89}
Brenan, K.E., Campbell, S.L., Petzold, L.R.: Numerical Solution of
  Initial-Value Problems in Differential-Algebraic Equations.
\newblock North-Holland, Amsterdam (1989)

\bibitem{BrueBast13}
Br\"uls, O., Bastos, G.J., Seifried, R.: A stable inversion method for
  feedforward control of constrained flexible multibody systems.
\newblock Journal of Computational and Nonlinear Dynamics \textbf{9}(1),
  011,014--011,014--9 (2013).

\bibitem{BurkSeifEber19}
Burkhardt, M., Seifried, R., Eberhard, P.: Experimental studies of control
  concepts for a parallel manipulator with flexible links.
\newblock Journal of Mechanical Science and Technology \textbf{29}(7),
  2685--2691 (2015).


\bibitem{Camp95b}
Campbell, S.L.: High-index differential algebraic equations.
\newblock Mechanics of Structures and Machines \textbf{23}(2), 199--222 (1995)

\bibitem{ChenPade96}
Chen, D., Paden, B.: Stable inversion of nonlinear non-minimum phase systems.
\newblock Int. J. Control \textbf{64}(1), 81--97 (1996)

\bibitem{DaiHe18}
Dai, S.-L., He, S., Lin, H., Wang, C.: Platoon formation control with prescribed performance guarantees for USVs.
\newblock {IEEE} Trans. Ind. Electronics \textbf{65}(5), 4237--4246 (2018)

\bibitem{Deva99}
Devasia, S.: Approximated stable inversion for nonlinear systems with
  nonhyperbolic internal dynamics.
\newblock {IEEE} Trans. Autom. Control \textbf{44}(7), 1419--1425 (1999)

\bibitem{DevaChen96}
Devasia, S., Chen, D., Paden, B.: Nonlinear inversion-based output tracking.
\newblock {IEEE} Trans. Autom. Control \textbf{41}(7), 930--942 (1996)

\bibitem{DrueSeif20}
Dr\"ucker, S., Seifried, R.: Stable inversion for flexible multibody systems
  using the {ANCF}.
\newblock submitted to Proceedings in Applied Mathematics and Mechanics 2020
  (2020)

\bibitem{FlieLevi95}
Fliess, M., Levine, J., Martin, P., Rouchon, P.: Flatness and defect of
  non-linear-systems: introductory theory and examples.
\newblock Int. J. Control \textbf{61}, 1327--1361 (1995)

\bibitem{Fran77}
Francis, B.A.: The linear multivariable regulator problem.
\newblock {SIAM} J. Control Optim. \textbf{15}, 486--505 (1977)

\bibitem{FumaMasa10}
Fumagalli, A., Masarati, P., Morandini, M., Mantegazza, P.: Control constraint
  realization for multibody systems.
\newblock J. Comput. Nonlinear Dynam. \textbf{6}(1), 011,002--011,002--8 (2010)

\bibitem{Gear88}
Gear, C.W.: Differential-algebraic equation index transformations.
\newblock {SIAM} J. Sci. Stat. Comput. \textbf{9}, 39--47 (1988)

\bibitem{Hack17}
Hackl, C.M.: Non-identifier Based Adaptive Control in Mechatronics--Theory and
  Application, \emph{Lecture Notes in Control and Information Sciences}, vol.
  466.
\newblock Springer-Verlag, Cham, Switzerland (2017)

\bibitem{HairWann96}
Hairer, E., Wanner, G.: Solving ordinary differential equations {II}: {S}tiff
  and differential-algebraic problems, \emph{Springer Series in Computational
  Mathematics}, vol.~14, 2nd edn.
\newblock Springer-Verlag, Berlin (1996)

\bibitem{IlchReis13b}
Ilchmann, A., Reis, T. (eds.): Surveys in {D}ifferential-{A}lgebraic
  {E}quations~{I}--{IV}.
\newblock Differential-Algebraic Equations Forum. Springer-Verlag, Berlin
  (2013--2019)

\bibitem{IlchRyan02b}
Ilchmann, A., Ryan, E.P., Sangwin, C.J.: Tracking with prescribed transient
  behaviour.
\newblock ESAIM: Control, Optimisation and Calculus of Variations \textbf{7},
  471--493 (2002)

\bibitem{IlchTren04}
Ilchmann, A., Trenn, S.: Input constrained funnel control with applications to
  chemical reactor models.
\newblock Syst. Control Lett. \textbf{53}(5), 361--375 (2004)

\bibitem{Isid95}
Isidori, A.: Nonlinear Control Systems, 3rd edn.
\newblock Communications and Control Engineering Series. Springer-Verlag,
  Berlin (1995)

\bibitem{Jin19}
Jin, X.: Adaptive fixed-time control for MIMO nonlinear systems with asymmetric output constraints using universal barrier functions.
\newblock {IEEE} Trans. Autom. Control \textbf{64}(7), 3046--3053 (2019)

\bibitem{KunkMehr06}
Kunkel, P., Mehrmann, V.: Differential-Algebraic Equations. Analysis and
  Numerical Solution.
\newblock EMS Publishing House, Z{\"u}rich, Switzerland (2006).

\bibitem{MorlMeye18}
Morlock, M., Meyer, N., Pick, M.A., Seifried, R.: Modeling and trajectory
  tracking control of a new parallel flexible link robot.
\newblock IEEE International Conference on Intelligent Robots and Systems pp.
  6484--6489 (2018)

\bibitem{MorlMeye21}
Morlock, M., Meyer, N., Pick, M.A., Seifried, R.: Real-time trajectory tracking
  control of a parallel robot with flexible links.
\newblock Mechanism and Machine Theory \textbf{158}, Article 104220 (2021)

\bibitem{OttoSeif18a}
Otto, S., Seifried, R.: Open-loop control of underactuated mechanical systems
  using servo-constraints: Analysis and some examples.
\newblock In: S.~Campbell, A.~Ilchmann, V.~Mehrmann, T.~Reis (eds.)
  Applications of Differential-Algebraic Equations: Examples and Benchmarks,
  Differential-Algebraic Equations Forum, pp. 81--122. Springer-Verlag, Cham
  (2018).

\bibitem{OttoSeif18}
Otto, S., Seifried, R.: Real-time trajectory control of an overhead crane using
  servo-constraints.
\newblock Multibody Sys. Dyn. \textbf{42}(1), 1--17 (2018).

\bibitem{PompWeye15}
Pomprapa, A., Weyer, S., Leonhardt, S., Walter, M., Misgeld, B.: Periodic
  funnel-based control for peak inspiratory pressure.
\newblock In: Proc. 54th~{IEEE} Conf. Decis. Control, Osaka, Japan, pp.
  5617--5622 (2015)

\bibitem{Seif14}
Seifried, R.: Dynamics of Underactuated Multibody Systems. Modeling, Control
  and Optimal Design.
\newblock No. 205 in Solid Mechanics and Its Applications. Springer-Verlag
  (2014)

\bibitem{SenfPaug14}
Senfelds, A., Paugurs, A.: Electrical drive {DC} link power flow control with
  adaptive approach.
\newblock In: Proc. 55th Int. Sci. Conf. Power Electr. Engg. Riga Techn. Univ.,
  Riga, Latvia, pp. 30--33 (2014)

\bibitem{ShamGlad03}
Shampine, L.F., Gladwell, I., Thompson, S.: Solving ODEs with MATLAB.
\newblock Cambridge University Press, Cambridge (2003).
\newblock \urlprefix\url{http://www.runet.edu/~thompson/webddes/}

\bibitem{Sime13}
Simeon, B.: Computational Flexible Multibody Dynamics.
\newblock Differential-Algebraic Equations Forum. Springer-Verlag,
  Heidelberg-Berlin (2013).

\bibitem{SkogPost05}
Skogestad, S., Postlethwaite, I.: Multivariable Feedback Control: Analysis and
  Design, 2nd edn.
\newblock John Wiley and Sons Inc., Chichester (2005)

\bibitem{Sont98a}
Sontag, E.D.: Mathematical Control Theory: Deterministic Finite Dimensional
  Systems, 2nd edn.
\newblock Springer-Verlag, New York (1998)

\bibitem{TeeRen11}
Tee, K.P., Ren, B., Ge, S.S.: Control of nonlinear systems with time-varying output constraints.
\newblock Automatica \textbf{47}, 2511--2516 (2011)

\end{thebibliography}

\end{document}